\documentclass{article}
\usepackage{graphicx}
\usepackage{amsmath}

\makeatletter
\renewcommand\paragraph{\@startsection{paragraph}{4}{\z@}%
                                    {1ex \@plus1ex \@minus.2ex}%
                                    {-1em}%
                                    {\normalfont\normalsize\bfseries}}
\makeatother

{\begin{list}               
    {$\bullet$ \hfill}{
        \setlength{\leftmargin}{\parindent}
        \setlength{\parsep}{0.04\baselineskip}
        \setlength{\itemsep}{0.5\parsep}
        \setlength{\labelwidth}{\leftmargin}
        \setlength{\labelsep}{0em}}
    }
{\end{list}}

\usepackage{amsmath}
\usepackage{amssymb}
\usepackage{amsthm}
\usepackage[margin=1in]{geometry}

\usepackage[nocompress]{cite}

\makeatletter
\renewcommand{\@biblabel}[1]{[#1]\hfill}
\makeatother
\makeatletter
\providecommand{\bstctlcite}{\@ifnextchar[{\@bstctlcite}{\@bstctlcite[@auxout]}}
\providecommand{\@bstctlcite}[2][@auxout]{\@bsphack
  \@for\@citeb:=#2\do{%
    \edef\@citeb{\expandafter\@firstofone\@citeb}%
    \if@filesw\immediate\write\csname #1\endcsname{\string\citation{\@citeb}}\fi}%
  \@esphack}
\makeatother
\AtBeginDocument{\bstctlcite{IEEErefcontrol}}
\newcommand{\bibdoi}[1]{%
  {\urlstyle{same}\href{https://doi.org/#1}{\nolinkurl{#1}}}%
}

\usepackage{color}
\usepackage{enumitem}
\usepackage{hyperref}
 \hypersetup{
     colorlinks=true,
     linkcolor=blue,
     filecolor=blue,
     citecolor = blue,      
     urlcolor=black,
     }
\usepackage{bm}
\usepackage{mathrsfs}
\usepackage[scr=boondox]{mathalpha}

\usepackage{extarrows}
\usepackage{graphicx}
\usepackage{subfig}
\usepackage{mathrsfs}
\usepackage{cmll}
\usepackage{mathtools}

\usepackage[cmyk]{xcolor} 
\usepackage[titletoc,title]{appendix}
\numberwithin{equation}{section}

\newtheorem{theorem}{Theorem}
\newtheorem{proposition}{Proposition}
\newtheorem{lemma}{Lemma}

\theoremstyle{definition}

\newtheorem{definition}{Definition}

\theoremstyle{remark}
\newtheorem{remark}{Remark}

\providecommand{\calR}{\mathcal{R}}

\renewcommand{\vec}[1]{\ensuremath{\boldsymbol{#1}}}
\providecommand{\mat}[1]{\ensuremath{\boldsymbol{#1}}}
\providecommand{\mA}{\mat{A}} \providecommand{\mB}{\mat{B}}
\providecommand{\mC}{\mat{C}} 
\providecommand{\mD}{\mat{D}}
\providecommand{\mF}{\mat{F}}
\providecommand{\mH}{\mat{H}}
\providecommand{\mI}{\mat{I}} 
 
\providecommand{\mK}{\mat{K}} 
 
\providecommand{\mM}{\mat{M}} 
\providecommand{\mP}{\mat{P}} 
\providecommand{\mQ}{\mat{Q}} 
\providecommand{\mR}{\mat{R}}
\providecommand{\mS}{\mat{S}} 
\providecommand{\mU}{\mat{U}} 
\providecommand{\mV}{\mat{V}}
\providecommand{\mW}{\mat{W}}
\providecommand{\mT}{\mat{T}}
\providecommand{\mZ}{\mat{Z}}
\providecommand{\mSigma}{\mat{\Sigma}}
 
\providecommand{\mG}{\mat{G}}

\providecommand{\mX}{\mat{X}}

\providecommand{\mY}{\mat{Y}}
\newcommand{\mLambda}{\mat{\Lambda}}
\providecommand{\mPhi}{\mat{\Phi}}
\providecommand{\va}{\vec{a}}

\providecommand{\vg}{\vec{g}}

\providecommand{\vq}{\vec{q}}

\providecommand{\vt}{\vec{t}}
\providecommand{\vu}{\vec{u}} 
\providecommand{\vw}{\vec{w}}
\providecommand{\vx}{\vec{x}} 
\providecommand{\vy}{\vec{y}}
\providecommand{\vz}{\vec{z}}

\providecommand{\vv}{\vec{v}}

\providecommand{\vone}{\vec{1}}

\providecommand{\vtheta}{\vec{\theta}}
\providecommand{\vphi}{\vec{\phi}}
\providecommand{\veta}{\boldsymbol{\eta}}

\DeclareMathOperator{\Tr}{Tr}

\renewcommand{\tilde}{\widetilde}
\renewcommand{\hat}{\widehat}
\DeclareFontFamily{U}{mathx}{\hyphenchar\font45}
\DeclareFontShape{U}{mathx}{m}{n}{
      <5> <6> <7> <8> <9> <10>
      <10.95> <12> <14.4> <17.28> <20.74> <24.88>
      mathx10
      }{}
\DeclareSymbolFont{mathx}{U}{mathx}{m}{n}
\DeclareFontSubstitution{U}{mathx}{m}{n}
\DeclareMathAccent{\widecheck}{0}{mathx}{"71}
\DeclareMathAccent{\wideparen}{0}{mathx}{"75}

\newcommand{\BE}{\begin{equation}}
\newcommand{\EE}{\end{equation}}
\newcommand{\BS}{\begin{subequations}}
\newcommand{\ES}{\end{subequations}}

\newcommand{\R}{\mathbb{R}}
\newcommand{\E}{\mathbb{E}}
\newcommand{\Pp}{\mathbb{P}}
\newcommand{\CN}{\mathcal{CN}}

\newcommand{\diag}{\mathrm{diag}}

\usepackage{indentfirst,csquotes}

\hypersetup{ colorlinks=true, linkcolor=black, filecolor=black, urlcolor=black }

\usepackage{lipsum}

\begin{document}
\author{%
  Songbin Liu%
  \thanks{Songbin Liu and Hanyu Yang contributed equally to this work.}%
  \textsuperscript{\,\,\,,}%
  \thanks{Department of Electrical Engineering, Columbia University.}%
  \quad Hanyu Yang\footnotemark[1]\textsuperscript{\,\,\,,}%
  \thanks{Department of Statistics and Data Science, School of Management, Fudan University.}%
  \quad Junren Chen%
  \thanks{Department of Statistics, Columbia University.}%
  \quad Arian Maleki\footnotemark[2]\textsuperscript{\,\,\,,}\footnotemark[4]%
}
\date{}
\title{Fundamental Limits of Coherent Imaging in The Presence of Atmospheric Turbulence} 
\date{\today}
\maketitle

\begin{abstract}
We study the fundamental limits of jointly estimating the signal amplitude and the atmospheric phase distortion in complex-valued multi-look coherent imaging, where the measurements are corrupted by multiplicative speckle and additive Gaussian noise, all looks share a known Gaussian sensing matrix, and both parameters belong to structured classes.
By analyzing a constrained maximum likelihood estimator, we establish an upper bound on the joint minimax risk that holds uniformly over the admissible parameter classes, while oracle reductions and Fano-type arguments yield worst-case lower bounds over the admissible parameter classes.
These results quantify the effects of structural complexity, signal and measurement dimensions, the number of looks, and the noise level, and establish matching minimax rates in several settings, with a remaining gap for phase estimation at low noise when the measurement dimension exceeds the signal dimension.
We also construct a two-stage estimator for signal recovery under arbitrary diagonal sign distortions
in the real-valued analogue. Under additional conditions, it recovers the signs up to a global sign with high probability and yields a uniform signal risk bound without a distortion-complexity term.
Numerical experiments illustrate the dependence of the joint, signal, and phase errors on the model parameters and provide further support for the theoretical analysis.
\end{abstract}

\section{Introduction}
\label{sec:Introduction}

\subsection{Main challenge and related work}
 Coherent imaging systems, such as synthetic aperture radar (SAR) and digital holography, are active sensing modalities that illuminate a scene with a controlled electromagnetic waveform and measure the reflected field. In many applications, a coherent laser source is used to illuminate the target, enabling the acquisition of both amplitude and phase information. By exploiting this phase information, coherent imaging systems can achieve high-resolution image reconstruction. However, when imaging through the atmosphere, random fluctuations in the refractive index caused by turbulence introduce spatially varying phase distortions that degrade image quality and complicate the reconstruction process. Consequently, a central challenge in coherent imaging is the recovery of the underlying scene from measurements corrupted by atmospheric turbulence and measurement noise. 
 
The main goal of this paper is to study the fundamental limits of coherent imaging systems in the presence of atmospheric turbulence. 
Mathematically speaking, the objective of such imaging systems is to recover an unknown complex-valued signal $\bm{x}_o \in \mathbb{C}^n$ from measurement $\{\bm{y}\}$, which is modeled as
\begin{equation}
    \bm{y} = \bm{\Phi}_o \bm{A} \bm{X}_o \bm{w} + \bm{z},
\end{equation}
where $\bm{X}_o = \operatorname{diag}(\bm{x}_o)$. $\bm{A} \in \mathbb{C}^{m \times n}$ represents the forward operator of the imaging system, while $\bm{\Phi}_o = \operatorname{diag}(\bm{\phi}_o) \in \mathbb{C}^{m \times m}$ models the unknown phase distortion induced by atmospheric turbulence. The vectors $\bm{w} \sim \mathcal{CN}(\bm{0}, \bm{I}_n),$ 
represent multiplicative noise, commonly referred to as speckle noise, which arises from surface roughness of the object relative to the wavelength of the illuminating light. Furthermore,
$\bm{z} \sim \mathcal{CN}(\bm{0}, \sigma_z^2 \bm{I}_m) $
denotes additive measurement noise, typically introduced by the imaging sensors. We should emphasize that the phase errors satisfy $    |\phi_{o,i}| = 1, \quad \forall i \in [m].$ The main objective of the imaging system is then to recover $\bm{x}_o$ from the measurement $\bm{y}$. 
    
In practice, the presence of significant noise (specially the speckle noise and phase errors) limits the quality of the reconstructed images, making single-shot estimates inadequate for many real-world applications. To improve reconstruction performance, practitioners commonly acquire multiple sets of measurements by varying certain imaging parameters, such as the illumination angle or frequency. Each acquisition produces a new set of measurements, referred to as a \emph{look} or \emph{shot}. This multi-look acquisition process can be modeled as
\begin{equation}\label{eqn:def_y_ell}
\bm{y}_l
=
\bm{\Phi}_o \bm{A} \bm{X}_o \bm{w}_l
+
\bm{z}_l,
\qquad l = 1, \ldots, L,
\end{equation}
where $L$ denotes the number of looks (or shots), and $\bm{y}_l \in \mathbb{C}^m$ represents the measurement vector collected during the $l^{\rm th}$ acquisition. The vectors $\{\bm{w}_l\}_{l=1}^L$ model independent realizations of the multiplicative speckle noise, while $\{\bm{z}_l\}_{l=1}^L$ represent independent additive measurement noise. 

A large body of work in signal processing, computational imaging, and physics communities, are dedicated to developing algorithms for recovering $\bm{x}_o$ from $\{\bm{y}_l\}_{\ell=1}^L$. Most of these classical methods attempt to address this scenario using a two-stage approach: (i) Estimate or compensate for the unknown phase distortion $\Phi_\ell$; (ii) Apply a transformation to convert the problem into a classical despeckling task, i.e. recovery image from measurements $\bm{y}_\ell = \bm{X}_o \bm{w}_\ell$. Two widely used strategies for phase estimation in this context are: (1) Adaptive Optics (AO) \cite{jang2015holographic, kim2012adaptive, roddier1999adaptive, booth2007adaptive, hampson2021adaptive}:AO systems utilize hardware components, such as wavefront sensors (e.g., Shack-Hartmann) and deformable mirrors, to measure and correct phase distortions in real time.  While effective under moderate turbulence, AO requires specialized, often expensive hardware and precise calibration. (ii) Digital Shearing Methods \cite{banet2023speckle, thurman2008phase}: These methods estimate phase distortion by analyzing the autocorrelation structure across shifted versions of the received measurements. Once estimated, $\Phi_\ell$ is inverted or compensated for, allowing the use of despeckling techniques. More recently researchers have started using optimization-based techniques for this purpose \cite{pellizzari2017phase,pellizzari2019imaging,pellizzari2023coherent,bate2021model,pellizzari2020coherent,pellizzari2017optically,chen2026monte,chen2026maximum,chen2025multilook}. In such approaches, the likelihood of $\{\bm{y}_\ell\}_{\ell=1}^L$ is written as the function of $\bm{\Phi}_o$ and $\bm{X}_o$. The likelihood  is often timed a nonconvex. However, researchers have developed different algorithms for finding stationary point of the likelihood functions. 
Related developments include model-based reconstruction through distributed-volume aberrations \cite{pellizzari2019imaging} and coherent plug-and-play reconstruction with learned despeckling priors \cite{pellizzari2023coherent}. CRISP jointly estimates optical aberrations and reconstructs speckle-free reflective scenes from pupil-modulated intensity measurements \cite{pellizzarispeckle}. For coherent despeckling, implicit neural representations have also been combined with complex Gaussian covariance likelihoods \cite{ziemann2026implicit}.

Many of these approaches explicitly estimate or compensate for phase distortions, either before recovering the signal or jointly with it. In many applications, however, only $\bm{x}_o$ is of interest, and $\bm{\Phi}_o$ is a nuisance parameter. A useful point of comparison is classical phase retrieval, where a fixed signal can be recovered, up to a global phase, from the magnitudes of its linear measurements under suitable conditions \cite{fienup1982phase,shechtman2015phase, balan2006signal,candes2014solving}. This motivates asking whether signal-only estimation can avoid a dependence on the complexity of the unknown phase class. However, existing phase retrieval results do not apply directly: each look contains a different latent speckle field $\bm{X}_o\bm{w}_\ell$, the common amplitude $\bm{x}_o$ enters the observation law through the covariance, and the same phase distortion is shared across all looks.

Despite substantial progress in the development of imaging systems and reconstruction algorithms for coherent imaging under atmospheric turbulence, many of the fundamental theoretical questions surrounding this problem remain unresolved. In particular, a central question concerns the fundamental limits of reconstruction accuracy as a function of the signal dimension, the number of measurements, the number of looks, the noise level, and the structural properties of both the underlying signal and the phase distortions induced by turbulence. Understanding these limits is important not only for characterizing the ultimate performance achievable by such imaging systems, but also for quantifying the gap between existing reconstruction algorithms and the information-theoretic optimum. Such insights can, in turn, guide the design of improved algorithms and system architectures that more effectively exploit the available measurements.

As discussed above, the fundamental limits of coherent imaging in the presence of atmospheric turbulence remain largely unexplored. The works most closely related to ours are a series of recent papers by some of the present authors and their collaborators \cite{zhou2022compressed,xing2025minimax,chen2025multilook,chen2024bagged}, which investigate the fundamental limits of coherent imaging systems described by
\[
\bm{y} = \bm{A}\bm{X}_o\bm{w} + \bm{z}.
\]

There are two key differences between these works and the present paper. First, the models considered in \cite{zhou2022compressed,xing2025minimax,chen2025multilook,chen2024bagged} are real-valued, with both the forward operator and the noise assumed to be real. In contrast, coherent imaging naturally involves complex-valued measurements and phase information. Second, the effects of atmospheric turbulence, and the resulting phase distortions, are not considered in these prior works. We address these limitations by studying joint signal and phase recovery in the complex-valued model. We further ask whether the complexity of the unknown phase distortion necessarily affects the accuracy of signal-only estimation, and provide a partial answer for a real-valued analogue with arbitrary sign distortions.

As we will demonstrate, incorporating turbulence-induced phase distortion introduces significant new challenges to the theoretical analysis. More importantly, these phase distortions fundamentally alter the statistical structure of the measurements and can substantially impact the achievable reconstruction accuracy. Consequently, the fundamental limits derived in the turbulence-free setting do not directly extend to the imaging model considered here. 

\paragraph{Contributions.}
Our main contributions are threefold.

\begin{enumerate}
    \item \textbf{Upper bound for joint signal and phase recovery.}
    We establish an upper bound on the joint minimax risk for estimating the signal amplitude and phase distortion in the complex-valued multi-look coherent imaging model \eqref{eqn:def_y_ell}, where all looks share the same known Gaussian sensing matrix and the same phase distortion.
    The bound applies to every pair of signal and phase classes whose covering numbers satisfy the polynomial bound of Definition~\ref{def:polynomial-complexity}, and it quantifies the dependence on their respective complexities, the signal and measurement dimensions, the number of looks, and the additive-noise variance.
    Our proof is based on an analysis of a constrained maximum likelihood estimator.

    \item \textbf{Lower bounds and sharpness of the rates.} Using oracle reductions and Fano-type arguments, we derive lower bounds that hold in the worst case over the families of signal and phase classes satisfying these covering-number bounds.
    The phase contribution exhibits different noise dependence depending on the relative sizes of the signal and measurement dimensions.
    Together with the upper bound, these results establish the minimax rates up to constant factors in several regimes, provided that the number of looks is small relative to the dimensions, as is typical in imaging applications.
    A gap remains for phase estimation when $m \ge n$ and the noise level is low.

    \item \textbf{Signal-only recovery under arbitrary sign distortions.}
    We also consider the real-valued analogue of \eqref{eqn:def_y_ell}, in which the sensing matrix, the speckle, and the noise are real-valued and the phase distortion reduces to an arbitrary unknown diagonal matrix with entries $\pm 1$.
    For this model, we establish a signal-only upper bound that depends on the complexity of the signal class, without a phase-complexity term. Under additional conditions on the dimensions and the number of looks, the bound holds uniformly over all diagonal sign matrices.
    Our proof uses a two-stage estimator that first recovers the unknown signs up to a global sign with high probability and then estimates the signal amplitude.
\end{enumerate}

\subsection{Notation}

Throughout the paper, we use bold uppercase letters, such as $\mA$, $\mX$, and $\mSigma$, to denote matrices, and bold lowercase letters, such as $\vx$, $\vy$, and $\vz$, to denote vectors. 
For a vector $\vx$, $x_i$ denotes its $i$-th entry, $|\vx|$ denotes the entrywise modulus, and $\operatorname{diag}(\vx)$ denotes the diagonal matrix with diagonal entries given by $\vx$. 
For a matrix $\mA$, $\mA^H$ denotes its conjugate transpose, $\sigma_{\max}(\mA)$ and $\sigma_{\min}(\mA)$ denote its largest and smallest singular values, and $\|\mA\|_2$ and $\|\mA\|_{\mathrm{HS}}$ denote its spectral and Hilbert--Schmidt(HS) norms, respectively. 
We write $\mathcal{CN}(0,\mSigma)$ for the circularly symmetric complex Gaussian distribution with mean zero and covariance $\mSigma$. 
Given two positive sequences ${a_n}$ and ${b_n}$, we write $a_n=O(b_n)$ if there exists a constant $C>0$ such that $a_n\le Cb_n$ for all sufficiently large $n$, and $a_n=\Omega(b_n)$ if $b_n=O(a_n)$. 
We write $a_n=\Theta(b_n)$ if both $a_n=O(b_n)$ and $a_n=\Omega(b_n)$ hold. 
Unless otherwise specified, the hidden constants may depend only on fixed model parameters and not on the problem dimensions or sample size.
Constants denoted by $c,C,c_1,C_1,\ldots$ may change from line to line.

\subsection{Organization}

The remainder of the paper is organized as follows.
Section~\ref{sec:main results} introduces the observation model, the structured signal and phase classes, and the joint loss function. It then presents the joint minimax upper and lower bounds, discusses their implications, and gives a signal-only upper bound for the real-valued model.
Section~\ref{sec:proof upper bound} proves the joint upper bound by analyzing a constrained maximum likelihood estimator, controlling likelihood fluctuations, and relating covariance discrepancies to the joint loss.
Section~\ref{sec:proof lower bound} establishes the joint lower bounds through oracle reductions and Fano-type arguments.
Section~\ref{sec:experiments} presents numerical experiments that illustrate the theoretical predictions and discusses the open questions they raise.
The appendices contain the identifiability result, auxiliary technical lemmas, and additional proofs, including the proof of the signal-only upper bound.

\section{Main results}
\label{sec:main results}

\subsection{Model, loss, and assumptions}
\label{sec:model and assump}
We now formalize the multi-look coherent imaging model introduced in Section~\ref{sec:Introduction}.

\subsubsection{Measurement model}\label{ssec:measurementmodel}

Let \(\mA\in\mathbb C^{m\times n}\) be the sensing matrix whose entries \(A_{ij}\) are drawn independently from \(\mathcal{CN}(0,1)\). The matrix \(\mA\) is known to the estimator and fixed across all looks.
Under the Gaussian speckle model specified below, coordinatewise phase rotations of the underlying complex-valued signal can be absorbed into the speckle vectors without changing their joint distribution.
The observation law therefore depends on the signal only through its entrywise amplitude. We use \(\vx_o\) to denote this amplitude vector and impose the bounds
\begin{equation}
\label{eq:amplitude_definition}
\vx_o\in[x_{\min},x_{\max}]^n,
\qquad
\mX_o:=\operatorname{diag}(\vx_o),
\end{equation}
where \(0<x_{\min}<x_{\max}<\infty\) are fixed constants.
For \(\ell=1,\ldots,L\), the observations are
\begin{equation}
\label{eq:multilook_model}
\vy_\ell = \mPhi_o\mA\mX_o\vw_\ell+\vz_\ell,
\end{equation}
where
\[
\mPhi(\vphi)
:=
\operatorname{diag}
\bigl(e^{i\phi_1},\ldots,e^{i\phi_m}\bigr),
\qquad
\mPhi_o:=\mPhi(\vphi_o),
\]
and \(\vphi_o\in\mathbb R^m\) parameterizes the unknown diagonal phase distortion. The speckle and additive-noise vectors satisfy
\[
\vw_\ell\sim\mathcal{CN}(0,\mI_n),
\qquad
\vz_\ell\sim\mathcal{CN}(0,\sigma_z^2\mI_m).
\]
All these vectors are mutually independent and independent of \(\mA\).
We assume that the noise variance \(\sigma_z^2>0\) is known. The unknown parameters are deterministic and do not depend on the realization of \(\mA\).

\subsubsection{Signal and phase error model}\label{ssec:signalmodel}
As we will discuss later, the existing structure in the signal $\bm{x}_o$ and phase error $\bm{\Phi}_o$ has a major impact on the quality of the reconstruction of $\bm{x}_o$. Following the framework of~\cite{xing2025minimax}, we describe this structure through covering numbers.
For a nonempty compact set $\mathcal{C}\subset\mathbb R^d$, let $N_\epsilon(\mathcal{C})$ denote its covering number under $\ell_2$ metric, namely the least number of $\ell_2-$balls covering $\mathcal{C}$.
\begin{definition}\label{def:polynomial-complexity}
We say that a compact set \(\mathcal{C}\subset\mathbb{R}^d\) has polynomial complexity of order \(k\) if there exist constants \(a>0\) and \(b\geq 0\), independent of \(k\) and \(d\), such that
\begin{equation}
N_\epsilon(\mathcal{C})
\le
\left(\frac{a d^b}{\epsilon}\right)^k,
\qquad \forall \epsilon\in(0,1].
\end{equation}
\end{definition}
The parameter \(k\) can be interpreted as the effective dimension of the class. We impose this condition separately on the signal-amplitude class and the phase-angle class.

\begin{itemize}

\item We assume that the signal amplitude belongs to a compact set $\vx_o\in\mathcal C_x\subset [x_{\min},x_{\max}]^n$ and that \(\mathcal C_x\) has polynomial complexity of order \(k_x\). Namely, there exist constants \(a_x>0\) and \(b_x\ge 0\) such that
\[
N_\epsilon(\mathcal C_x)
\le
\left(\frac{a_x n^{b_x}}{\epsilon}\right)^{k_x},
\qquad \forall \epsilon\in(0,1].
\]

\item We parameterize the unknown diagonal phase distortion by a real-valued angle vector
\(\vec{\phi}_o=(\phi_{o,1},\ldots,\phi_{o,m})\in\mathbb R^m\), and the corresponding phase distortion matrix is 
$ \mPhi(\vec{\phi}_o) =
\operatorname{diag}\big(e^{i\phi_{o,1}},\ldots,e^{i\phi_{o,m}}\big)$. 
Since multiplying \(\mPhi(\vec{\phi}_o)\) by a global phase does not change the covariance model, the phase angles are identifiable only up to a global additive shift. We remove this ambiguity by fixing the gauge \(\phi_{o,1}=0\). Choosing the principal representatives of the phase angles in \([0,2\pi)\), we assume that
\[\vec{\phi}_o\in
\mathcal C_\phi
\subset
\{\vec{\phi}\in[0,2\pi)^m:\phi_1=0\},\]
where \(\mathcal C_\phi\) is compact and has polynomial complexity of order \(k_\phi\). Namely, there exist constants \(a_\phi>0\) and \(b_\phi\ge 0\) such that
\[
N_\epsilon(\mathcal C_\phi)
\le
\left(\frac{a_\phi m^{b_\phi}}{\epsilon}\right)^{k_\phi},
\qquad \forall \epsilon\in(0,1].
\]

\end{itemize}

Accordingly, for \(a>0\), \(b\ge0\), and \(k\ge1\), let \(\mathcal F_x(a,b,k,n)\) denote the family of nonempty compact sets \(\mathcal C\subset[x_{\min},x_{\max}]^n\) that have polynomial complexity of order \(k\) with constants \(a\) and \(b\), and let \(\mathcal F_\phi(a,b,k,m)\) denote the family of nonempty compact sets \(\mathcal C\subset\{\vphi\in[0,2\pi)^m:\phi_1=0\}\) with the same property.
Thus \(\mathcal C_x\in\mathcal F_x(a_x,b_x,k_x,n)\) and \(\mathcal C_\phi\in\mathcal F_\phi(a_\phi,b_\phi,k_\phi,m)\).
We assume the parameter classes are known to the estimator and are fixed independently of \(\mA\).
The constants \(a_x,b_x,a_\phi,b_\phi\) do not depend on the dimensions or complexity orders.

\begin{remark}
The polynomial-complexity assumption covers a broad range of structured parameter classes. Examples include bounded piecewise constant and piecewise polynomial classes, bounded transform-sparse signal classes, and low-dimensional generative parameterizations with controlled parameter domains and Lipschitz constants. Further details and proofs of the corresponding covering-number bounds can be found in \cite{xing2025minimax}.
\end{remark}

\begin{remark}
Intuitively speaking, $k_x$ denotes the complexity of the signal class. If a signal class is not structured at all, then we expect $k_x \approx n$. On the other hand, for natural images one expects $k_x \ll n$.  The growth conditions on the complexity parameters are discussed in more detail in Section~\ref{subsec:minimax_bounds}.
\end{remark}

\begin{remark}
The complexity parameter $k_\phi$ associated with the phase errors reflects the strength of the atmospheric turbulence. Under weak turbulence, one would expect $k_\phi \ll m$, whereas under very strong turbulence, $k_\phi$ may approach $m$.\end{remark}

With these two parameter classes, we define the joint parameter space by
\[
\Theta:=\mathcal C_x\times\mathcal C_\phi,
\]
and write \(\vtheta_o:=(\vx_o,\vphi_o)\in\Theta\) for the true parameter.

\subsubsection{Joint loss for reconstruction} 

We now formulate the joint estimation problem and define the loss function used throughout the paper.  For any candidate parameter \(\vtheta=(\vx,\vec{\phi})\in\Theta\), define the precision matrix
\[
\mSigma(\vtheta)
:=
\left(
\mPhi(\vphi)\mA\operatorname{diag}(\vx^2)\mA^H\mPhi(\vphi)^H
+
\sigma_z^2\mI_m
\right)^{-1},
\]
where \(\vx^2\) denotes entrywise squaring. The inverse exists because \(\sigma_z^2>0\). Conditioned on \(\mA\), the observations are independent and identically distributed as 
\[
\vy_\ell\mid\mA
\sim
\mathcal{CN}\left(0,\mSigma(\vtheta_o)^{-1}\right),
\qquad
\ell=1,\ldots,L.
\]
Thus, conditioned on \(\mA\), the signal amplitude and phase distortion enter the statistical model through the covariance matrix \(\mSigma(\vtheta_o)^{-1}\). 

We now define the loss function and the corresponding minimax risk for joint estimation of the signal amplitude and phase distortion. 
For the signal amplitude, we use the normalized squared Euclidean loss 
\begin{equation} \label{eq:signal_loss} 
d_x^2(\hat\vx,\vx_o) := \frac{1}{n}\|\hat\vx-\vx_o\|_2^2. 
\end{equation} 
For a phase estimate \(\widehat{\mPhi}\) taking values in the set of diagonal unitary matrices, define
\begin{equation} \label{eq:phase_loss}
d_\Phi^2(\hat\mPhi,\mPhi_o) := \min_{\gamma\in[0,2\pi)} \left\| \hat\mPhi-e^{i\gamma}\mPhi_o \right\|_F^2 = 2\left( m-\left|\operatorname{Tr}(\hat\mPhi^H\mPhi_o)\right| \right).
\end{equation}
This loss is invariant to a global phase shift and avoids the dependence on a particular representation of the phase angles modulo \(2\pi\). 
We combine the two components into the normalized joint loss 
\begin{equation} \label{eq:joint_loss} 
d_{\mathrm{joint}}^2 \bigl( (\hat\vx,\hat\mPhi), (\vx_o,\mPhi_o) \bigr) := \frac{1}{n}\|\hat\vx-\vx_o\|_2^2 + \frac{1}{m}d_\Phi^2(\hat\mPhi,\mPhi_o). 
\end{equation}
For parameters \(\vtheta,\vtheta'\in\Theta\), we also write \(d_{\mathrm{joint}}^2(\vtheta,\vtheta')\) for this loss evaluated at their corresponding amplitude--phase matrix pairs.

The main theoretical question we would like to address in this paper is the following:

\begin{center} 

\textbf{Given the measurement model described in Section~\ref{ssec:measurementmodel}, together with the signal and phase-error models introduced in Section~\ref{ssec:signalmodel}, what is the best performance that any algorithm can achieve in reconstructing $\bm{\Phi}_o$ and $\bm{x}_o$?}

\end{center}

To address this challenge we study the following minimax risk: \begin{equation} \label{eq:joint_minimax_risk} 
R_{\mathrm{joint}} \bigl( \mathcal C_x,\mathcal C_\phi; m,n,L,\sigma_z \bigr) := \inf_{\hat\vx,\hat\mPhi} \sup_{\substack{ \vx_o\in\mathcal C_x\\ \vec{\phi}_o\in\mathcal C_\phi }} \mathbb E \left[ d_{\mathrm{joint}}^2 \bigl( (\hat\vx,\hat\mPhi), (\vx_o,\mPhi(\vec{\phi}_o)) \bigr) \right], 
\end{equation} 
where the infimum is over all measurable estimators based on \((\mA,\vy_1,\ldots,\vy_L)\), with \(\widehat{\vx}\in[x_{\min},x_{\max}]^n\) and \(\widehat{\mPhi}\) a diagonal unitary matrix.
The expectation is taken over the sensing matrix and all speckle and additive-noise vectors, with the true parameter held fixed.

In words, under this minimax framework, the performance of an estimator is evaluated based on its performance for the most challenging combination of signal and phase error.

\subsection{Minimax upper and lower bounds}
\label{subsec:minimax_bounds}

We now present the upper and lower bounds for joint minimax risk of the signal amplitude and phase distortion. Their interpretation and sharpness are discussed in Subsection~\ref{subsec:interpretation_rates}.
We use the class families \(\mathcal F_x(a_x,b_x,k_x,n)\) and \(\mathcal F_\phi(a_\phi,b_\phi,k_\phi,m)\) of Section~\ref{ssec:signalmodel}. Throughout this subsection, \(a_x,b_x,a_\phi,b_\phi\) are fixed constants.

\paragraph{Upper bound.}
We begin with an upper bound that holds uniformly over all admissible signal and phase classes.

\begin{theorem}[Joint minimax upper bound]
\label{thm:joint_upper_bound}
Let $ \mathcal C_x\in\mathcal F_x(a_x,b_x,k_x,n)$ and $ \mathcal C_\phi\in\mathcal F_\phi(a_\phi,b_\phi,k_\phi,m) $. Assume that \(m,n\ge2\), \(k_x,k_\phi\ge1\), and $mL\le n^4\log n$, $L\le m^4\log m$.
Then there exists a constant \(C>0\), depending only on \(x_{\min},x_{\max},a_x,b_x,a_\phi,b_\phi\), such that
\begin{equation}
\label{eq:joint-upper-bound}
R_{\mathrm{joint}}(\mathcal C_x,\mathcal C_\phi;m,n,L,\sigma_z)
\le
C\min\left\{
(k_x\log n+k_\phi\log m)
\left[
\frac{\max\{\sigma_z^4,m^2,n^2\}}{m^2nL}
+
\frac{m+n}{m^2n}
\right], 1
\right\}.
\end{equation}
\end{theorem}
The proof is given in Section~\ref{sec:proof upper bound} and is based on an analysis of a constrained maximum likelihood estimator.

\paragraph{Lower bounds.}
We next establish worst-case lower bounds over the admissible families of signal and phase classes. Both lower-bound results use the following common assumptions.
Fix \(\varepsilon\in(0,1/2)\), and suppose that
\[
\log n=\Theta(\log m),\quad
\log L=O(\log m),\quad
1\le k_x\le n^{1-2\varepsilon},\quad
1\le k_\phi\le m^{1-2\varepsilon}.
\]
We also assume that $a_x\ge x_{\max}-x_{\min}$, $b_x\ge1$, $a_\phi\ge a_{\phi,0}$, $b_\phi\ge b_{\phi,0}$ where \(a_{\phi,0},b_{\phi,0}>0\) are absolute constants, and impose the rate conditions
\begin{equation}
\label{eq:common-lower-bound-conditions}
\max\{\sigma_z^4,m^2,n^2\}\,k_x\log n
\le m^2n^{1-\varepsilon}L, \quad
\max\{\sigma_z^4,m^2,n^2\}\,k_\phi\log m
\le m^{2-\varepsilon}nL.
\end{equation}

\begin{theorem}[Lower bound when \(n\ge(1+\tau)m\)]
\label{thm:joint-lb-n-ge-m}
Suppose that the common assumptions hold. Fix \(\tau>0\) and assume \(n\geq(1+\tau)m\). Then there exists a constant \(C=C(\tau,\varepsilon,x_{\min},x_{\max})>0\) such that, for all sufficiently large \(m,n\),
\begin{equation}
\label{eq:joint_lower_bound_n_ge_m}
\sup_{\substack{
\mathcal C_x\in\mathcal F_x(a_x,b_x,k_x,n)\\
\mathcal C_\phi\in\mathcal F_\phi(a_\phi,b_\phi,k_\phi,m)
}}
R_{\mathrm{joint}}
(\mathcal C_x,\mathcal C_\phi;m,n,L,\sigma_z)
\geq
C\,
\frac{\max\{\sigma_z^4,n^2\}}{m^2nL}
\max\{k_x\log n,k_\phi\log m\}.
\end{equation}
\end{theorem}

\begin{theorem}[Lower bound when \(m\ge n\)]
\label{thm:joint-lb-m-ge-n}
Suppose that the common assumptions hold and that \(m\geq n\).
Then there exists a constant \(C=C(\varepsilon,x_{\min},x_{\max})>0\) such that, for all sufficiently large \(m,n\),
\begin{equation}
\label{eq:joint_lower_bound_m_ge_n}
\begin{aligned}
&\sup_{\substack{
\mathcal C_x\in\mathcal F_x(a_x,b_x,k_x,n)\\
\mathcal C_\phi\in\mathcal F_\phi(a_\phi,b_\phi,k_\phi,m)
}}
R_{\mathrm{joint}}
(\mathcal C_x,\mathcal C_\phi;m,n,L,\sigma_z)
\geq
C\,
\frac{\max\{\sigma_z^4,m^2\}}{m^2nL}
\max\left\{
k_x\log n,\,
\frac{\sigma_z^2}{\max\{\sigma_z^2,m\}}
k_\phi\log m
\right\}.
\end{aligned}
\end{equation}
\end{theorem}
The proofs are given in Section~\ref{sec:proof lower bound}. We reduce joint estimation to signal- and phase-oracle problems and derive the corresponding lower bounds using Fano-type arguments.

\begin{remark} \label{rem:near-square}
For a fixed \(\tau>0\), Theorems~\ref{thm:joint-lb-n-ge-m} and~\ref{thm:joint-lb-m-ge-n} provide lower bounds for \(n\ge(1+\tau)m\) and \(m\ge n\), respectively.
The first condition requires the ratio \(n/m\) to stay a fixed positive distance above one, while the second includes the square case \(m=n\). For this choice of \(\tau\), the intermediate range \(m<n<(1+\tau)m\) is not covered by these results.
\end{remark}

\paragraph{Discussion of the assumptions.}
Before interpreting these rates, we discuss the assumptions and their roles in the proofs.

\begin{enumerate}
    \item \emph{Conditions on the number of looks.} 
    The upper-bound conditions \(mL\le n^4\log n\) and \(L\le m^4\log m\) are mild when the signal and measurement dimensions are large and comparable.
    As discussed in \cite{xing2025minimax}, acquiring more than \(100\) independent looks is often difficult in practice, whereas the signal dimension can reach hundreds of thousands or millions.
    In such settings, the two upper bounds on \(L\) are therefore well above the number of independent looks that are typically available.

    \item \emph{Growth of the dimensions and sample size.}
    The condition \(\log n=\Theta(\log m)\) in lower bound allows a substantially wider range of dimensions than a fixed aspect-ratio assumption.
    For example, it holds whenever \(n=m^\alpha\) or any fixed \(\alpha>0\), covering both \(n>m\) and \(m>n\).
    The condition \(\log L=O(\log m)\) is also mild in applications where the number of independent looks is limited, as discussed above.
    It holds automatically when \(L\) is bounded and also permits \(L\) to increase polynomially with \(m\).

    \item \emph{Signal and phase complexities.}
    The restrictions \(k_x\le n^{1-2\varepsilon}\) and \(k_\phi\le m^{1-2\varepsilon}\) reflect the setting in which the signal and phase have low-dimensional descriptions relative to their ambient dimensions.
    For example, block-constant signals and phases with a fixed number of blocks satisfy these conditions as the dimensions increase.
    The conditions also allow the complexity orders to grow logarithmically or as fixed sublinear powers of the dimensions, provided that \(\varepsilon>0\) is chosen sufficiently small.
    They are therefore mild for parameter classes whose structural complexity grows slowly with resolution.

    \item \emph{Covering-number constants.}
    The lower thresholds on \(a_x,b_x,a_\phi,b_\phi\) are sufficient conditions for the classes constructed in the lower-bound proofs to belong to the admissible families.
    The signal thresholds agree with those in \cite[Theorem~2.15]{xing2025minimax}. For the phase class, the construction in Appendix~\ref{pf:lem:phase_class_construction} satisfies
    \[
    N_\rho(\Gamma_\phi)
    \le \left(\frac{2em}{\rho}\right)^{k_\phi},
    \qquad 0<\rho\le1.
    \]
    Consequently, one may take \(a_{\phi,0}=2e\) and \(b_{\phi,0}=1\) in the lower-bound assumptions, so
    the conditions \(a_\phi\ge a_{\phi,0}\) and \(b_\phi\ge b_{\phi,0}\) then ensure that \(\Gamma_\phi\) belongs to the admissible phase-class family.
    Once these constants are fixed, they do not change the stated dependence of the rates on the dimensions, the number of looks, or the noise level.

    \item \emph{Rate conditions.}
    The signal rate condition in \eqref{eq:common-lower-bound-conditions} is inherited from \cite[Theorem~2.15]{xing2025minimax}, which we use to establish the signal-oracle lower bound. The phase rate condition ensures that the packing radius chosen in our proof remains within the prescribed phase range.
    Together, these conditions specify a sufficient parameter range for the lower-bound arguments.

    These conditions are mild in several settings. For example, when \(m=n\) and \(\sigma_z^2\le n\), they reduce to
    \[
    k_x\log n\le n^{1-\varepsilon}L,
    \qquad
    k_\phi\log n\le n^{1-\varepsilon}L.
    \]
    The complexity assumptions \(k_x,k_\phi\le n^{1-2\varepsilon}\) imply both inequalities for all sufficiently large \(n\) and every \(L\ge1\).
    Thus, in this setting, the rate conditions impose no additional restriction on the number of looks.

\end{enumerate}

\subsection{Interpretation and comparison of the rates}
\label{subsec:interpretation_rates}

We now discuss how the bounds in Theorems~\ref{thm:joint_upper_bound}--\ref{thm:joint-lb-m-ge-n} depend on the model parameters and compare their orders.

\paragraph{Interpretation of the upper bound.}
When \(\sigma_z^2\leq\max\{m,n\}\), the dimensional terms dominate the additive-noise contribution in \(\max\{\sigma_z^4,m^2,n^2\}\).
Theorem~\ref{thm:joint_upper_bound} then gives, up to constant factors,
\begin{equation}
\label{eq:upper-bound-low-noise-simplified}
R_{\mathrm{joint}}
(\mathcal C_x,\mathcal C_\phi;m,n,L,\sigma_z)
\leq
\begin{cases}
\displaystyle
C\min\left\{
(k_x\log n+k_\phi\log m)
\left(\frac{1}{nL}+\frac{1}{mn}\right),1
\right\},
& m\geq n,
\\[10pt]
\displaystyle
C\min\left\{
(k_x\log n+k_\phi\log m)
\left(\frac{n}{m^2L}+\frac{1}{m^2}\right),1
\right\},
& n\geq m.
\end{cases}
\end{equation}
For fixed \(m\), the dimensional factor in the \(L\)-dependent term decreases with \(n\) when \(n\leq m\) and increases when \(n\geq m\). This factor is therefore smallest at \(n=m\), with the dependence through \(k_x\log n+k_\phi\log m\) considered separately.

\begin{itemize}

\item \emph{Signal and phase complexities.}
The signal and phase classes contribute through \(k_x\log n\) and \(k_\phi\log m\), respectively.
Compared with the signal-only bounds for the shared-sensing-matrix model in \cite{xing2025minimax}, it contains the additional term \(k_\phi\log m\) associated with estimating the unknown phase distortion.
This dependence raises the question of whether signal recovery alone must incur the same phase-complexity cost. The next subsection gives a partial answer by showing that this cost can be avoided in a real-valued analogue under additional conditions.
Section~\ref{sec:experiments} further examines the dependence of the signal and phase errors on phase complexity in the complex-valued model.

\item \emph{Additive-noise level.}
The \(L\)-dependent term in \eqref{eq:joint-upper-bound} has no explicit noise dependence when \(\sigma_z^2\le\max\{m,n\}\), and when \(\sigma_z^2 > \max\{m,n\}\) it is proportional to \(\sigma_z^4\). Heuristically, the noise makes nearby parameters harder to distinguish only when \(\sigma_z^2\) is comparable to or larger than the scale of the signal part of the observation covariance.

Conditioned on \(\mA\), each look \(\vy_\ell\) is a zero-mean complex Gaussian vector with covariance \(\mSigma(\vtheta_o)^{-1}=\mPhi_o\mA\mX_o^2\mA^H\mPhi_o^H+\sigma_z^2\mI_m\), so the parameters affect the observations only through this covariance.
For a scalar variance \(v\) and a small perturbation \(\delta\),
\[
D_{\mathrm{KL}}\bigl(\mathcal{CN}(0,v+\delta)\,\|\,\mathcal{CN}(0,v)\bigr)\approx\frac{\delta^2}{2v^2}.
\]
Thus, the information induced by a covariance change depends on its scale relative to the total variance. Analogously, in the vector case, the eigenvalues of the covariance play the role of \(v\).
The nonzero eigenvalues of the signal part \(\mPhi_o\mA\mX_o^2\mA^H\mPhi_o^H\) are typically of order \(\max\{m,n\}\), and changing the amplitude and phase parameters can induce covariance changes on the same scale, whereas the additive noise contributes \(\sigma_z^2\) to every eigenvalue. 
Hence, the noise does not change the leading order of total covariance eigenvalues when \(\sigma_z^2\lesssim\max\{m,n\}\), 
but reduces this information by a factor of order \(\max\{m,n\}^2/\sigma_z^4\) when \(\sigma_z^2\gg\max\{m,n\}\), consistent with the factor \(\max\{\sigma_z^4,m^2,n^2\}\) in the bound.
The same threshold appears in the model without phase distortion \cite[Remark~2.12]{xing2025minimax}.

When \(m>n\), however, the signal covariance has rank \(n\), and its range is the column space of \(\mPhi_o\mA\), which does not depend on \(\vx_o\).
Hence the variance of \(\vy_\ell\) is only \(\sigma_z^2\) in the \(m-n\) directions orthogonal to this space, whereas in this space it is typically of order \(\max\{\sigma_z^2,m\}\).
If \(\mPhi_o\) is replaced by another phase distortion matrix \(\mPhi\), the range of the signal covariance becomes the column space of \(\mPhi\mA\), so in general the covariance also changes in these \(m-n\) directions.
When \(\sigma_z^2\ll m\), these directions have much smaller variance under the true parameter than the signal-dominated directions, making such phase-induced covariance changes more distinguishable. This additional sensitivity to the noise level may help explain the factor \(\sigma_z^2/\max\{\sigma_z^2,m\}\) in the phase lower bound in Theorem~\ref{thm:joint-lb-m-ge-n}, as will discussed later.

\item \emph{Number of looks.}
The upper bound in \eqref{eq:joint-upper-bound} contains an \(L\)-dependent contribution and an \(L\)-independent remainder
\[
\frac{\max\{\sigma_z^4,m^2,n^2\}(k_x\log n+k_\phi\log m)}{m^2nL}, \qquad
\frac{(m+n)(k_x\log n+k_\phi\log m)}{m^2n}.
\]
The former dominates the latter when
\begin{equation}
\label{eq:upper-statistical-term-dominates}
L(m+n)\lesssim\max\{\sigma_z^4,m^2,n^2\}.
\end{equation}
As discussed above, the number of independent looks is typically limited and much smaller than the signal and measurement dimensions. In such settings, this condition is readily satisfied, and the \(L\)-dependent term dominates the remainder.

To interpret the \(L\)-independent term, recall that all looks share the same sensing matrix \(\mA\). Additional looks provide an increasingly accurate estimate of the same population covariance, while parameter recovery still depends on the inverse problem associated with this fixed matrix. This offers a heuristic explanation for the error floor in the current upper bound; see \cite{chen2025multilook} for a more detailed discussion.
However, the population covariance uniquely determines the parameters in identifiable settings, and our lower bounds contain no analogous remainder. Hence, whether the \(L\)-independent term is unavoidable therefore remains unclear, even in the signal-only model \cite{xing2025minimax}.
\end{itemize}

\paragraph{Comparison of the upper and lower bounds.}
We next compare the uniform upper bound with the worst-case lower bounds over the admissible class families, focusing on settings where both results apply. To compare their dependence on complexity, we use
\[
\max\{k_x\log n,k_\phi\log m\}
\leq
k_x\log n+k_\phi\log m
\leq
2\max\{k_x\log n,k_\phi\log m\}.
\]
Thus, combining the two complexity terms through their sum or maximum changes the rate only by a constant factor.

\begin{itemize}

\item \emph{Underdetermined regime (\(n\geq(1+\tau)m\)).}
For fixed \(\tau>0\), we have $\max\{\sigma_z^4,m^2,n^2\}=\max\{\sigma_z^4,n^2\}$. Hence the lower bound in Theorem~\ref{thm:joint-lb-n-ge-m} has the same order as the \(L\)-dependent term in the upper bound, including its dependence on the dimensions, noise level, and both complexity orders. Consequently, the bounds match in order whenever \eqref{eq:upper-statistical-term-dominates} holds, as is typical in imaging applications where the number of independent looks is small relative to the signal and measurement dimensions.
In particular, \(L\lesssim n\) is sufficient in this regime, while higher noise levels can allow matching for larger \(L\).

\item \emph{Overdetermined regime (\(m\geq n\)).}
In this regime, the phase contribution to the lower bound contains the additional factor
\[
\frac{\sigma_z^2}{\max\{\sigma_z^2,m\}},
\]
which leads to different behavior in the high- and low-noise settings. When \(\sigma_z^2 \gtrsim m\), this factor is of constant order, so the lower bound matches the \(L\)-dependent term in the upper bound. The full bounds therefore match in order whenever \eqref{eq:upper-statistical-term-dominates} holds.

At low noise, \(\sigma_z^2\ll m\), the factor becomes \(\sigma_z^2/m\). The phase contribution to the lower bound is then of order \(\sigma_z^2 k_\phi\log m/(mnL)\), compared with \(k_\phi\log m/(nL)\) in the \(L\)-dependent term of the upper bound. The signal contributions, however, retain the same order \(k_x\log n/(nL)\).
Consequently, the joint bounds still match when \(k_x\log n\gtrsim k_\phi\log m\) and \eqref{eq:upper-statistical-term-dominates} holds, because the signal-complexity term controls both bounds.
A gap can remain when \(k_\phi\log m\gg k_x\log n\), so that the phase-complexity term dominates the upper bound.

One limitation of the current upper-bound analysis comes from the comparison across noise levels.
For \(m\geq n\) and \(0<\sigma_z^2<m\), the proof transfers a risk bound from noise variance \(m\) to the smaller variance \(\sigma_z^2\).
The resulting bound therefore retains the larger-noise scale. Although this scale is consistent with the signal-only rates in \cite{xing2025minimax}, it may overlook a possible improvement in phase estimation at low noise, as suggested by the factor \(\sigma_z^2/m\) in the phase lower bound.

A heuristic interpretation comes from the additional phase constraints imposed by the rank-deficient signal covariance when \(m>n\). Its range is \(\mPhi_o\operatorname{range}(\mA)\), independently of the unknown positive signal amplitudes. To exploit this structure, let
\[
Q_A
=
\mI_m-\mA(\mA^H\mA)^{-1}\mA^H
\]
denote the orthogonal projector onto \(\operatorname{range}(\mA)^\perp\). Every look then satisfies
\[
Q_A\mPhi_o^H\vy_\ell
=
Q_A\mA\mX_o\vw_\ell + Q_A\mPhi_o^H\vz_\ell
=
Q_A\mPhi_o^H\vz_\ell.
\]
In the absence of additive noise, these relations give exact constraints on the phase distortion without requiring knowledge of the amplitudes or the speckle noise.
Hence these constraints suggest a phase-only criterion,
\[
\mathcal J(\vphi)
=
\frac1L\sum_{\ell=1}^L
\left\|
Q_A\mPhi(\vphi)^H\vy_\ell
\right\|_2^2.
\]
Under the true phase, the residual $\mathcal J(\vphi_o)$ depends only on the  additive noise, whereas an incorrect phase can introduce a signal contribution to \(\mathcal J(\vphi)\).
In particular, 
\[\mathbb E[\mathcal J(\vphi_o)\mid\mA]=\sigma_z^2\operatorname{Tr}(Q_A^2)
=\sigma_z^2(m-n).\]
Reducing the additive-noise level therefore makes these constraints more accurate, while amplitude estimation remains subject to speckle fluctuations. This provides a possible explanation for improved phase accuracy at low noise. The extent of this improvement depends on whether the available looks and the phase class provide enough constraints to distinguish different phase distortions stably.

\item \emph{Behavior near \(n=m\).}
Theorem~\ref{thm:joint-lb-m-ge-n} covers \(m\geq n\), including the square case \(m=n\). For \(n>m\), Theorem~\ref{thm:joint-lb-n-ge-m} requires \(n/m\geq1+\tau\) for a fixed \(\tau>0\). Thus, for a given \(\tau\), the intermediate range \(m<n<(1+\tau)m\) is not covered, as noted in Remark~\ref{rem:near-square}.
Choosing a smaller fixed \(\tau\) extends the covered range toward \(n=m\), but the lower-bound constant also depends on \(\tau\). The current results therefore do not provide a uniform rate comparison as \(n/m\) approaches one from above.
The experiments in Section~\ref{sec:experiments} also show a rapid change in phase error near \(n=m\) at low noise. Hence a sharper characterization of this transition would require more precise control of how the bounds depend jointly on the noise level and the distance of the aspect ratio from one, including the behavior of the constants as \(n/m\) approaches one from above.

\end{itemize}

\subsection{Signal-only estimation}
\label{sec:signal-main}

While the main focus of this paper is on the joint estimation error, we ask in this section the following  question: if we are only interested in the recovery of the signal $\vx_o$ (as in some real applications), can we attain a sharper error rate? While the joint upper bound \eqref{eq:joint-upper-bound} in
Theorem~\ref{thm:joint_upper_bound} depends on both $k_x$ and $k_\phi$, our intuition is that the signal-only estimation error  $d_x^2(\widehat{\vx},\vx_o)$ should be free of the complexity of the phase error class $k_{\phi}$. 
We give a rigorous result in this direction for the corresponding \emph{real-valued}
model, showing that, information-theoretically, an $k_\phi$-free upper bound on $d_x^2(\widehat{\vx},\vx_o)$ is attainable in some regime.

We start with the real-valued model in which the unknown phase distortion becomes an arbitrary diagonal sign matrix. Retain the compact amplitude class
$\mathcal C_x\in\mathcal F_x(a_x,b_x,k_x,n)$ from
Section~\ref{subsec:minimax_bounds}, with
$\mathcal C_x\subset[x_{\min},x_{\max}]^n$ and
$0<x_{\min}<x_{\max}<\infty$ as in \eqref{eq:amplitude_definition}.
Polynomial complexity is understood in the covering-number sense of
Section~\ref{sec:model and assump}. In the real-valued analogue of
\eqref{eq:multilook_model}, let $\mA\in\R^{m\times n}$ have independent
$\mathcal N(0,1)$ entries and be known and shared across all looks. We observe
\begin{equation}
 \vy_\ell=\mD_o\mA\mX_o\vw_\ell+\vz_\ell,\qquad
 \vw_\ell\sim \mathcal N(0,\mI_n),\qquad
 \vz_\ell\sim \mathcal N(0,\sigma_z^2\mI_m),\qquad \ell=1,\ldots,L,
 \label{eq:signal-model}
\end{equation}
where $\mX_o=\diag(\vx_o)$, $\sigma_z^2>0$, and all Gaussian noise
vectors are mutually independent and independent of $\mA$.
Write $\mD_o=\diag(\vec{\varepsilon}_o)\in
\mathcal D_m:=\{\diag(\vec{\varepsilon}):\vec{\varepsilon}\in\{\pm1\}^m\}$.
Only the pair $\{\mD_o,-\mD_o\}$ is identifiable, and no structural
restriction is imposed on the signs. Using the normalized signal loss in \eqref{eq:signal_loss}, write
$d_x^2(\widehat{\vx},\vx_o)=n^{-1}\|\widehat{\vx}-\vx_o\|_2^2$.

\paragraph{A two-stage estimator.} We find it technically challenging to derive sharper upper bound on $d_x^2(\widehat{\vx},\vx_o)$ by analyzing the MLE. As such, to sharpen the (information-theoretic) upper bound, we switch to a novel two-stage estimator. The idea is to separately estimate the phase error $\mD_o$ and the signal $\vx_o$.  Form
\begin{equation}
 \widehat{\mC}=\frac1L\sum_{\ell=1}^L\vy_\ell\vy_\ell^\top,\qquad
 \mG=\mA\mA^\top,\qquad
 \widehat W_{ij}=G_{ij}\widehat C_{ij}\ (i\ne j),\quad
 \widehat W_{ii}=0.
 \label{eq:signal-weights}
\end{equation}
For $i\ne j$, independence and Gaussian second moments give
$\E\widehat W_{ij}=\varepsilon_{o,i}\varepsilon_{o,j}\|\vx_o\|_2^2$.
Here the expectation averages over the sensing matrix and observations.
Thus these weights favor the true relative signs on average. Estimate them by
\begin{equation}
 \widehat{\vec{\varepsilon}}\in
 \operatorname*{arg\,max}_{\vec{\varepsilon}\in\{\pm1\}^m}
 \vec{\varepsilon}^\top\widehat{\mW}\vec{\varepsilon},\qquad
 \widehat{\mD}=\diag(\widehat{\vec{\varepsilon}}).
 \label{eq:signal-sync}
\end{equation}
In the real-valued model, our key insight is that the above program {\it exactly} recovers $\mD_o$ in some regime. As such, we are able to completely get rid of the impact of the unknown phase $\pm\mD_o$ (up to a sign), in the light of 
\[\widehat{\mD}\vy_\ell = \eta\big(\mA\mX_o\vw_\ell+\mD_o\vz_\ell\big),\qquad \forall \ell=1,...,L,\]
for some $\eta=\pm 1$. Therefore,
for the second stage we can let $\mathcal T_{\mA}$ be a signal estimator in  \cite[Theorem~2.14]{xing2025minimax}.
We choose it to take values in $[x_{\min},x_{\max}]^n$ and satisfy
$\mathcal T_{\mA}(-\mZ)=\mathcal T_{\mA}(\mZ)$ for every input
$\mZ$; Appendix~\ref{sec:signal-appendix} justifies this choice without increasing its risk.
With $\mY=(\vy_1,\ldots,\vy_L)$, set
\begin{equation}
 \widehat{\vx}=\mathcal T_{\mA}(\widehat{\mD}\mY).
 \label{eq:signal-estimator}
\end{equation}

\begin{theorem}[Signal-only upper bound]
\label{thm:signal-only}
Let $a_x,b_x$ be fixed and put $V:=\max\{\sigma_z^4,m^2,n^2\}$.
For every $B>0$, there exist constants $C_B,C>0$, depending only on
$B,x_{\min},x_{\max},a_x,b_x$, such that if
\begin{equation}
 \begin{aligned}
 &m\ge C_B,\qquad n\ge C_B\log(em),\\
 &L\ge C_B\frac{V}{mn}\log(em),\qquad
 mL\le n^4k_x\log n,
 \end{aligned}
 \label{eq:signal-conditions}
\end{equation}
then the estimator \eqref{eq:signal-estimator} satisfies
\begin{equation}
 \sup_{\substack{\vx_o\in\mathcal C_x\\ \mD_o\in\mathcal D_m}}
 \E d_x^2(\widehat{\vx},\vx_o)
 \le
 C\min\left\{
 \frac{Vk_x\log n}{m^2nL}
 +\frac{k_x\log m\log n}{m^2},\,1
 \right\}
 +C_Bm^{-B}.
 \label{eq:signal-bound}
\end{equation}
If $V\ge nL\log m$, the second term inside the minimum can be omitted
after adjusting $C$.
\end{theorem}

To see the improvement over previous upper bound, we ignore logarithmic factors and provide a   detailed comparison: while  the joint estimation upper bound in Theorem \ref{thm:joint_upper_bound} reads 
\[(k_x+k_\phi) \bigg[\frac{V}{m^2nL}+\frac{m+n}{m^2n}\bigg],\]
the signal-only upper bound in Theorem \ref{thm:signal-only}, up to a   higher order term $O(m^{-B})$ (we set $B$ large enough), is given by 
\[k_x\bigg[\frac{V}{m^2nL}+ \frac{1}{m^2}\bigg].\]
It is clear that we remove $k_\phi$ in the leading factor, along with a slight improvement of $\frac{m+n}{m^2n}$ to $\frac{1}{m^2}$.

Nonetheless, Theorem \ref{thm:signal-only} requires some additional constraints. The most notable one is
the lower bound on $L\ge C_B\frac{V\log(em)}{mn}$, which is needed so that exact recovery of $\pm\mD_o$ is achieved in the first stage. Here, we only comment that
it reduces to the mild condition $L\gtrsim\log m$ in the regime $m\asymp n$ and $\sigma_z^2\lesssim n$. 
 
Interestingly, 
the numerical experiments in the next section (Figure~\ref{fig:MSE_vs_Kphi}) also indicate that the signal-only estimation error may not depend on $k_\phi$.  

\section{Proof of Upper Bound}
\label{sec:proof upper bound}

In this section, we prove Theorem~\ref{thm:joint_upper_bound} by analyzing the constrained maximum likelihood estimator. We start from an inequality that relates the population likelihood gap to an empirical fluctuation term, then bound the two sides of this inequality and relate the resulting covariance bound to the joint estimation error.

\subsection{Likelihood and Basic Inequality}
\label{subsec:Basic inequality}

Let \(\vtheta_o=(\vx_o,\vphi_o)\) denote the true parameter we want to estimate. Then conditioned on \(\mA\), the observations \(\vy_1,\ldots,\vy_L\) under the model in Section~\ref{sec:model and assump} are i.i.d. complex Gaussian with covariance
\[
\vy_{\ell} | \mA \sim \mathcal{CN}(0,\mSigma(\vtheta_o)^{-1}), \qquad
\mSigma(\vtheta)^{-1}
=
\mPhi\mA\mX^2\mA^H\mPhi^H+\sigma_z^2\mI_m.
\]
Therefore, up to a constants, the negative log-likelihood per look is
\[
\mathcal L(\vtheta)
=
\log\det\big(\mSigma(\vtheta)^{-1}\big)
+
\frac1L\sum_{\ell=1}^L \vy_\ell^H\mSigma(\vtheta)\vy_\ell.
\]
Write \(\mSigma_o=\mSigma(\vtheta_o)\). Since
$
\E\left[ \frac1L\sum_{\ell=1}^L
\vy_\ell\vy_\ell^H \,\middle|\,\mA
\right] = \mSigma_o^{-1},
$
the conditional expected likelihood is
\begin{align}\label{def:exp-liklihood}
\overline{\mathcal L}(\vtheta)
:= \E[\mathcal L(\vtheta)\mid\mA]
= \log\det\big(\mSigma(\vtheta)^{-1}\big)
+ \Tr\big(\mSigma(\vtheta)\mSigma_o^{-1}\big).
\end{align}
Hence
\begin{align}
\label{eq:likelihood_decomp}
\mathcal L(\vtheta)
=
\overline{\mathcal L}(\vtheta)
+
\frac1L\sum_{\ell=1}^L
\vy_\ell^H\mSigma(\vtheta)\vy_\ell
-
\Tr\big(\mSigma(\vtheta)\mSigma_o^{-1}\big).
\end{align}

\begin{remark}[Uniqueness]\label{rem:uniqueness}
The population objective is minimized at the true parameter \(\vtheta_o\).
We give a sufficient condition for this minimizer to be unique in Appendix~\ref{app:proof_uniqueness}.
\end{remark}

Let $ \widehat{\vtheta} = (\widehat{\vx},\widehat{\vphi}) \in \arg\min_{\vtheta\in\Theta}\mathcal L(\vtheta) $ be a measurable constrained maximum likelihood estimator, and write \(\widehat{\mSigma}=\mSigma(\widehat{\vtheta})\).
Such a selection exists by compactness of $\Theta$ and continuity of the likelihood in $\vtheta$, since $\sigma_z^2>0$. We fix the true parameter $\vtheta_o$ independently of $\mA$ throughout the argument.
Using \(\mathcal L(\widehat{\vtheta})\leq\mathcal L(\vtheta_o)\) and \eqref{eq:likelihood_decomp}, we obtain
\begin{equation}\label{eqn:BasicIneq}
    \bar{\mathcal{L}}(\hat{\vtheta}) - \bar{\mathcal{L}}(\vtheta_o) 
    \le 
    \Tr\left( (\hat{\mSigma} - \mSigma_o)\mSigma_o^{-1}\right)
    - \frac{1}{L}\sum_{l=1}^L\vy_l^H(\hat{\mSigma} - \mSigma_o)\vy_l.
\end{equation}

The inequality above relates the population likelihood gap to the empirical fluctuation term and serves as the starting point of the proof. 
The main argument proceeds in three steps.
\begin{itemize}
    \item 
    Control the right-hand side of \eqref{eqn:BasicIneq} uniformly over the
    parameter class using discretization and concentration, in terms of $\| \mSigma_o^{-1/2} \bigl(\mSigma(\vtheta)-\mSigma_o\bigr)
    \mSigma_o^{-1/2} \|_{\mathrm{HS}}^2.$

    \item Lower bound the left-hand side of \eqref{eqn:BasicIneq} in terms of the same quantity.

    \item  Relate the resulting matrix bound to the difference between the corresponding covariance matrices, and then lower bound this difference
    in terms of the joint errors.
\end{itemize}

We first provide the detailed proof of these three steps under the conditions $n\geq 4m$ or $n<4m,\ \sigma_z^2\geq m$. The remaining regime \(n<4m\) and \(0<\sigma_z^2<m\) is handled at the end of the proof using the monotonicity in the noise variance.



\subsection{Main argument for the upper bound}
\label{sec:proof_upper_bound}

Throughout this subsection, we assume that either \(n\geq 4m\) or \(\sigma_z^2\geq m\). We first introduce the spectral quantities and covering nets used in the proof. Define
\begin{align}\label{eq:def-Emax-norm2A}
E_{\max}:=\lambda_{\max}(\mA\mA^H),
\end{align}
and define
\begin{align}\label{eq:def-Lambda-of-Sigma}
\Lambda_{\min}
:=
\inf_{\vtheta\in\Theta}
\lambda_{\min}\!\left(\mSigma(\vtheta)^{-1}\right),
\qquad
\Lambda_{\max}
:=
\sup_{\vtheta\in\Theta}
\lambda_{\max}\!\left(\mSigma(\vtheta)^{-1}\right).
\end{align}
We further define
\begin{equation}
C_{\sigma,m,n}
:=
\max\left\{
\frac{\Lambda_{\max}}{\Lambda_{\min}},
\frac{\Lambda_{\max}E_{\max}}{\Lambda_{\min}^2}
\right\}.
\label{eqn:def-C-sigma-m-n}
\end{equation}

Fix \(\delta_x,\delta_\phi\in(0,1]\). Let \(\mathcal N_x\subset\mathcal C_x\) and \(\mathcal N_\phi\subset\mathcal C_\phi\) be deterministic \(\delta_x\)- and \(\delta_\phi\)-nets, respectively, and define $\mathcal N:=\mathcal N_x\times\mathcal N_\phi$. The nets are independent of $\mA$ and the observations. Then by the polynomial-complexity assumptions in Section~\ref{sec:model and assump}, we have
\begin{align} \label{eq:number-of-logN}
\log|\mathcal N|
\leq
k_x\log\left(\frac{a_xn^{b_x}}{\delta_x}\right)
+ k_\phi\log\left(\frac{a_\phi m^{b_\phi}}{\delta_\phi}\right).
\end{align}

\paragraph{Upper bounding the right-hand side of \eqref{eqn:BasicIneq}.}
For each \(\vtheta\in\Theta\), define
\begin{align}\label{eq:def-G}
\mathbb{G}(\vtheta)
:=
\Tr\left((\mSigma(\vtheta)-\mSigma_o)\mSigma_o^{-1}\right)
-
\frac{1}{L}\sum_{\ell=1}^L
\vy_\ell^H(\mSigma(\vtheta)-\mSigma_o)\vy_\ell .
\end{align}
Then the inequality in \eqref{eqn:BasicIneq} can be written as
\[
\bar{\mathcal L}(\hat\vtheta)-\bar{\mathcal L}(\vtheta_o)
\leq \mathbb G(\hat\vtheta).
\]
For any fixed \(\vtheta\), the quantity \(\mathbb G(\vtheta)\) is centered conditional on \(\mA\), and can be controlled by concentration inequalities.
However, since the MLE \(\hat\vtheta\) depends on the observations, a pointwise concentration bound cannot be applied directly to \(\mathbb G(\hat\vtheta)\). We therefore use the covering net \(\mathcal N\) introduced above.
Let \(\tilde\vtheta=(\tilde\vx,\tilde\vphi)\in\mathcal N\) be a nearest net point to \(\hat\vtheta=(\hat\vx,\hat\vphi)\), so that
\[
\|\hat\vx-\tilde\vx\|_2\leq\delta_x,
\qquad
\|\hat\vphi-\tilde\vphi\|_2\leq\delta_\phi.
\]
We then decompose
\begin{equation}
|\mathbb G(\hat\vtheta)|
\leq
|\mathbb G(\tilde\vtheta)|
+
|\mathbb G(\hat\vtheta)-\mathbb G(\tilde\vtheta)|.
\label{eq:G_decomposition}
\end{equation}

We first consider the first term. For a fixed \(\tilde{\vtheta}\in\mathcal N\) and conditional on \(\mA\), \(\mathbb G(\tilde{\vtheta})\) is an average of independent centered quadratic forms. Applying a concentration inequality at each point of \(\mathcal N\) and then taking a union bound over the net gives the following result.


\begin{lemma}
\label{lem:quartic-y-upbound-1}
Conditioned on \(\mA\), define
\begin{subequations}
\begin{align}
\label{eqn:def_R}
\kappa_{\delta}
&:= k_x\log\left(\frac{a_xn^{b_x}}{\delta_x}\right)
+
k_\phi\log\left(\frac{a_\phi m^{b_\phi}}{\delta_\phi}\right), \qquad
\mathcal R_{\delta}
:= \sqrt{\frac{c_1\kappa_{\delta}}{L}},
\end{align}
where \(c_1>0\) is a sufficiently large constant. Then there exists a constant \(c>0\) such that, with probability at least
\(1-2\exp(-c\kappa_{\delta})\),
\begin{align}
|\mathbb G(\tilde{\vtheta})|
\leq\;&
\mathcal R_{\delta}
\left\|
\mSigma_o^{-1/2}\bigl(\mSigma(\tilde{\vtheta})-\mSigma_o\bigr)\mSigma_o^{-1/2}
\right\|_{\mathrm{HS}}
+
C_{\sigma,m,n} \mathcal R_{\delta}^2
\label{eqn:net-fluctuation}
\end{align}
\end{subequations}
holds simultaneously for all \(\tilde{\vtheta}\in\mathcal N\).
\end{lemma}
The proof is given in Appendix~\ref{PF:quartic-y-upbound-1}.
The second term in \eqref{eq:G_decomposition} is the error induced by replacing \(\hat{\vtheta}\) with net point \(\tilde{\vtheta}\), which is controlled by the next lemma.
\begin{lemma}
\label{lem:quartic-y-upbound-2}
The proof is given in Appendix~\ref{PF:quartic-y-upbound-2}. Let $C_{\sigma,m,n}$ be as defined in \eqref{eqn:def-C-sigma-m-n}.
There exist an absolute constant $c_2>0$ and a constant $C_2>0$, depending only on $x_{\max}$, such that, conditional on $\mA$, with probability at least $1-2\exp(-c_2mL)$,
\begin{equation}\label{eq:lem:quartic-y-upbound-2}
|\mathbb{G}(\hat{\vtheta}) - \mathbb{G}(\tilde{\vtheta})| 
\le C_2 C_{\sigma,m,n} m(\delta_{x} + \delta_{\phi}) .
\end{equation}
\end{lemma}
Combining Lemmas~\ref{lem:quartic-y-upbound-1} and \ref{lem:quartic-y-upbound-2} with \eqref{eq:G_decomposition}, we obtain
\begin{equation}\label{eqn:upperbound of Ghat}
|\mathbb G(\hat{\vtheta})|
\le
R_\delta
\left\|
\mSigma_o^{-1/2}\bigl(\mSigma(\tilde{\vtheta})-\mSigma_o\bigr)\mSigma_o^{-1/2}
\right\|_{\mathrm{HS}}
+ C_{\sigma,m,n}R_\delta^2
+ C_2 C_{\sigma,m,n}m(\delta_x+\delta_\phi).
\end{equation}
Conditional on $\mA$, this inequality holds with probability at least $1-2\exp(-c\kappa_\delta)-2\exp(-c_2mL)$.

\paragraph{Lower bounding the left-hand side of \eqref{eqn:BasicIneq}.}
Write $\tilde{\mSigma}:=\mSigma(\tilde{\vtheta})$. The following lemma bounds the population likelihood gap in \eqref{eqn:BasicIneq} and the proof is given in Appendix~\ref{PF:lowBd-fbar-trace}.
\begin{lemma}
\label{lem:population-gap-lower-bound}
Let \(C_{\sigma,m,n}\) be defined as in
\eqref{eqn:def-C-sigma-m-n}.
There exists a constant \(C>0\), depending only on \(x_{\max}\), such that
\begin{align}
\overline{\mathcal L}(\hat{\vtheta})
-
\overline{\mathcal L}(\vtheta_o)
\geq\;&
\frac{1}{2C_{\sigma,m,n}}
\left\|
\mSigma_o^{-1/2}\bigl(\tilde{\mSigma}-\mSigma_o\bigr)\mSigma_o^{-1/2}
\right\|_{\mathrm{HS}}^2
- C C_{\sigma,m,n}m(\delta_x+\delta_\phi).
\label{eqn:population-gap-lower-bound}
\end{align}
\end{lemma}
Combining this lower bound with the upper bound in \eqref{eqn:upperbound of Ghat} and the basic inequality \eqref{eqn:BasicIneq}, we obtain, with high probability,
\begin{align}\label{eq:combine-upbound-lowbound}
\frac{1}{2C_{\sigma,m,n}} \|\mM(\tilde{\vtheta})\|_{HS}^2
\le \calR_\delta\|\mM(\tilde{\vtheta})\|_{HS}
+ C_{\sigma,m,n}\calR_\delta^2
+ C_3C_{\sigma,m,n}m(\delta_x+\delta_\phi),
\end{align}
where
\[
\mM(\tilde{\vtheta})
:=
\mSigma_o^{-1/2}(\tilde{\mSigma}-\mSigma_o)\mSigma_o^{-1/2},
\]
with $C_3>0$ depending only on $x_{\max}$. Set
\begin{align}\label{eqn:def-zabc}
z = \|\mM(\tilde{\vtheta})\|_{HS} ,\quad 
a = \frac{1}{2C_{\sigma,m,n}} , \quad 
b = \calR_{\delta} , \quad 
c = C_{\sigma,m,n} \calR^2_{\delta} + C_3 C_{\sigma,m,n} m (\delta_x + \delta_\phi), 
\end{align}
then the inequality \eqref{eq:combine-upbound-lowbound} reduces to $az^2 - bz - c \le 0$, so $\frac{b - \sqrt{b^2+4ac}}{2a} \le z \le \frac{b+\sqrt{b^2+4ac}}{2a}$.
Hence we have
\[ z^2 \le 2\left( \frac{b^2}{4a^2} + \frac{b^2 + 4ac}{ 4a^2} \right) = \frac{b^2}{a^2}+\frac{2c}{a} .\]
Substituting the definitions in \eqref{eqn:def-zabc} and simplifying gives
\begin{equation}
\|\mM(\tilde{\vtheta})\|_{HS}^2
\le
8C_{\sigma,m,n}^2\calR_\delta^2 + 4C_3C_{\sigma,m,n}^2m(\delta_x+\delta_\phi).
\label{eqn:trace-term-up-bound}
\end{equation}


\paragraph{Lower bounding $\|\mM(\tilde{\vtheta})\|_{HS}^2$.} 
We next relate $\|\mM(\tilde{\vtheta})\|_{HS}^2$ to the joint estimation error.
By Lemma~\ref{lem:low_and_up_bound_Tr}, we have the following deterministic lower bound,
\begin{equation} \label{eq:lemma_trace_lbd_result}
\left\|
\mSigma_o^{-1/2}\bigl(\tilde{\mSigma}-\mSigma_o\bigr)\mSigma_o^{-1/2}
\right\|_{\mathrm{HS}}^2
\ge \frac{1}{\Lambda_{\max}^2} 
\left\|\tilde{\mPhi} \mA|\tilde{\mX}|^2\mA^H\tilde{\mPhi}^H - \mPhi_o \mA|\mX_o|^2\mA^H\mPhi^H_o \right\|_{\text{HS}}^2.
\end{equation}
The following lemma provides a lower bound for the right-hand side in terms of the joint estimation error and the proof is given in Appendix~\ref{PF:lowerbound-HSnorm}
\begin{lemma}\label{lem:lowerbound-HSnorm}
Let $\kappa_\delta$ be as defined in \eqref{eqn:def_R}. There exist constants $c,c_4,C>0$, depending only on $x_{\min}$ and $x_{\max}$, such that, with probability at least $1-Ce^{-\kappa_\delta}-Ce^{-c(m+n)}$,
\begin{equation}
\left\| \tilde{\mPhi}\mA\tilde{\mX}^2\mA^H\tilde{\mPhi}^H
-
\mPhi_o\mA\mX_o^2\mA^H\mPhi_o^H \right\|_{HS}^2
\ge c_4\bigl(m^2\|\tilde\vx-\vx_o\|_2^2
+
nm\bigl(m-|\Tr(\tilde{\mPhi}^H\mPhi_o)|\bigr)\bigr)
-C(m+n)\kappa_\delta
\label{eq:lowerbound-HSnorm}
\end{equation}
holds simultaneously for all $\tilde\vtheta\in\mathcal N$.
\end{lemma}
Combining \eqref{eq:lowerbound-HSnorm}, \eqref{eq:lemma_trace_lbd_result} and \eqref{eqn:trace-term-up-bound}, we obtain, with probability at least $1-C(e^{-c\kappa_\delta}+e^{-cmL}+e^{-c(m+n)})$,
\begin{align}\label{eq:Etildevtheta-upbound}
m^2\|\tilde\vx-\vx_o\|_2^2
+
nm\bigl(m-|\Tr(\tilde{\mPhi}^H\mPhi_o)|\bigr)
\le
C_5(m+n)\kappa_\delta
+ C_6 \Lambda_{\max}^2 C_{\sigma,m,n}^2
\left( \calR_\delta^2+ m(\delta_x+\delta_\phi)\right),
\end{align}
where $C_5,C_6>0$ depend only on $x_{\min}$ and $x_{\max}$.

\paragraph{Combining the results.}
The previous inequality still depends on the spectrum of $\mA$. By Lemma~\ref{lem:svd-gaussian}, with probability at least $1-Ce^{-c(m+n)}$ for some absolute constants $c,C>0$, the spectral quantities defined in \eqref{eq:def-Lambda-of-Sigma} satisfy
\[
E_{\max}\le C(m+n),
\qquad
\Lambda_{\max}
\le \sigma_z^2+C x_{\max}^2(m+n).
\]
On the same event, when $n\ge4m$, we also have $\Lambda_{\min}\ge\sigma_z^2+c x_{\min}^2n$. For $n<4m$, we use the deterministic bound $\Lambda_{\min}\ge\sigma_z^2$ together with $\sigma_z^2\ge m$. Thus, in both regimes considered in this subsection, we have
\[
\Lambda_{\max}^2
\le C_7\max\{\sigma_z^4,m^2,n^2\},
\qquad
C_{\sigma,m,n}\le C_8,
\]
where $C_7,C_8>0$ depend only on $x_{\min}$ and $x_{\max}$. Substituting these bounds and $\calR_\delta^2=c_1\kappa_\delta/L$ into \eqref{eq:Etildevtheta-upbound} and dividing both sides by $m^2n$, we obtain
\begin{equation}
\begin{aligned}
&\frac{1}{n}\|\tilde{\vx}-\vx_o\|_2^2
+\frac{1}{m}
\bigl(m-|\Tr(\tilde{\mPhi}^H\mPhi_o)|\bigr)
\le
C_{11}\left[ \frac{(m+n)\kappa_\delta}{m^2n}
+
\frac{\max\{\sigma_z^4,m^2,n^2\}}{m^2n}
\left(\frac{\kappa_\delta}{L}+m(\delta_x+\delta_\phi)\right)
\right].
\end{aligned}
\label{eq:joint-loss-tilde-theta}
\end{equation}
We now pass from the net point $\tilde{\vtheta}$ to the estimator $\hat{\vtheta}$.  By the covering property and the triangle inequality,
\[
\frac{1}{n}\|\hat{\vx}-\vx_o\|_2^2
\le
\frac{2}{n}\|\tilde{\vx}-\vx_o\|_2^2+\frac{2\delta_x^2}{n}, 
\qquad
\frac{1}{m}\bigl(m-|\Tr(\hat{\mPhi}^H\mPhi_o)|\bigr)
\le
\frac{2}{m}\bigl(m-|\Tr(\tilde{\mPhi}^H\mPhi_o)|\bigr)+\frac{\delta_\phi^2}{m}.
\]
We choose $\delta_x=n^{-8}$ and $\delta_\phi=m^{-8}$. Under the assumptions $mL\le n^4\log n$ and $L\le m^4\log m$, we have
\[
mL(\delta_x+\delta_\phi)
\le \frac{\log n}{n^4} +\frac{\log m}{m^3} \lesssim \kappa_\delta.
\]
Moreover, $\delta_x^2/n+\delta_\phi^2/m= \frac{1}{n^{17}}+\frac{1}{m^{17}}$ is much smaller than $\frac{(m+n)\kappa_\delta}{m^2n}$. Combining these bounds with
\eqref{eq:joint-loss-tilde-theta} yields
\begin{align}
&\frac{1}{n}\|\hat{\vx}-\vx_o\|_2^2
+\frac{1}{m}\bigl(m-|\Tr(\hat{\mPhi}^H\mPhi_o)|\bigr)
\le
C_7\kappa_\delta
\left[\frac{\max\{\sigma_z^4,m^2,n^2\}}{m^2nL}+\frac{m+n}{m^2n}\right].
\label{eq:joint-loss-final-kappa}
\end{align}
Here, the covering-number bounds give
\begin{align*}
\kappa_\delta
=k_x\log(a_xn^{b_x+8})+k_\phi\log(a_\phi m^{b_\phi+8}) 
\le C_8(k_x\log n+k_\phi\log m).
\end{align*}
The bound holds with probability at least $1-C\left( e^{-c\kappa_\delta} +e^{-cmL} +e^{-c(m+n)} \right)$.

\subsection{Proof of Theorem~\ref{thm:joint_upper_bound}}

We first consider $n\ge4m$, and $n<4m$ with $\sigma_z^2\ge m$. Since the estimator takes values in $\Theta$, its joint loss is bounded by a constant $C_1= (x_{\max}-x_{\min})^2+2$. The conclusion of Section~\ref{sec:proof_upper_bound} therefore gives
\begin{align}
\E\!\left[d_{\rm joint}^2
\bigl((\hat\vx,\hat\mPhi),(\vx_o,\mPhi_o)\bigr)\right]
\le{}&C\kappa_\delta
\left[\frac{\max\{\sigma_z^4,m^2,n^2\}}{m^2nL}
+\frac{m+n}{m^2n}\right]
+C_1\bigl(e^{-\kappa_\delta}+e^{-cmL}+e^{-c(m+n)}\bigr),
\label{eq:joint-risk-expectation}
\end{align}
where the expectation includes both $\mA$ and the observations. To control the exceptional-event contribution, note that the conditions $mL\le n^4\log n$ and $L\le m^4\log m$ imply $mL \le\min\{m,n\}^{6}$. 
Also, the choices $\delta_x=n^{-8}$ and $\delta_\phi=m^{-8}$ give
\[
\kappa_\delta \ge8(k_x\log n+k_\phi\log m) \ge8\log\min\{m,n\}.
\]
Using the inequalities $e^{-t}\le t^{-1}$ and $e^{-t}\le6!\,t^{-6}$ we obtain
\begin{align*}
e^{-\kappa_\delta}+e^{-cmL}+e^{-c(m+n)}
\le \min\{m,n\}^{-8} +\frac{C}{mL} +\frac{C}{(m+n)^6}
\le \frac{C}{mL}
\le C\frac{\max\{\sigma_z^4,m^2,n^2\}}{m^2nL},
\end{align*}
Substituting this estimate into \eqref{eq:joint-risk-expectation} and using $\kappa_\delta\le C(k_x\log n+k_\phi\log m)$ together with the boundedness of the joint loss, we obtain
\begin{align}
&\sup_{\vtheta_o\in\Theta}
\E\!\left[d_{\rm joint}^2
\bigl((\hat\vx,\hat\mPhi),(\vx_o,\mPhi_o)\bigr)\right]
\le C\min\!\left\{
(k_x\log n+k_\phi\log m)
\left[\frac{\max\{\sigma_z^4,m^2,n^2\}}{m^2nL}
+\frac{m+n}{m^2n}\right],1\right\}.
\label{eq:joint-risk-direct-regimes}
\end{align}
Since the minimax risk is at most the risk of this estimator, the same
bound holds for
$R_{\rm joint}(\mathcal C_x,\mathcal C_\phi;m,n,L,\sigma_z)$.

To cover the remaining noise regime, we use the following comparison.
\begin{lemma}\label{lem:noise-comparison}
Let $\Theta=\mathcal C_x\times\mathcal C_\phi$ be the compact parameter space. For $s>0$, let $\E_{\vtheta,s}$ denote expectation under parameter $\vtheta$ and additive-noise covariance $s\mI_m$, including the randomness of $\mA$.
Fix $0<s_1\le s_2$. For every measurable $\Theta$-valued estimator $\widehat\vtheta_2$ at noise variance $s_2$, there exists a measurable $\Theta$-valued estimator $\widehat\vtheta_1$ at noise variance $s_1$ such that
\begin{equation}
\E_{\vtheta_o,s_1}
\bigl[d_{\rm joint}^2(\widehat\vtheta_1,\vtheta_o)\bigr]
\le4\,\E_{\vtheta_o,s_2}
\bigl[d_{\rm joint}^2(\widehat\vtheta_2,\vtheta_o)\bigr],
\qquad \vtheta_o\in\Theta.
\label{eq:noise-comparison-risk}
\end{equation}
\end{lemma}
The proof is given in Appendix~\ref{app:noise-comparison}.
Now suppose $n<4m$ and $0<\sigma_z^2<m$. Let $\widehat\vtheta_m$ be the measurable constrained MLE for additive-noise variance $m$. Applying Lemma~\ref{lem:noise-comparison} with $s_1=\sigma_z^2$ and $s_2=m$ transfers its uniform risk bound to a deterministic $\Theta$-valued estimator at variance $\sigma_z^2$. By the definition of the minimax risk and
\eqref{eq:joint-risk-direct-regimes},
\begin{align*}
R_{\rm joint}(\mathcal C_x,\mathcal C_\phi;m,n,L,\sigma_z)
&\le4\sup_{\vtheta_o\in\Theta}
\E_{\vtheta_o,m}\!\left[
 d_{\rm joint}^2(\widehat\vtheta_m,\vtheta_o)\right]\\
&\le C\min\!\left\{
(k_x\log n+k_\phi\log m)
\left[\frac{\max\{m^2,n^2\}}{m^2nL}
+\frac{m+n}{m^2n}\right],1\right\}.
\end{align*}
Since $\sigma_z^2<m$ implies
$\max\{\sigma_z^4,m^2,n^2\}=\max\{m^2,n^2\}$,
this yields the same upper bound in the remaining regime
and completes the proof.

\begin{remark}
Xing et al.~\cite{xing2025minimax} compare noise levels by adding independent Gaussian noise and then taking the conditional mean of the resulting estimator. For Euclidean squared loss, Jensen's inequality ensures that this averaging does not increase the risk.
We use a similar argument, but cannot directly apply the averaging step since averaging phase matrices may violate the unit-modulus constraint. Instead, Lemma~\ref{lem:noise-comparison} chooses an estimate in $\Theta$ that minimizes the expected joint loss between it and the estimate obtained after adding noise. This yields a risk comparison up to a factor of four, which is absorbed into the constant.
\end{remark}

\section{Proof of Lower Bound}
\label{sec:proof lower bound}

We prove the lower bound in Theorems~\ref{thm:joint-lb-n-ge-m} and~\ref{thm:joint-lb-m-ge-n} by reducing joint estimation to signal- and phase-oracle problems.
Subsection~\ref{sec: lower bound main sec} introduces the oracle risks, states their lower bounds, and combines these bounds to prove the two theorems. Subsection~\ref{sec: lower bound phase} provides the proofs of the phase-oracle lower bounds.

\subsection{Oracle bounds and proofs of Theorems~\ref{thm:joint-lb-n-ge-m} and~\ref{thm:joint-lb-m-ge-n}} \label{sec: lower bound main sec}

Recall that the joint loss is
\[
d_{\mathrm{joint}}^2\bigl((\hat{\vx},\hat{\mPhi}),(\vx_o,\mPhi_o)\bigr)
=
\frac{1}{n}\|\hat{\vx}-\vx_o\|_2^2
+
\frac{1}{m}d_\Phi^2(\hat{\mPhi},\mPhi_o).
\]
Both loss components are nonnegative. Moreover, fixing one parameter at an admissible value and revealing it to the estimator can only make estimation of the other parameter easier. Consequently, a lower bound for either oracle problem yields a lower bound for the corresponding joint estimation problem.
We introduce the following two oracle minimax risks.
\begin{itemize}
    \item \textbf{Signal-oracle risk.}
    Fix an admissible phase matrix $\mPhi_o$ and reveal it
    to the estimator, so that only the signal amplitude remains
    unknown. Define
    \[
    R_{X}(\mathcal C_x;m,n,L,\sigma_z)
    :=
    \inf_{\widehat{\vx}}\sup_{\vx_o\in\mathcal C_x}
    \mathbb E_{\vx_o,\mPhi_o}
    \left[ \frac{1}{n}\|\widehat{\vx}-\vx_o\|_2^2 \right],
    \]
    where the infimum is over all measurable estimators based on $\{\vy_\ell\}_{\ell=1}^L$, $\mA$ and $\mPhi_o$. 
    Applying $\mPhi_o^H$ to each observation removes the phase distortion without changing the additive-noise distribution. Thus this risk does not depend on the fixed phase matrix $\mPhi_o$, which is omitted from its notation.

    \item \textbf{Phase-oracle risk.}
    Fix a signal amplitude $\vx_o\in[x_{\min},x_{\max}]^n$ and reveal it to the estimator, so that only the phase distortion remains unknown. Define
    \[
    R_\Phi(\mathcal C_\phi;m,n,L,\sigma_z)
    :=
    \inf_{\widehat{\mPhi}}\sup_{\vphi_o\in\mathcal C_\phi}
    \mathbb E_{\vx_o,\vphi_o}
    \left[ \frac{1}{m} d_\Phi^2 \bigl(\widehat{\mPhi},\mPhi(\vphi_o)\bigr) \right],
    \]
    where the infimum is over all measurable estimators taking values in the set of diagonal unitary matrices and based on $\{\vy_\ell\}_{\ell=1}^L$, $\mA$ and $\vx_o$. 
    For notation convenience, dependence of this risk on the fixed signal $\vx_o$ is suppressed in the notation.
\end{itemize}

\subsubsection{Oracle lower bounds}

We first give lower bounds for the risks corresponding to the two oracle models.
\paragraph{Signal-oracle Model.}
Since $\mPhi_o$ is known, multiplying each observation by $\mPhi_o^H$ removes the phase distortion:
\begin{equation}
\label{eq:model-oracle-signal}
\widetilde{\vy}_\ell
:=\mPhi_o^H\vy_\ell
=\mA\mX_o\vw_\ell+\widetilde{\vz}_\ell,
\qquad \ell=1,\ldots,L,
\end{equation}
where $\widetilde{\vz}_\ell:=\mPhi_o^H\vz_\ell$.
The transformed additive noise remains $\mathcal{CN}(\mathbf0,\sigma_z^2\mI_m)$ because $\mPhi_o$ is unitary and the original noise is isotropic complex Gaussian.
Thus \eqref{eq:model-oracle-signal} is a complex-valued version of the fixed-sensing-matrix model studied in \cite{xing2025minimax}.
We state the complex-valued version of \cite[Theorem~2.15]{xing2025minimax} below.
\begin{lemma}\label{lem:signal-oracle-lb}
Consider the model \eqref{eq:model-oracle-signal}, where the known matrix $\mA\in\mathbb C^{m\times n}$ has independent $\mathcal{CN}(0,1)$ entries. Assume that $\vw_\ell\sim\mathcal{CN}(\mathbf0,\mI_n)$ and $\widetilde{\vz}_\ell\sim\mathcal{CN}(\mathbf0,\sigma_z^2\mI_m)$, with all these noise vectors mutually independent and independent of $\mA$. Fix $\varepsilon\in(0,1/2)$. Suppose that (i) $a_x\ge x_{\max}-x_{\min}$, $b_x\ge1$, $\log m=\Theta(\log n)$, and $\log L=O(\log n)$; (ii) $k_x\le n^{1-2\varepsilon}$; and (iii) $\max\{\sigma_z^4,m^2,n^2\}k_x\log n \le m^2n^{1-\varepsilon}L.$
Then
\[
\sup_{\mathcal C_x\in\mathcal F_{a_x,b_x,k_x,n}}
R_X(\mathcal C_x;m,n,L,\sigma_z)
=
\Omega_{\varepsilon,x_{\min},x_{\max}}
\left( \frac{\max\{\sigma_z^4,m^2,n^2\}k_x\log n}{m^2nL} \right).
\]
\end{lemma}
The proof is given in Appendix~\ref{pf:lem:signal-oracle-lb}, where we show that the real-valued lower bound in \cite[Theorem~2.15]{xing2025minimax} extends to the complex model up to constant factors. For the phase-oracle problem below, we follow the same Fano strategy used in \cite{xing2025minimax}, with some different arguments on the KL bounds which depends on the covariance structure.

\paragraph{Phase-oracle Model.}
For phase estimation, the lower bound depends on the spectral of the covariance matrix $\mSigma^{-1}$, which have different behavior depends on the model dimension $n$ and $m$. Hence we gives the results in two regimes $n\ge(1+\tau)m$ and $m\ge n$. 
\begin{lemma}[$n\ge(1+\tau)m$]
\label{lem:phase-only lower bound n ge m}
Fix $\tau>0$, $\varepsilon\in(0,1/2)$, and a known signal
$\vx_o\in[x_{\min},x_{\max}]^n$. Suppose that
(i) $a_\phi\ge a_{\phi,0}$ and $b_\phi\ge b_{\phi,0}$, where
$a_{\phi,0},b_{\phi,0}>0$ are the absolute constants specified in
Subsection~2.2;
(ii) $n\ge(1+\tau)m$, $\log n=\Theta(\log m)$, and
$\log L=O(\log m)$;
(iii) $k_\phi\le m^{1-2\varepsilon}$; and
(iv)
\[
\max\{\sigma_z^4,n^2\}k_\phi\log m
\le m^{2-\varepsilon}nL.
\]
Then we have,
\[
\sup_{\mathcal C_\phi\in\mathcal F_{a_\phi,b_\phi,k_\phi,m}}
R_\Phi(\mathcal C_\phi;m,n,L,\sigma_z)
=
\Omega_{\tau,\varepsilon,x_{\min},x_{\max}}
\left(
\frac{\max\{\sigma_z^4,m^2,n^2\}k_\phi\log m}{m^2nL}
\right).
\]
\end{lemma}

\begin{lemma}[$m\ge n$]
\label{lem:phase-only lower bound n le m}
Fix $\varepsilon\in(0,1/2)$ and a known signal
$\vx_o\in[x_{\min},x_{\max}]^n$. Suppose that
(i) $a_\phi\ge a_{\phi,0}$ and $b_\phi\ge b_{\phi,0}$;
(ii) $m\ge n$, $\log n=\Theta(\log m)$, and $\log L=O(\log m)$;
(iii) $k_\phi\le m^{1-2\varepsilon}$; and
(iv)
\[
\max\{\sigma_z^4,m^2\}k_\phi\log m \le m^{2-\varepsilon}nL.
\]
For $\sigma_z^2>0$, uniformly over the fixed known signal,
\[
\sup_{\mathcal C_\phi\in\mathcal F_{a_\phi,b_\phi,k_\phi,m}}
R_\Phi(\mathcal C_\phi;m,n,L,\sigma_z)
=
\Omega_{\varepsilon,x_{\min},x_{\max}}
\left(
\frac{\max\{\sigma_z^4,m^2,n^2\}k_\phi\log m}{m^2nL}
\cdot\frac{\sigma_z^2}{\max\{\sigma_z^2,m,n\}}
\right).
\]
At $\sigma_z^2=0$, the corresponding inequality has a zero
right-hand side and follows from nonnegativity of the risk.
\end{lemma}


The proofs of both phase-oracle bounds are given in Subsection~\ref{sec: lower bound phase}. The second lower bound contains an additional factor depending on the additive-noise level and the model dimensions. When $\sigma_z^2=\Theta(m)$, this factor is constant so the two lower-bound formulas have the same form up to constant factors. When $\sigma_z^2\ll m$, however, it equals $\sigma_z^2/m$ and makes the second lower bound much smaller.

This factor can be understood in terms of the condition number of the oracle covariance matrix $\mSigma^{-1}$. When $n\ge(1+\tau)m$, the eigenvalues of $\mA\mA^H$ are of order $n$ with high probability, so the covariance condition number is bounded by a constant and is absorbed into the constants.
When $m>n$ and $\sigma_z^2>0$, $\mA\mA^H$ is rank deficient: the smallest eigenvalue of the oracle covariance equals $\sigma_z^2$, while its largest eigenvalue is of order $\max\{\sigma_z^2,m\}$ with high probability. Thus the additional factor is of the same order as the inverse condition number. At $m=n$, this factor gives a lower bound on the inverse condition number up to a constant.

\subsubsection{Proof of Theorems~\ref{thm:joint-lb-n-ge-m}
and~\ref{thm:joint-lb-m-ge-n}}
\label{sec:proof_lb_thm}

Fix $\vx_o\in[x_{\min},x_{\max}]^n$ and $\vphi_o\in\{\vphi\in[0,2\pi)^m:\phi_1=0\}$. 
Since singleton classes have covering number one, both $\{\vx_o\}$ and $\{\vphi_o\}$ satisfy the polynomial-complexity conditions and hence belong to $\mathcal F_{a_x,b_x,k_x,n}$ and $\mathcal F^\phi_{a_\phi,b_\phi,k_\phi,m}$ respectively for all sufficiently large $n$ and $m$. 
Consequently, the supremum over admissible joint classes includes both $(\mathcal C_x,\{\vphi_o\})$ and $(\{\vx_o\},\mathcal C_\phi)$. Hence, by the oracle reductions and nonnegativity of the two losses, for any $\mathcal C_x\in\mathcal F_{a_x,b_x,k_x,n}$ and $\mathcal C_\phi\in\mathcal F_{a_\phi,b_\phi,k_\phi,m}$, we obtain.
\begin{align*}
R_{\mathrm{joint}}(\mathcal C_x,\{\vphi_o\};m,n,L,\sigma_z)
&\ge R_X(\mathcal C_x;m,n,L,\sigma_z),\\
R_{\mathrm{joint}}(\{\vx_o\},\mathcal C_\phi;m,n,L,\sigma_z)
&\ge R_\Phi(\mathcal C_\phi;m,n,L,\sigma_z).
\end{align*}
Taking suprema yields
\begin{equation}
\label{eq:joint-oracle-reduction}
\begin{aligned}
&\sup_{\substack{
\mathcal C_x\in\mathcal F_{a_x,b_x,k_x,n}\\
\mathcal C_\phi\in\mathcal F^\phi_{a_\phi,b_\phi,k_\phi,m}
}}
R_{\mathrm{joint}}(\mathcal C_x,\mathcal C_\phi;m,n,L,\sigma_z)
\ge
\max\left\{
\begin{aligned}
&\sup_{\mathcal C_x\in\mathcal F_{a_x,b_x,k_x,n}}
R_X(\mathcal C_x;m,n,L,\sigma_z),\\
&\sup_{\mathcal C_\phi\in\mathcal F^\phi_{a_\phi,b_\phi,k_\phi,m}}
R_\Phi(\mathcal C_\phi;m,n,L,\sigma_z)
\end{aligned}
\right\}.
\end{aligned}
\end{equation}

We first consider the regime $n\ge(1+\tau)m$. Applying Lemmas~\ref{lem:signal-oracle-lb} and~\ref{lem:phase-only lower bound n ge m} in \eqref{eq:joint-oracle-reduction} and $\max\{\sigma_z^4,m^2,n^2\}=\max\{\sigma_z^4,n^2\}$, we obtain
\begin{align*}
\sup_{\substack{
\mathcal C_x\in\mathcal F_{a_x,b_x,k_x,n}\\
\mathcal C_\phi\in\mathcal F_{a_\phi,b_\phi,k_\phi,m}
}}
R_{\mathrm{joint}}(\mathcal C_x,\mathcal C_\phi;m,n,L,\sigma_z)
\ge
c\frac{\max\{\sigma_z^4,n^2\}}{m^2nL}
\max\{k_x\log n,k_\phi\log m\},
\end{align*}
where $c>0$ depends only on
$\tau,\varepsilon,x_{\min},x_{\max}$.
This proves Theorem~\ref{thm:joint-lb-n-ge-m}.

We next consider the regime $m\ge n$. Applying Lemmas~\ref{lem:signal-oracle-lb} and~\ref{lem:phase-only lower bound n le m} in \eqref{eq:joint-oracle-reduction} and $\max\{\sigma_z^2,m,n\}=\max\{\sigma_z^2,m\}$,  we obtain
\begin{align*}
\sup_{\substack{
\mathcal C_x\in\mathcal F_{a_x,b_x,k_x,n}\\
\mathcal C_\phi\in\mathcal F_{a_\phi,b_\phi,k_\phi,m}
}}
R_{\mathrm{joint}}(\mathcal C_x,\mathcal C_\phi;m,n,L,\sigma_z)
\ge
c\frac{\max\{\sigma_z^4,m^2\}}{m^2nL}
\max\left\{
k_x\log n,\,
\frac{\sigma_z^2}{\max\{\sigma_z^2,m\}}
k_\phi\log m
\right\},
\end{align*}
where $c>0$ depends only on $\varepsilon,x_{\min},x_{\max}$.
This proves Theorem~\ref{thm:joint-lb-m-ge-n}.

\subsection{Lower Bound of Phase-Only Model}
\label{sec: lower bound phase}

In this subsection, we prove
Lemmas~\ref{lem:phase-only lower bound n ge m}
and~\ref{lem:phase-only lower bound n le m} for the phase-oracle model.

\paragraph{Oracle phase model.}

Fix the known signal $\vx_o\in[x_{\min},x_{\max}]^n$, and write
\[
\mX_o:=\operatorname{diag}(\vx_o),\qquad
\mS_o:=\mA\mX_o^2\mA^H.
\]
We parameterize the unknown phase by an angle vector $\vphi_o\in\mathcal C_\phi\subset \{\vphi\in[0,2\pi)^m:\phi_1=0\}$ and write
\[
\mPhi_o:=\mPhi(\vphi_o),\qquad
\mPhi(\vphi):=\diag(e^{i\phi_1},\ldots,e^{i\phi_m}).
\]
Whenever $\sigma_z^2\mI_m+\mS_o\succ0$, conditional on $\mA$ the observations are independent and satisfy
\[
\vy_\ell\mid\mA
\sim\mathcal{CN}(\mathbf0,\mSigma(\mPhi_o)^{-1}),
\qquad \ell=1,\ldots,L,
\]
where, for each candidate phase matrix $\mPhi=\mPhi(\vphi)$,
\[
\mSigma(\mPhi)^{-1} :=\sigma_z^2\mI_m+\mPhi\mS_o\mPhi^H.
\]
For each $\vphi\in\mathcal C_\phi$, denote the single-look conditional law by $Q_{\mPhi}:=\mathcal{CN}(\mathbf0,\mSigma(\mPhi)^{-1})$ and the $L$-look conditional law by $P_{\mPhi}:=Q_{\mPhi}^{\otimes L}$. All KL divergences below are conditional on $\mA$.

\paragraph{Proof outline.}
As in the signal-oracle problem consider in \cite{xing2025minimax}, our proof relies on the generalized Fano lemma in Lemma~\ref{lem:generalized_fano}. The proof proceed in three steps.
\begin{itemize}
\item \textbf{Step 1: Construction of a phase class and packing.}
We construct a compact piecewise-constant angle class satisfying the polynomial-complexity condition and select a finite packing. We then obtain a separation lower bound $\alpha_r$ for the corresponding phase matrices under the normalized phase pseudo-metric $d:=m^{-1/2}d_\Phi$.

\item \textbf{Step 2: KL upper bound.}
We bound the maximum pairwise KL divergence conditional on $\mA$, denoted by $\beta_r(\mA)$, using different estimates for the regimes $n\ge(1+\tau)m$ and $m\ge n$. These bounds allow us to keep the observation distributions sufficiently close, making the separated phase parameters difficult to distinguish from the data.

\item \textbf{Step 3: Applying Fano's theorem.}
We choose the packing radius so that $\beta_r(\mA)$ is small relative to $\log r$. Fano's inequality and the packing separation then gives a squared-loss lower bound of order $\alpha_r^2$, yielding the two phase-oracle lower bounds.
\end{itemize}

\paragraph{Phase packing and Fano reduction}
We first construct a phase class satisfying the polynomial-complexity condition and a finite packing with large cardinality, separation in the phase metric, and controlled pairwise Frobenius distances.
\begin{lemma}[Construction of the phase class]
\label{lem:phase_class_construction}
Fix $\eta\in(0,1)$, and let $N_{\mathrm{div}}^\phi:=\lfloor k_\phi m^\eta\rfloor\le m-1$. Then there exists a compact angle class
\[
\Gamma_\phi\subset \{\vphi\in[0,\pi/8]^m:\phi_1=0\}
\]
with polynomial complexity of order $k_\phi$: for absolute constants $a_{\phi,0},b_{\phi,0}>0$,
\[
N_\epsilon(\Gamma_\phi)
\le
\left(\frac{a_{\phi,0}m^{b_{\phi,0}}}{\epsilon}
\right)^{k_\phi}, \qquad \epsilon\in(0,1].
\]
Moreover, for every $\delta_r\in(0,\pi/8)$, one can construct a finite subset $U_{\mathrm{sep}}=\{\vphi^{(1)},\ldots,\vphi^{(r)}\}\subset\Gamma_\phi$ with the following properties:
\begin{enumerate}
\item[(i)] $\log r=\Theta(k_\phi\log m)$.
\item[(ii)] For all $i\ne j$,
\[
c\frac{m k_\phi}{N_{\mathrm{div}}^\phi}\delta_r^2
\le d_\Phi^2(\mPhi_i,\mPhi_j)
\le C\frac{m k_\phi}{N_{\mathrm{div}}^\phi}\delta_r^2.
\]
\item[(iii)] For all $i\ne j$,
\[
\|\mPhi_i-\mPhi_j\|_F^2
\le C\frac{m k_\phi}{N_{\mathrm{div}}^\phi}\delta_r^2.
\]
\end{enumerate}
Here $c,C>0$ are absolute constants.
\end{lemma}

The proof is given in Appendix~\ref{pf:lem:phase_class_construction}.

Let $\eta:=\varepsilon/2$ and $N_{\mathrm{div}}^\phi:=\lfloor k_\phi m^\eta\rfloor$ . Since $k_\phi\le m^{1-2\varepsilon}$, we have $N_{\mathrm{div}}^\phi\le m^{1-3\varepsilon/2}$, so Lemma~\ref{lem:phase_class_construction} applies for sufficiently large $m$. For a deterministic radius $\delta_r\in(0,\pi/8)$ be chosen later, let $U_{\mathrm{sep}}=\{\vphi^{(1)},\ldots,\vphi^{(r)}\}\subset\Gamma_\phi$ be the packing constructed in Lemma~\ref{lem:phase_class_construction}, independently of $\mA$. Then $\log r=\Theta(k_\phi\log m)$, and the corresponding phase matrices constructed by the phase in class $\Gamma_\phi$ satisfy
\begin{align}\label{eq:packing-separation-alpha}
d(\mPhi_i,\mPhi_j) :=m^{-1/2}d_\Phi(\mPhi_i,\mPhi_j) \ge\alpha_r,
\qquad
\alpha_r^2:=c\frac{k_\phi}{N_{\mathrm{div}}^\phi}\delta_r^2.
\end{align}
Define the maximum pairwise KL divergence conditional on $\mA$ by
\begin{equation}\label{eq:def-beta-maxKL}
\beta_r(\mA):=
\max_{1\le i\ne j\le r}D_{\rm KL}(P_{\mPhi_i}\|P_{\mPhi_j}).
\end{equation}
We next bound $\beta_r(\mA)$ in the two regimes and choose $\delta_r$ so that $\beta_r(\mA)\le (\log r)/8$ holds with high probability. The minimax lower bounds will then follow from the Fano's theorem in Lemma~\ref{lem:generalized_fano}.

\paragraph{KL upper bound for $n\ge(1+\tau)m$.}
Following a similar approach to \cite{xing2025minimax}, we control the conditional KL divergence by following lemma.
\begin{lemma}
\label{lem:phase-oracle-kl}
Suppose $n\ge(1+\tau)m$ and let $\mA\in\mathbb C^{m\times n}$ have independent $\mathcal{CN}(0,1)$ entries. Fix the known signal $\vx_o\in[x_{\min},x_{\max}]^n$, and set $\mX_o=\operatorname{diag}(\vx_o)$ and $\mS_o=\mA\mX_o^2\mA^H$.
There exists a constant $C>0$, depending only on $x_{\max}$, such that, with probability at least $1-me^{-m}-me^{-n}$, the following bound holds simultaneously for all diagonal unitary phase matrices $\mPhi_i,\mPhi_j$:
\[
D_{\rm KL}(P_{\mPhi_i}\|P_{\mPhi_j})
\le
CL\frac{mn}{ \lambda_{\min}^2(\sigma_z^2\mI_m+\mS_o) }\|\mPhi_i-\mPhi_j\|_F^2.
\]
\end{lemma}


The proof is given in Appendix~\ref{Pf:lem-phase-oracle-kl}.
Since $n\ge(1+\tau)m$, by Lemma~\ref{lem:svd-gaussian} and $\mX_o^2\succeq x_{\min}^2\mI_n$, with probability at least $1-2e^{-c_\tau n}$, we have
\[
\lambda_{\min}(\sigma_z^2\mI_m+\mS_o)
\ge c_{\tau,x_{\min}}\max\{\sigma_z^2,n\}.
\]
Combining this with Lemmas~\ref{lem:phase-oracle-kl} and \ref{lem:phase_class_construction}(iii), we obtain, with probability at least $1-2e^{-c_\tau n}-me^{-m}-me^{-n}$,
\[
\beta_r(\mA)
\le C_{\tau,x_{\min},x_{\max}}
\frac{Lm^2n}{\max\{\sigma_z^4,n^2\}}
\frac{k_\phi}{N_{\mathrm{div}}^\phi}\delta_r^2.
\]
Choose the deterministic radius
\begin{align}
\label{eq:choice-deltar-ngem}
\delta_r^2:=c_0
\frac{\max\{\sigma_z^4,n^2\}}{m^2nL} N_{\mathrm{div}}^\phi\log m,
\end{align}
where $c_0>0$ depends only on $\tau,\varepsilon,x_{\min},x_{\max}$. Then the rate assumption (iv) in Lemma~\ref{lem:phase-only lower bound n ge m} and $N_{\mathrm{div}}^\phi\le k_\phi m^{\varepsilon/2}$ give
\[
\delta_r^2
\le c_0m^{\varepsilon/2}
\frac{\max\{\sigma_z^4,n^2\}k_\phi\log m}{m^2nL}
\le c_0m^{-\varepsilon/2}.
\]
Thus, the chosen radius \(\delta_r\) satisfies the requirement $0<\delta_r<\pi/8$  of Lemma~\ref{lem:phase_class_construction} for sufficiently large \(m\). Since $\log r=\Theta(k_\phi\log m)$, choosing $c_0$ sufficiently small yields
\begin{equation}\label{eq:phase-kl-prob-n-ge-m}
\mathbb P\!\left(\beta_r(\mA)\le\frac{\log r}{8}\right)
\ge 1-2e^{-c_\tau n}-me^{-m}-me^{-n}.
\end{equation}

\paragraph{KL upper bound for $n\le m$.}
In this setting, the covariance matrix may have a small minimum eigenvalue when $\sigma_z^2$ is small. Then the bound in Lemma~\ref{lem:phase-oracle-kl} will be too loose. The following lemma avoids this loss.

\begin{lemma}
\label{lem:refined-kl-tall-phase}
Let $\mA\in\mathbb C^{m\times n}$ have independent $\mathcal{CN}(0,1)$ entries and assume that \(n\le m\). Fix the known signal $\vx_o\in[x_{\min},x_{\max}]^n$, let $\mX_o=\operatorname{diag}(\vx_o)$, and set $\mS_o:=\mA\mX_o^2\mA^H$, and assume $\sigma_z^2>0$.
For two phase matrices $\mPhi_i,\mPhi_j$, let $P_{\mPhi_i}$ and $P_{\mPhi_j}$ denote the corresponding $L$-look phase-oracle laws.
Then there exists a constant $C>0$, depending only on $x_{\max}$, such that, with probability at least $1-2e^{-m/2}-me^{-n}$, the following bound holds simultaneously for all diagonal unitary phase matrices $\mPhi_i,\mPhi_j$:
\[
D_{\rm KL}(P_{\mPhi_i}\|P_{\mPhi_j})
\le CL\frac{mn}{\sigma_z^2\max\{\sigma_z^2,m\}}
\|\mPhi_i-\mPhi_j\|_F^2.
\]
\end{lemma}

The proof is given in Appendix~\ref{pf:lem:refined-kl-tall-phase}.
For $\sigma_z^2=0$, the lower bound in Lemma~\ref{lem:phase-only lower bound n le m} follows immediately from the nonnegativity of the risk, since its right-hand side vanishes. We therefore assume $\sigma_z^2>0$ in what follows. By Lemmas~\ref{lem:refined-kl-tall-phase} and~\ref{lem:phase_class_construction}(iii), with probability at least $1-2e^{-m/2}-me^{-n}$, we have
\[
\beta_r(\mA)
\le C_{x_{\max}}
\frac{Lm^2n}{\sigma_z^2\max\{\sigma_z^2,m\}}
\frac{k_\phi}{N_{\mathrm{div}}^\phi}\delta_r^2.
\]
Choose the deterministic radius
\begin{align}
\label{eq:choice-deltar-mgen}
\delta_r^2:=c_0
\frac{\sigma_z^2\max\{\sigma_z^2,m\}}{m^2nL} N_{\mathrm{div}}^\phi\log m,
\end{align}
where $c_0>0$ depends only on $\varepsilon,x_{\max}$. The rate assumption (iv) in Lemma~\ref{lem:phase-only lower bound n le m} and $N_{\mathrm{div}}^\phi\le k_\phi m^{\varepsilon/2}$ imply
\begin{align}
\delta_r^2
\le c_0m^{\varepsilon/2}
\frac{\sigma_z^2\max\{\sigma_z^2,m\}k_\phi\log m}{m^2nL}
\le c_0\frac{\sigma_z^2}{\max\{\sigma_z^2,m\}}m^{-\varepsilon/2}
\le c_0m^{-\varepsilon/2}.
\end{align}
Thus, the chosen radius \(\delta_r\) satisfies the requirement $0<\delta_r<\pi/8$  of Lemma~\ref{lem:phase_class_construction} for sufficiently large \(m\). Since $\log r=\Theta(k_\phi\log m)$, choosing $c_0$ sufficiently small gives
\begin{equation}\label{eq:phase-kl-prob-n-le-m}
\mathbb P\!\left(\beta_r(\mA)\le\frac{\log r}{8}\right)
\ge 1-2e^{-m/2}-me^{-n}.
\end{equation}
Under the standing assumption $\log n=\Theta(\log m)$, the failure probability in this bound tends to zero.

\paragraph{Applying Fano's inequality.}
In both regimes, the KL bounds and the choices of $\delta_r$ imply
\[
\mathbb P\!\left(\beta_r(\mA)\le\frac{\log r}{8}\right)\ge1-\epsilon_{m,n},
\qquad \epsilon_{m,n}=o(1),
\]
where $\epsilon_{m,n}$ is given by the failure-probability bound in \eqref{eq:phase-kl-prob-n-ge-m} or \eqref{eq:phase-kl-prob-n-le-m}, respectively. Define
\[
\mathcal E
:=\left\{\mA:\beta_r(\mA)\le\frac{\log r}{8}\right\}.
\]
For each $j\in[r]$, let $\widetilde P_j$ be the joint law of $\mA$ and the $L$ observations under $\mPhi_j$, conditional on $\mathcal E$. Since $\mathcal E$ depends only on $\mA$ (with $\vx_o$ fixed), and the distribution of $\mA$ is the same under all hypotheses, the chain rule for KL divergence yields
\begin{align}
D_{\rm KL}(\widetilde P_i\|\widetilde P_j)
&=\mathbb E_{\mA\mid\mathcal E}
\left[D_{\rm KL}(P_{\mPhi_i}\|P_{\mPhi_j})\right]
\le\frac{\log r}{8},\qquad i\ne j.
\end{align}
By Lemma~\ref{lem:generalized_fano} and the packing separation $\alpha_r$, for every estimator $\widehat\mPhi$,
\[
\max_{1\le j\le r}
\mathbb E_{\widetilde P_j}
\left[d(\widehat\mPhi,\mPhi_j)\right]
\ge\frac{\alpha_r}{2}
\left(1-\frac{\log r/8+\log2}{\log r}\right).
\]
Using Jensen's inequality and $\log2\le(\log r)/8$
for sufficiently large $m$, we obtain
\[
\max_{1\le j\le r}
\mathbb E_{\widetilde P_j}
\left[d(\widehat\mPhi,\mPhi_j)^2\right]
\ge\frac{9\alpha_r^2}{64}.
\]
By the nonnegativity of the loss and the fact that $\mA$ has the same distribution under all hypotheses,
\begin{align}
\max_{1\le j\le r}\mathbb E_{\mPhi_j}
\left[d(\widehat\mPhi,\mPhi_j)^2\right]
\ge\mathbb P(\mathcal E)
\max_{1\le j\le r}\mathbb E_{\widetilde P_j}
\left[d(\widehat\mPhi,\mPhi_j)^2\right]
\ge\frac{9\alpha_r^2}{64}(1-\epsilon_{m,n}).
\end{align}
Since $U_{\mathrm{sep}}\subset\Gamma_\phi$, taking the infimum over estimators, using $\epsilon_{m,n}=o(1)$ and substituting the expression of $\alpha_r$ from \eqref{eq:packing-separation-alpha} yields
\[
R_\Phi(\Gamma_\phi;m,n,L,\sigma_z)
\ge c\frac{k_\phi}{N_{\mathrm{div}}^\phi}\delta_r^2.
\]

For $n\ge(1+\tau)m$, substituting the corresponding choice of $\delta_r$ in \eqref{eq:choice-deltar-ngem} and using $\Gamma_\phi\in\mathcal F_{a_\phi,b_\phi,k_\phi,m}$ gives
\[
\sup_{\mathcal C_\phi\in\mathcal F^\phi_{a_\phi,b_\phi,k_\phi,m}}
R_\Phi(\mathcal C_\phi;m,n,L,\sigma_z)
\ge c_{\tau,\varepsilon,x_{\min},x_{\max}}
\frac{\max\{\sigma_z^4,n^2\}k_\phi\log m}{m^2nL}.
\]
This proves Lemma~\ref{lem:phase-only lower bound n ge m}. For $n\le m$ and $\sigma_z^2>0$, the corresponding choice of $\delta_r$ in \eqref{eq:choice-deltar-mgen} gives
\[
\sup_{\mathcal C_\phi\in\mathcal F^\phi_{a_\phi,b_\phi,k_\phi,m}}
R_\Phi(\mathcal C_\phi;m,n,L,\sigma_z)
\ge c_{\varepsilon,x_{\max}}
\frac{\sigma_z^2\max\{\sigma_z^2,m\}k_\phi\log m}{m^2nL}.
\]
This proves Lemma~\ref{lem:phase-only lower bound n le m}.

\section{Simulation and discussion}
\label{sec:experiments}
In this section we present numerical experiments to illustrate the dependence of the joint loss and its signal and phase components on the problem parameters. We compare the observed behavior with the bounds in
Section~\ref{sec:main results} and discuss questions concerning the sharpness of the current theory.

\subsection{Experiment setup}
\label{subsec:synthetic-model}

\paragraph{Data generation.}
We generate the observations according to the model in Section~\ref{sec:model and assump}. For each Monte Carlo trial we draw a sensing matrix \(\mA\in\mathbb C^{m\times n}\) with i.i.d.\ \(\mathcal{CN}(0,1)\) entries and keep it fixed across all \(L\) looks.
The signal amplitude \(\vx_o\) is block-constant with \(k_x\) contiguous nearly equal-length blocks of approximately equal length. The block amplitudes are drawn independently from \([0.1,4]\)  and then rescaled so that $\|\vx_o\|_2^2/n=1$. Similarly, the phase vector \(\vphi_o\) is block-constant with \(k_\phi\) contiguous blocks, whose values are drawn independently from \([-\pi/3,\pi/3]\). The first phase block is fixed at zero to remove the global phase ambiguity. 
For the fixed block partitions used here, these signal and phase parameterizations satisfy the polynomial covering-number bounds in Section~\ref{sec:model and assump}, with complexity orders at most \(k_x\) and \(k_\phi\), respectively.
Both speckle and additive noise are generated independently according to Section~\ref{sec:model and assump} and form \(\vy_\ell\) using \eqref{eqn:def_y_ell}.

\paragraph{Estimator.}
From the observations, we form the sample covariance
\[
\widehat{\mC}
=
\frac{1}{L}\sum_{\ell=1}^L
\vy_\ell\vy_\ell^H .
\]
For a candidate signal--phase pair \(\vtheta=(\vx,\vphi)\), let \(\mSigma(\theta)\) denote the precision matrix defined in Section~\ref{sec:model and assump}.
Then conditional on \(\mA\), the negative log-likelihood per look, up to constants independent of \(\vtheta\), is
\begin{equation}
\label{eq:empirical_nll}
\mathcal{L}(\vtheta;\widehat{\mC})
=
-\log\det\mSigma(\vtheta)
+
\operatorname{Tr}\!\left(
\mSigma(\vtheta)\widehat{\mC}
\right).
\end{equation}
We use Adam~\cite{kingma2015adam} to approximately minimize \eqref{eq:empirical_nll} over the block-constant signal and phase classes with known block partitions, and refer to the resulting estimator as the Adam MLE. 
We also use a covariance-based initialization. The squared block amplitudes are estimated by least squares using the diagonal of \(\widehat{\mC}-\sigma_z^2 I_m\), clipped to $[0.1^2,4^2]$, and square-rooted to obtain the initial amplitudes. Given this signal initialization, the phase is initialized using eigenvector-based angular synchronization~\cite{singer2011angular}, and circularly averaged within each block, with the first block fixed at zero and the remaining block phases clipped to $[-\pi/3,\pi/3]$. Adam then jointly updates the $k_x$ block amplitudes and $k_\phi-1$ free block phases, using PyTorch~\cite{paszke2019pytorch}'s default moment-decay parameters $\beta_1 = 0.9, \beta_2 = 0.999$ and numerical-stability parameter $\epsilon=10^{-8}$. The separate learning-rate schedules for amplitudes and phases, together with the fixed iteration budget for each experimental setting, are chosen heuristically based on preliminary experiments. After each update, the block amplitudes and free block phases are clipped to their respective bounds, and the iterate with the smallest likelihood objective is returned.

For each setting, we perform $100$ independent Monte Carlo trials, redrawing the sensing matrix $\mA$, signal amplitude $\vx_o$, phase distortion $\mPhi_o$, speckle noise $\bm{w}$ and additive noise $\bm{z}$ in each trial. Within a trial, $\mA$, $\vx_o$ and $\mPhi_o$ are drawn only once and fixed across all $L$ looks, while $\bm{w}$ and $\bm{z}$ are independently redrawn for each look. 
\subsection{Results}
\label{subsec:numerical-results}

We report the empirical means of the joint loss and its normalized signal and phase components defined in Section~\ref{sec:model and assump}. The experimental parameters are specified in the figure captions.

\paragraph{Number of looks and noise level.}

\begin{figure}[htbp]
    \centering
    \includegraphics[width=0.8\textwidth]{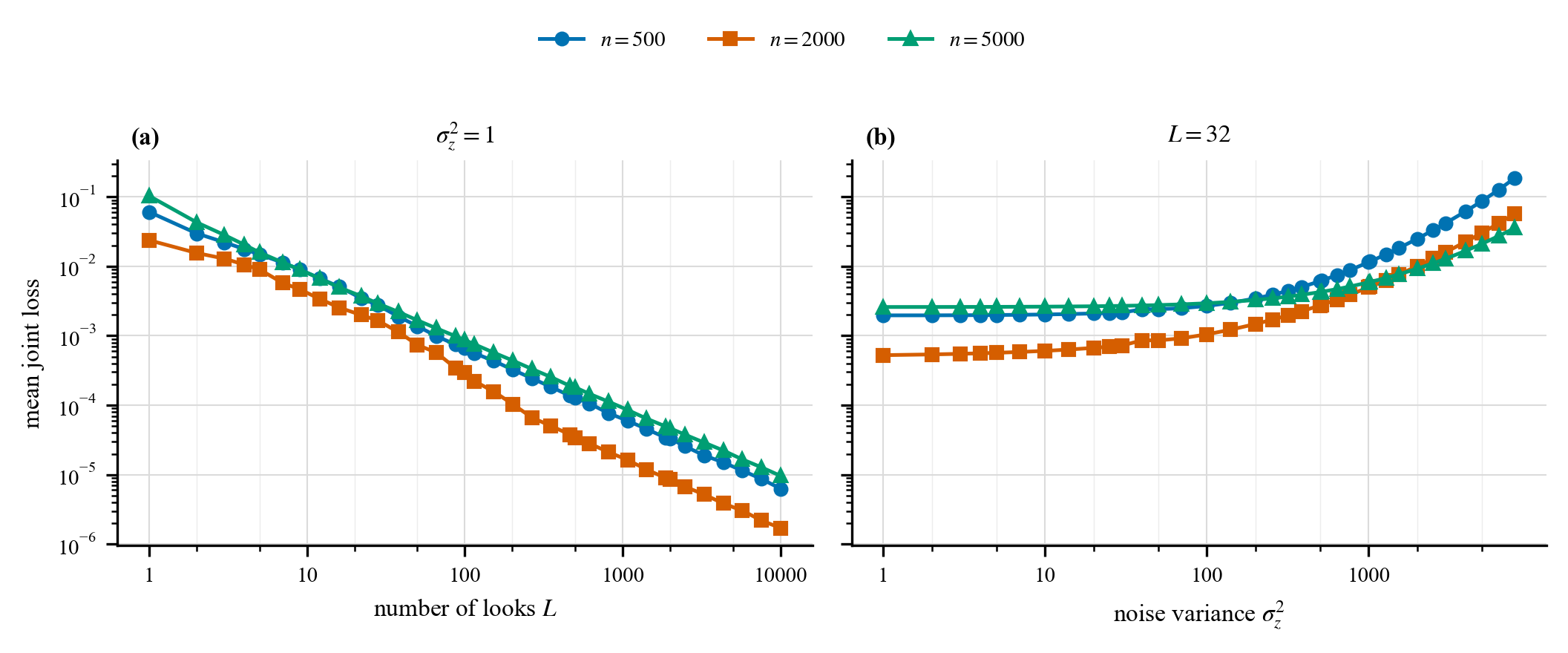}
    \caption{Performance of the Adam-based approximate MLE, measured by the mean joint loss. In both panels, we fix \(m=2000\) and \((k_x,k_\phi)=(128,200)\), and compare \(n\in\{500,2000,5000\}\), corresponding to the overdetermined, square, and underdetermined settings, respectively.
    (a) Joint loss versus the number of looks \(L\) with \(\sigma_z^2=1\). 
    (b) Joint loss versus the additive-noise variance \(\sigma_z^2\) with \(L=32\). 
    Results are averaged over independent Monte Carlo trials.}
    \label{fig:MSE_vs_L_sigma}
\end{figure}

Figure~\ref{fig:MSE_vs_L_sigma} shows the dependence of the joint loss on the number of looks and the additive-noise level. We fix \(m=2000\) and \((k_x,k_\phi)=(128,200)\), and consider \(n\in\{500,2000,5000\}\), corresponding to the overdetermined, square, and underdetermined settings, respectively.

Figure~\ref{fig:MSE_vs_L_sigma}(a) shows that the joint loss decreases with the number of looks in all three regimes.
This behavior is qualitatively consistent with the \(1/L\) dependence of the statistical term in our bounds.
The loss continues to decrease throughout the tested range, with no clear error floor.
This observation supports the conjecture that the \(L\)-independent term in \eqref{eq:joint-upper-bound} may be reduced or removed by a sharper analysis.

Figure~\ref{fig:MSE_vs_L_sigma}(b) shows that the joint loss changes little at low noise levels and increases as the additive noise becomes stronger.
This behavior agrees qualitatively with the dependence on \(\max\{\sigma_z^4,m^2,n^2\}\) in the theoretical bounds: the dimensional terms dominate at low noise, whereas the additive-noise term becomes important at higher noise levels.

\paragraph{Aspect ratio.}

\begin{figure}[htbp]
    \centering
    \includegraphics[width=0.8\textwidth]{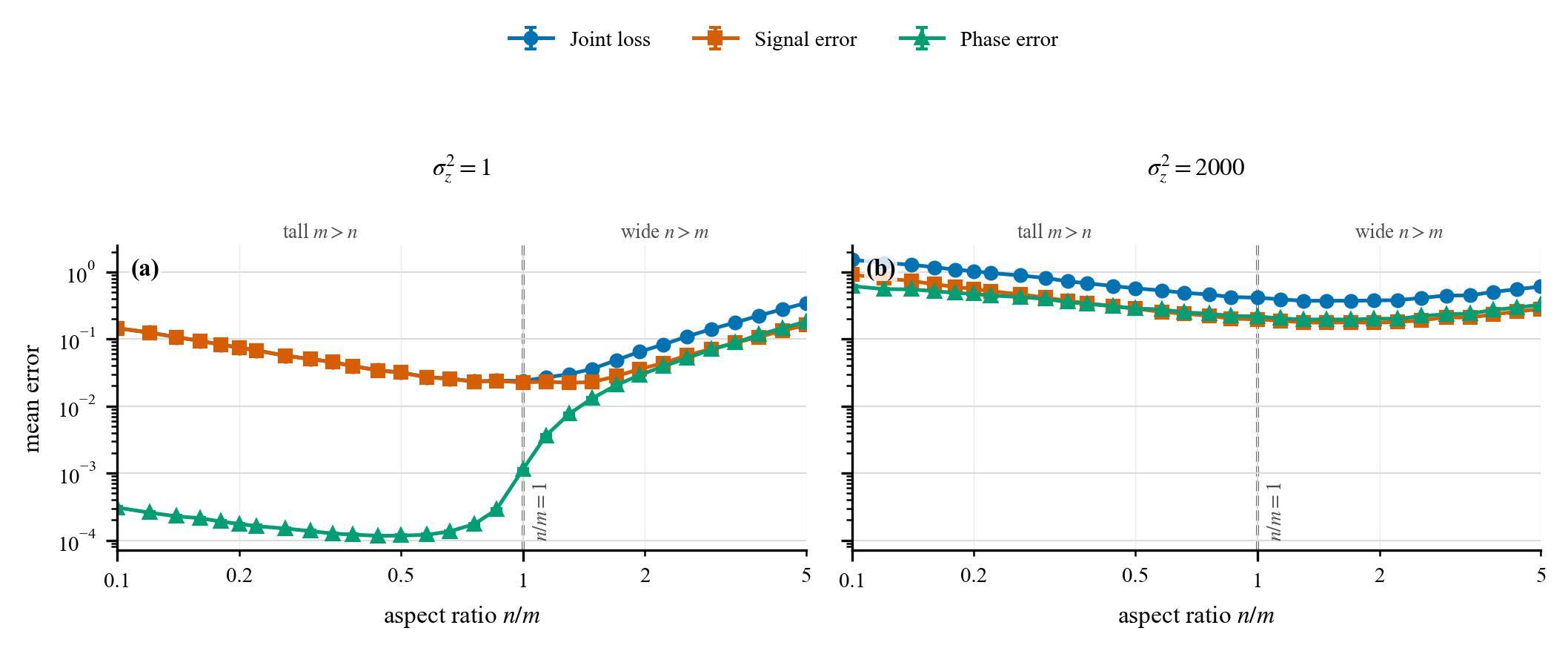}
    \caption{Performance of the Adam-based approximate MLE, measured by the mean joint loss and its signal and phase components, as a function of the aspect ratio \(n/m\).
    In both panels, we fix \(m=2000\) and \((L,k_x,k_\phi)=(1,128,200)\).
    (a) Low-noise setting with \(\sigma_z^2=1\).
    (b) Higher-noise setting with \(\sigma_z^2=2000\).
    The dashed vertical line marks the square case \(n=m\). Results are averaged over independent Monte Carlo trials.}
    \label{fig:MSE_vs_delta}
\end{figure}

Figure~\ref{fig:MSE_vs_delta} compares the joint loss and its two components as the aspect ratio \(n/m\) varies.
At the low noise level in panel~(a), the phase error is much smaller than the signal error in the overdetermined setting and increases rapidly as \(n\) passes \(m\).
At the higher noise level in panel~(b), the signal and phase errors are closer in scale, and the phase error varies more smoothly near \(n=m\).

This result phenomenon is qualitatively consistent with the different phase contributions in the lower bounds of Theorems~\ref{thm:joint-lb-m-ge-n} and~\ref{thm:joint-lb-n-ge-m}.
In the overdetermined regime, the phase contribution to the lower bound in \eqref{eq:joint_lower_bound_m_ge_n} contains the additional factor
\[
\frac{\sigma_z^2}{\max\{\sigma_z^2,m\}}.
\]
When \(\sigma_z^2\ll m\), the factor is of order \(\sigma_z^2/m\ll1\), substantially reducing the phase contribution to the bound. This is consistent with the much smaller phase error relative to the signal error in Figure~\ref{fig:MSE_vs_delta}(a).
When \(\sigma_z^2=m\), the factor equals one and this reduction disappears, consistent with the signal and phase errors being closer in scale in Figure~\ref{fig:MSE_vs_delta}(b).
The current joint upper bound does not capture the possible improvement in phase estimation at low noise.
Obtaining a matching phase upper bound and characterizing the rapid change near the square regime remain open.

The signal component varies more smoothly with the aspect ratio in both panels, with its smallest values near the square regime. This broad trend resembles the aspect-ratio dependence of the joint upper bound in \eqref{eq:upper-bound-low-noise-simplified}.

\paragraph{Phase complexity.}

\begin{figure}[htbp]
    \centering
    \includegraphics[width=0.8\textwidth]{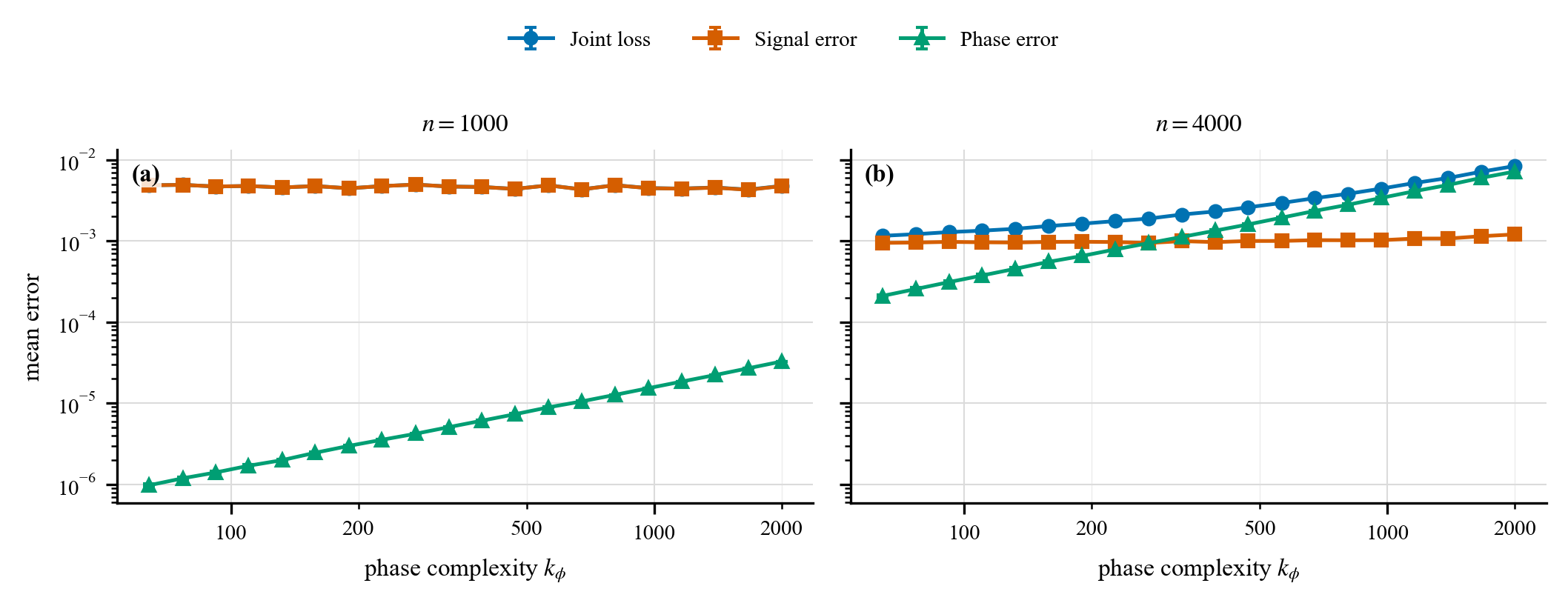}
    \caption{Performance of the Adam-based approximate MLE, measured by the mean joint loss and its signal and phase components, as a function of the phase complexity \(k_\phi\).
    In both panels, we fix \(m=2000\) and \((L,\sigma_z^2,k_x)=(32,1,128)\).
    (a) Overdetermined setting with \(n=1000\).
    (b) Underdetermined setting with \(n=4000\).
    Results are averaged over independent Monte Carlo trials.}
    \label{fig:MSE_vs_Kphi}
\end{figure}

Figure~\ref{fig:MSE_vs_Kphi} examines the effect of the phase complexity in the overdetermined and underdetermined settings.
In both panels, the phase error increases with \(k_\phi\), whereas the signal error remains nearly unchanged.
Consequently, the joint loss is nearly constant in panel~(a), where the signal component dominates, but increases in panel~(b) as the phase component becomes larger.

The increase in phase error agrees qualitatively with the phase-complexity dependence in our bounds.
The phase error is also substantially smaller in the low-noise overdetermined setting, consistent with the behavior in Figure~\ref{fig:MSE_vs_delta}. The increasing trend persists up to \(k_\phi=m\), beyond the complexity range covered by our lower-bound results.

The weak dependence of the signal error on \(k_\phi\) suggests that signal recovery may admit an error bound independent of the phase complexity.
Our signal-only upper bound provides a partial answer for a real-valued analogue under additional conditions. Extending this guarantee to the complex-valued model considered here, and determining whether it can be attained without explicit phase estimation, remain open.

\subsection{Discussion and open questions}
\label{subsec:discussion-open-questions}

The experiments illustrate the dependence of the joint error on the number of looks and the noise level, as well as the different behavior of signal and phase estimation.
Together with the gaps between our upper and lower bounds, these observations motivate several questions for further study.

\paragraph{The dependence on the number of looks.}
The current joint upper bound contains an \(L\)-independent term, whereas the available lower bounds have no analogous term.
Figure~\ref{fig:MSE_vs_L_sigma}(a) shows that the empirical loss continues to decrease throughout the tested range, with no clear saturation.
An open question is whether the \(L\)-independent term can be removed under suitable conditions, yielding a \(1/L\)-type risk bound when all looks share the same sensing matrix.

\paragraph{Sharper guarantees for phase estimation.}
The phase error is substantially smaller in the low-noise overdetermined regime, while the current joint upper bound does not capture this improvement. One limitation of our upper-bound analysis is the noise-level comparison used to handle \(0<\sigma_z^2<m\): we transfer a risk bound from the larger noise variance \(m\) to the smaller noise variance \(\sigma_z^2\).
This approach yields the correct bound for the signal-only problem in \cite{xing2025minimax}. In our model, however, it produces a bound that does not improve as \(\sigma_z^2\) decreases below \(m\), and therefore cannot capture any additional gain in phase accuracy at low noise.

The additional noise-dependent factor in the phase lower bound suggests a potentially smaller phase-error scale.
Whether this scale is attainable when the signal is also unknown remains open.
A related question concerns the behavior near \(n=m\), where the numerical results exhibit a rapid change in phase error and the current bounds do not provide a sharp characterization.
A more precise analysis of these regimes would clarify
the attainable phase accuracy and its relationship to
signal recovery.

\paragraph{Signal recovery with unknown phases.}
Figure~\ref{fig:MSE_vs_Kphi} shows that the signal error is nearly unchanged as the phase complexity increases. This observation raises the question of whether the complex-valued model admits a signal-only error bound that is independent of \(k_\phi\), while retaining rates comparable to those available when the phases are known.
Our signal-only result provides a partial answer for a real-valued analogue under additional conditions. Extending this guarantee to continuous phase distortions remains open. A further question is whether comparable signal recovery rates can be attained without explicitly estimating the phase distortion.

\section*{Acknowledgments}

The authors thank Sarah Anne Hagan and Professor Shirin Jalali for helpful discussions on this work.

\appendix

\section{Supplementary results}

This appendix contains the identifiability result mentioned in Remark~\ref{rem:uniqueness} and the auxiliary lemmas used in the proofs.
In Appendix~\ref{app:proof_uniqueness}, we show that if $m(m+1)/2\ge 2n-1$ and $\sigma_z^2>0$, then, with probability one over $\mA$, the true parameter $\vtheta_o$ is the unique minimizer of the population negative log-likelihood over $\Theta$.
Appendix~\ref{app:auxiliary_lemmas} first collects standard tools: bounds on the singular values of Gaussian matrices, concentration inequalities for quadratic forms and chi-square variables, a decoupling inequality, and the Kullback--Leibler divergence between complex Gaussian distributions.
It then gives two inequalities that compare the covariance matrices of two parameters, which are used in the proof of the upper bound.

\subsection{Identifiability Result}
\label{app:proof_uniqueness}

Recall that admissible amplitudes are strictly positive, and admissible
phase vectors satisfy $\phi_1=0$ and $\vphi\in[0,2\pi)^m$. Write 
\[
\mC(\vtheta)
=\mPhi(\vphi)\mA\operatorname{diag}(\vx^2)\mA^H\mPhi(\vphi)^H
+\sigma_z^2\mI_m.
\]
Then conditioned on $\mA$, the conditional population negative log-likelihood, up to constants, is
\[
\bar{\mathcal L}(\vtheta;\vtheta_o)
=\log\det\mC(\vtheta)
+\operatorname{Tr}\!\left(\mC(\vtheta)^{-1}\mC(\vtheta_o)\right).
\]

\begin{proposition}
\label{prop:population_identifiability}
Fix $\vtheta_o=(\vx_o,\vphi_o)\in\Theta$, and suppose that
$M:=m(m+1)/2\geq 2n-1$ and $\sigma_z^2>0$.
Let $\mA\in\mathbb C^{m\times n}$ have i.i.d.\ $\mathcal{CN}(0,1)$
entries. Then, with probability one over $\mA$, for every
$\vtheta\in\Theta$,
\begin{equation}
\bar{\mathcal L}(\vtheta;\vtheta_o)
-\bar{\mathcal L}(\vtheta_o;\vtheta_o)\geq0,
\label{eq:population_gap_nonnegative}
\end{equation}
with equality if and only if $\vtheta=\vtheta_o$.
\end{proposition}

\begin{proof}
\textbf{Step 1. Reduction to covariance equality.}
Both $\mC(\vtheta)$ and $\mC(\vtheta_0)$ are positive definite because $\sigma_z^2>0$. Thus,  $\mT=\mC(\vtheta)^{-1/2}\mC(\vtheta_o)\mC(\vtheta)^{-1/2}$ is Hermitian positive definite.
Let $\lambda_1,\ldots,\lambda_m$ be the eigenvalues of
$\mT$. Then
\begin{equation}
\bar{\mathcal L}(\vtheta;\vtheta_o)
-\bar{\mathcal L}(\vtheta_o;\vtheta_o)
=\sum_{j=1}^m(\lambda_j-\log\lambda_j-1)\geq0,
\label{eq:eigenvalue_gap}
\end{equation}
with equality holds exactly when
all $\lambda_j=1$, or equivalently $\mT=\mI_m$. This is the same as
\begin{equation}
\mC(\vtheta)=\mC(\vtheta_o).
\label{eq:covariance_equality}
\end{equation}

\medskip
\noindent\textbf{Step 2. Identifying the amplitudes.}
Set  
$\mS(\vx^2)=\mA\operatorname{diag}(\vx^2)\mA^H$, and
$\mD=\mPhi(\vphi_o)^H\mPhi(\vphi)$. Then equation \eqref{eq:covariance_equality} is equivalent to
\[
\mD\mS(\vx^2)\mD^H=\mS(\vx_o^2).
\]
Since $\mD$ is diagonal and unitary,
\[
|S_{pq}(\vx^2)|^2=|S_{pq}(\vx_o^2)|^2,
\qquad 1\leq p\leq q \leq m.
\]
Define
\[
\vv_{pq}
=(\overline{A_{p1}}A_{q1},\ldots,
  \overline{A_{pn}}A_{qn})^T,
\qquad
\mQ_{pq}(\mA)=\operatorname{Re}(\vv_{pq}\vv_{pq}^H).
\]
Each $\mQ_{pq}(\mA)$ is real symmetric. Since
$S_{pq}(\vx^2)=\vv_{pq}^H\vx^2$ and $\vx^2$ is real, 
$$|S_{pq}(\vx^2)|^2=(\vx^2)^T\mQ_{pq}(\mA)\vx^2.$$
Consequently,
\begin{equation}
(\vx^2-\vx_o^2)^T\mQ_{pq}(\mA)(\vx^2+\vx_o^2)=0,
\qquad p\leq q.
\label{eq:bilinear_ambiguity}
\end{equation}
We claim that, for almost every $\mA$, the following
holds for all $\vu,\vv\in\mathbb C^n$:
\begin{equation}
\vu^T\mQ_{pq}(\mA)\vv=0
\ \text{for all }p\leq q
\quad\Longrightarrow\quad
\vu=0\ \text{or}\ \vv=0.
\label{eq:generic_bilinear_identifiability}
\end{equation}
The entries of $\mQ_{pq}(\mA)$ are real polynomials in the real
and imaginary parts of $\mA$. By
Lemma~\ref{lem:polynomial_bilinear_generalize}, it therefore suffices
to verify \eqref{eq:generic_bilinear_identifiability} at a single
matrix $\mA_0$.

Choose distinct primes $t_1,\ldots,t_m$ and set
$(A_0)_{pk}=t_p^{k-1}$. Then the $M$
numbers $s_{pq}=t_pt_q$, $p\leq q$, are distinct.
For $\vu\in\mathbb C^n$, define
\[
f_{\vu}(z)=\sum_{k=1}^n u_kz^{k-1}.
\]
For this choice of $\mA_0$,
\[
\vu^T\mQ_{pq}(\mA_0)\vv
=f_{\vu}(s_{pq})f_{\vv}(s_{pq}).
\]
If these expressions vanish for every $p\leq q$, then the polynomial
$f_{\vu}f_{\vv}$ has at least $M\geq2n-1$ distinct roots.
Its degree is at most $2n-2$, so $f_{\vu}f_{\vv}\equiv0$.
Thus $f_{\vu}\equiv0$ or $f_{\vv}\equiv0$, which gives
$\vu=0$ or $\vv=0$.
This proves \eqref{eq:generic_bilinear_identifiability} at $\mA_0$.
Lemma~\ref{lem:polynomial_bilinear_generalize} extends it to almost every $\mA$, and hence to a
probability-one event under the complex Gaussian law of  $\mA$.

Applying \eqref{eq:generic_bilinear_identifiability} to
\eqref{eq:bilinear_ambiguity} gives $\vx^2=\vx_o^2$, since
$\vx^2+\vx_o^2\neq0$. Since amplitudes are positive, it yields
$\vx=\vx_o$.

\medskip
\noindent\textbf{Step 3. Identifying the phases.}
For the fixed vector $\vx_o^2>0$, each $|S_{1q}(\vx_o^2)|^2$,
$q=2,\ldots,m$, is a real polynomial in the real and imaginary
parts of $\mA$ that is not identically zero. Indeed, if every entry
of $\mA$ equals one, then
$S_{1q}(\vx_o^2)=\sum_{k=1}^n x_{o,k}^2>0$.
Since the zero set of each of these polynomials is Lebesgue-null
and $\mA$ has a density,
\[
S_{1q}(\vx_o^2)\neq0,\qquad q=2,\ldots,m,
\]
simultaneously almost surely.

We consider the intersection of this event with the probability-one event
from Step~2. The event from Step~2 gives $\vx^2=\vx_o^2$ and hence
\[
\mD\mS(\vx_o^2)\mD^H=\mS(\vx_o^2).
\]
The gauge constraint gives $D_{11}=1$. Comparing the $(1,q)$
entries yields
\[
S_{1q}(\vx_o^2)\overline{D_{qq}}=S_{1q}(\vx_o^2).
\]
Since $S_{1q}(\vx_o^2)\neq0$, we obtain $D_{qq}=1$ for
$q=2,\ldots,m$. Thus $\mD=\mI_m$ and
$\mPhi(\vphi)=\mPhi(\vphi_o)$. Both phase vectors lie in
$[0,2\pi)^m$, so $\vphi=\vphi_o$.

For every $\mA$ in the intersection of the two probability-one
events above, equality in \eqref{eq:eigenvalue_gap} implies
$\vtheta=\vtheta_o$ for every candidate $\vtheta\in\Theta$.
Conversely, substituting $\vtheta=\vtheta_o$ makes the likelihood
gap zero. Thus $\vtheta_o$ is the unique population minimizer
with probability one.
\end{proof}

\medskip
\begin{lemma}[]
\label{lem:polynomial_bilinear_generalize}
Let $Q_1(\vt),\ldots,Q_M(\vt)$ be $n\times n$ matrices whose entries are
real polynomials in $\vt\in\mathbb R^d$. Suppose that, for some
$\vt_0\in\mathbb R^d$,
\[
\vu^T\mQ_j(\vt_0)\vv=0\quad(j=1,\ldots,M),
\qquad \vu,\vv\in\mathbb C^n,
\]
implies $\vu=0$ or $\vv=0$. Then the same implication holds for
Lebesgue-almost every $\vt\in\mathbb R^d$.
\end{lemma}

\begin{proof}
We show that the real parameters for which the implication fails
lie in the zero set of a nonzero real polynomial.

Evaluate the polynomial entries of $Q_j$ at $\vt\in\mathbb C^d$.
For $\vt\in\mathbb C^d$ and $\vv\in\mathbb C^n$, and define
\[
\mA(\vt,\vv)
=
\begin{bmatrix}
\mQ_1(\vt)\vv & \cdots & \mQ_M(\vt)\vv
\end{bmatrix}
\in\mathbb C^{n\times M}.
\]
For fixed $\vt$ and $\vv$, the equations
$\vu^T\mQ_j(\vt)\vv=0$, $j=1,\ldots,M$, are equivalent to
$\mA(\vt,\vv)^T\vu=0$. They therefore admit a nonzero solution
$\vu$ if and only if $\operatorname{rank}\mA(\vt,\vv)<n$.
The hypothesis gives
\[
\operatorname{rank}\mA(\vt_0,\vv)=n
\qquad\text{for every }\vv\neq0,
\]
and in particular $M\ge n$.

The rank condition is unchanged when $\vv$ is multiplied by
a nonzero scalar. We therefore define
\[
Y
=
\left\{
(\vt,[\vv])\in
\mathbb C^d\times\mathbb P^{n-1}(\mathbb C):
\operatorname{rank}\mA(\vt,\vv)<n
\right\},
\]
where $[\vv]$ denotes the equivalence class of a nonzero vector
$\vv$ under multiplication by nonzero complex scalars.
The condition $\operatorname{rank}\mA(\vt,\vv)<n$ is equivalent
to the vanishing of all $n\times n$ minors of $\mA(\vt,\vv)$.
Each minor is polynomial in $\vt$ and $\vv$ and homogeneous
of degree $n$ in $\vv$. Hence $Y$ is Zariski closed.

Let $B$ be the image of $Y$ under the projection
$(\vt,[\vv])\mapsto\vt$. By the rank characterization above,
$B$ consists exactly of the complex parameters for which
the implication fails. The Projective Extension Theorem
\cite[Chapter~8, Section~5, Theorem~6]{cox2015ideals}
shows that $B$ is Zariski closed in $\mathbb C^d$.

Since $\vt_0\notin B$, there is a polynomial $g$ with complex
coefficients that vanishes on $B$ and satisfies
$g(\vt_0)\neq0$. Write $g=g_{\mathrm R}+i g_{\mathrm I}$,
where $g_{\mathrm R}$ and $g_{\mathrm I}$ are obtained by
taking the real and imaginary parts of the coefficients of $g$.
For real $\vt$, both $g_{\mathrm R}(\vt)$ and
$g_{\mathrm I}(\vt)$ are real, so both polynomials vanish
on $B\cap\mathbb R^d$. Since $g(\vt_0)\neq0$, at least one
of them, denoted by $h$, is nonzero at $\vt_0$.

Thus $B\cap\mathbb R^d$ is contained in the zero set of
the nonzero real polynomial $h$, which has Lebesgue measure
zero. The desired implication therefore holds for
Lebesgue-almost every $\vt\in\mathbb R^d$.
\end{proof}

\subsection{Auxiliary technical lemmas}
\label{app:auxiliary_lemmas}

\begin{lemma}[\cite{davidson2001local}, Theorem~II.13] \label{lem:svd-gaussian}
Let $A$ be an $m \times n$ random matrix with elements drawn i.i.d.\ from
$\mathcal{CN}(0,1)$. Then
\[
\mathbb{P}\!\left(
\sigma_{\max}(A) \le \sqrt{n} + \sqrt{m} + t
\right)
\ge 1 - 2 e^{-t^2/2},
\qquad t > 0.
\]
Moreover, if $m < n$, then for any $t > 0$,
\[
\mathbb{P}\!\left(
\sqrt{n} - \sqrt{m} - t
\le \sigma_{\min}(A)
\le \sigma_{\max}(A)
\le \sqrt{n} + \sqrt{m} + t
\right)
\ge 1 - 2 e^{-t^2/2}.
\]
\end{lemma}

\begin{lemma}[Complex Hanson Wright Inequality, \cite{bayraktar2020mass}, Theorem~3.6]\label{Thm:ComplexHW}
Let $\mA$ be a $N \times N$ Hermitian matrix and $\mX =
(X_1,\ldots,X_N) \in \mathbb{C}^N$ be a random vector whose components
$X_j$ are independent complex subgaussian random variables such that
$\max_j\bigl( \|\mathrm{Re}(X_j)\|_{\psi_2}, \|\mathrm{Im}(X_j)\|_{\psi_2} \bigr) \le K .$
Then there exists an absolute constant $c>0$ such that
\[
\mathbb{P}\bigl[|\mX^*\mA\mX - \mathbb{E}[\mX^*\mA\mX]| > t\bigr]
\le 2 \exp\Bigl( -c \min\Bigl\{\frac{t^2}{K^4 \|\mA\|_{HS}^2},
\frac{t}{K^2 \|\mA\|}\Bigr\} \Bigr)
\]
for each $t \ge 0$.
\end{lemma}


\begin{lemma}[Decoupling inequality]
\label{thm:decoupling-lemma}
Let $X_1,\ldots,X_N$ be independent random variables taking values in a measurable space $S$, and let $\tilde X_1,\ldots,\tilde X_N$ be an independent copy.
For each $r=1,\ldots,q$ and $i\ne j$, let $h_{ij}^{(r)}:S\times S\to\mathbb R$ be measurable.
Then there exists an absolute constant $D\ge1$ such that, for any positive thresholds $t_1,\ldots,t_q$,
\[
\begin{aligned}
&\mathbb P\left(
\exists\,r:
\left|
\sum_{i\ne j}h_{ij}^{(r)}(X_i,X_j)
\right|>t_r
\right)
\le
D\,\mathbb P\left(
\exists\,r:
\left|
\sum_{i\ne j}h_{ij}^{(r)}(X_i,\tilde X_j)
\right|>\frac{t_r}{D}
\right).
\end{aligned}
\]
In particular, $q=1$ gives the scalar decoupling inequality.
\end{lemma}

\begin{proof}
Apply \cite[Theorem~1]{de1995decoupling}, with order two, to the vector of kernels $(h_{ij}^{(1)},\ldots,h_{ij}^{(q)})\in\mathbb R^q$, using the norm
\[
\|(z_1,\ldots,z_q)\|
:= \max_{1\le r\le q}\frac{|z_r|}{t_r}
\]
and threshold $1$.
\end{proof}


\begin{lemma}[Laurent--Massart inequality{\cite[Equation~(4.3)]{laurent2000adaptive}}] \label{lem:laurent-massart}
Let $Z\sim\chi_d^2$ be a chi-square random variable with $d$ degrees of freedom. Then, for every $t>0$,
\[
\mathbb P\left( Z-d\ge 2\sqrt{dt}+2t \right)\le e^{-t}.
\]
\end{lemma}

\begin{lemma}\label{lem: bound of chi2}
Let \(\vg_1,\dots,\vg_L \overset{\mathrm{i.i.d.}}{\sim} \mathcal{CN}(0,I_m)\). Then there exists a universal constant \(c>0\) such that for any \(t>0\),
\begin{equation}
\mathbb P \Big(
\frac{1}{L}\sum_{\ell=1}^L \|\vg_\ell\|_2^2 \ge m+t
\Big)
\le
\exp\!\Big(
-cL\min\!\Big(\frac{t^2}{m},\, t\Big)
\Big).
\label{eqn:lem bound of chi2}
\end{equation}
\end{lemma}

\begin{proof}
By independence and the standard complex Gaussian distribution, we have
\[
2\sum_{\ell=1}^L\|\vg_\ell\|_2^2
\sim \chi^2_{2mL}.
\]
Applying Lemma~\ref{lem:laurent-massart} and dividing the threshold by $2L$, we obtain, for every $u>0$,
\[
\mathbb P\left(
\frac{1}{L}\sum_{\ell=1}^L\|\vg_\ell\|_2^2
\ge
m+\sqrt{\frac{2mu}{L}}+\frac{u}{L}
\right)
\le e^{-u}.
\]
For a given $t>0$, choose
\[
u:=\frac{L}{8}\min\left\{\frac{t^2}{m},\,t\right\}.
\]
Then
\[
\sqrt{\frac{2mu}{L}}+\frac{u}{L}
\le \frac{t}{2}+\frac{t}{8}
\le t.
\]
Consequently, the preceding probability bound implies
\eqref{eqn:lem bound of chi2} with $c=1/8$.
\end{proof}

\begin{lemma}[Complex Gaussian KL divergence, {\cite[Equation~(55)]{ferrante2012time}}]\label{lem:KL-div-Gaussian}
If $\mathbb{Q}_j \sim \mathcal{CN}(0,\mLambda_j)$, $j=1,2$, are two zero-mean circularly symmetric complex Gaussian distributions on $\mathbb{C}^d$, then
\[
D_{\mathrm{KL}}(\mathbb{Q}_1\|\mathbb{Q}_2)
=
\log\frac{\det\mLambda_2}{\det\mLambda_1} +
\operatorname{Tr}(\mLambda_2^{-1}\mLambda_1)-d.
\]
\end{lemma}

\begin{lemma}\label{lem:low_and_up_bound_Tr}
Let $\mSigma_o=(\sigma_z^2 I_m+\mPhi_o\,\mA|\mX_o|^2\mA^H\mPhi_o^H)^{-1}$ and $\mSigma_1=(\sigma_z^2 I_m+\mPhi_1\,\mA|\mX_1|^2\mA^H\mPhi_1^H)^{-1}$. Define
\[
\bm{\Delta}
:=
\mSigma_o^{-1}-\mSigma_1^{-1}
=
\mPhi_o\mA|\mX_o|^2\mA^H\mPhi_o^H
-
\mPhi_1\mA|\mX_1|^2\mA^H\mPhi_1^H.
\]
Then
\begin{equation}
\Tr\!\left[
\mSigma_o^{-1}(\mSigma_1-\mSigma_o)
\mSigma_o^{-1}(\mSigma_1-\mSigma_o)
\right]
\ge
\frac{\|\bm{\Delta}\|_{HS}^2}{\Lambda_{\max}^2}.
\label{eqn:lem Z lower bound}
\end{equation}
\end{lemma}

\begin{proof}
Since
$\mSigma_o^{-1}(\mSigma_1-\mSigma_o)
=\bm{\Delta}\mSigma_1$
and $\bm{\Delta}$ is Hermitian, cyclicity of trace gives
\begin{align*}
&\Tr\!\left[
\mSigma_o^{-1}(\mSigma_1-\mSigma_o)
\mSigma_o^{-1}(\mSigma_1-\mSigma_o)
\right]
= \Tr(\bm{\Delta}\mSigma_1\bm{\Delta}\mSigma_1)
= \bigl\| \mSigma_1^{1/2}\bm{\Delta}\mSigma_1^{1/2} \bigr\|_{HS}^2
\ge \lambda_{\min}(\mSigma_1)^2\|\bm{\Delta}\|_{HS}^2
\ge \frac{\|\bm{\Delta}\|_{HS}^2}{\Lambda_{\max}^2},
\end{align*}
where the last inequality follows from
$\lambda_{\min}(\mSigma_1)
=1/\lambda_{\max}(\mSigma_1^{-1})
\ge1/\Lambda_{\max}$.
\end{proof}

\begin{lemma}\label{lem: cov trace and 2 upbound}
Let $\mSigma_o=(\sigma_z^2 I_m+\mPhi_o\,\mA|\mX_o|^2\mA^H\mPhi_o^H)^{-1}$ and $\mSigma_1=(\sigma_z^2 I_m+\mPhi_1\,\mA|\mX_1|^2\mA^H\mPhi_1^H)^{-1}$, and assume $\Lambda_{\min}^2>0$.
Then
\begin{subequations}
\begin{equation}
\|\mSigma_1-\mSigma_o\|_2
\le
\frac{E_{\max}}{\Lambda_{\min}^2}
\left(
\bigl\||\vx_1|^2-|\vx_o|^2\bigr\|_\infty
+
2x_{\max}^2\|\mPhi_1-\mPhi_o\|_2
\right).
\label{eqn:lem preci error upbound 1}
\end{equation}
Moreover,
\begin{equation}
\left|
\Tr\!\left(
\mSigma_o^{-1/2}(\mSigma_1-\mSigma_o)\mSigma_o^{-1/2}
\right)
\right|
\le
\frac{mE_{\max}}{\Lambda_{\min}}
\left(
\bigl\||\vx_1|^2-|\vx_o|^2\bigr\|_\infty
+
2x_{\max}^2\|\mPhi_1-\mPhi_o\|_2
\right),
\label{eqn:lem cov upbound1}
\end{equation}
and
\begin{equation}
\left\|
\mSigma_o^{-1/2}(\mSigma_1-\mSigma_o)\mSigma_o^{-1/2}
\right\|_2
\le
\frac{E_{\max}}{\Lambda_{\min}}
\left(
\bigl\||\vx_1|^2-|\vx_o|^2\bigr\|_\infty
+
2x_{\max}^2\|\mPhi_1-\mPhi_o\|_2
\right).
\label{eqn:lem cov upbound2}
\end{equation}
\end{subequations}
\end{lemma}

\begin{proof}
Write $\bm{\Delta}:=\mSigma_o^{-1}-\mSigma_1^{-1}$. The inverse identity $\mSigma_1-\mSigma_o=\mSigma_1\bm{\Delta}\mSigma_o$ and the definition of $\Lambda_{\min}^2$ in \eqref{eq:def-Lambda-of-Sigma} implies
\begin{equation}
\|\mSigma_1-\mSigma_o\|_2
\le
\|\mSigma_1\|_2\|\mSigma_o\|_2\|\bm{\Delta}\|_2
\le
\frac{\|\bm{\Delta}\|_2}{\Lambda_{\min}^2}.
\label{eqn:lem preci error upbound 2}
\end{equation}

To bound $\bm{\Delta}$, set $\bm S_1:=\mA|\mX_1|^2\mA^H$ and separate the amplitude and phase differences:
\begin{align*}
\bm{\Delta}
&=
\mPhi_o\mA
\bigl(|\mX_o|^2-|\mX_1|^2\bigr)
\mA^H\mPhi_o^H
+
(\mPhi_o-\mPhi_1)\bm S_1\mPhi_o^H
+
\mPhi_1\bm S_1(\mPhi_o^H-\mPhi_1^H).
\end{align*}
Since the phase matrices are unitary, $\|\mA\|_2^2=E_{\max}$, and $\|\bm S_1\|_2\le x_{\max}^2E_{\max}$ from definition of $E_{\max}$ in \eqref{eq:def-Emax-norm2A}, we obtain
\begin{equation}
\|\bm{\Delta}\|_2
\le E_{\max} \left(
\bigl\||\vx_1|^2-|\vx_o|^2\bigr\|_\infty
+ 2x_{\max}^2\|\mPhi_1-\mPhi_o\|_2
\right).
\label{eqn:lem preci error upbound 3}
\end{equation}
Together with \eqref{eqn:lem preci error upbound 2}, this proves \eqref{eqn:lem preci error upbound 1}.

For the remaining bounds, define
\[
\bm M
:=
\mSigma_o^{-1/2}(\mSigma_1-\mSigma_o)\mSigma_o^{-1/2},
\qquad
\bm B:=\mSigma_o^{-1/2}\mSigma_1^{1/2}.
\]
Then
\[
\bm M=\bm B\bm B^H-\mI_m,
\qquad
\mSigma_1^{1/2}\bm{\Delta}\mSigma_1^{1/2}
=
\bm B^H\bm B-\mI_m.
\]
These two Hermitian matrices have the same eigenvalues and
therefore the same spectral norm. Hence
\begin{equation}
\|\bm M\|_2
=
\bigl\|
\mSigma_1^{1/2}\bm{\Delta}\mSigma_1^{1/2}
\bigr\|_2
\le
\|\mSigma_1\|_2\|\bm{\Delta}\|_2
\le
\frac{\|\bm{\Delta}\|_2}{\Lambda_{\min}}.
\label{eqn:lem cov upbound3}
\end{equation}
Also,
\begin{equation}
|\Tr(\bm M)|
\le
m\|\bm M\|_2
\le
\frac{m}{\Lambda_{\min}}\|\bm{\Delta}\|_2.
\label{eqn:lem cov upbound4}
\end{equation}
Substituting \eqref{eqn:lem preci error upbound 3}
into \eqref{eqn:lem cov upbound3}
and \eqref{eqn:lem cov upbound4}
proves \eqref{eqn:lem cov upbound2}
and \eqref{eqn:lem cov upbound1}, respectively.
\end{proof}

\section{High-probability part for Upper Bound}

The following lemma gives the concentration of quadratic forms, the proof is similar to the one in \cite{xing2025minimax} but replace with complex model. We give the proof in following for reading convenience.

\begin{lemma}\label{lem:quad-yt}
Let $\mSigma_o\in\mathbb C^{m\times m}$ be a fixed Hermitian positive definite matrix, and let $\hat{\mSigma},\tilde{\mSigma}\in\mathbb C^{m\times m}$ be fixed Hermitian matrices.
Suppose that $\vy_1,\ldots,\vy_L\overset{\mathrm{i.i.d.}}{\sim} \mathcal{CN}(0,\mSigma_o^{-1})$, and define
\[
\mM
:= \mSigma_o^{-1/2}(\hat{\mSigma}-\tilde{\mSigma})\mSigma_o^{-1/2}.
\]
Then there exists an absolute constant $c>0$ such that, for every $t>0$,
\begin{equation}
\begin{aligned}
&\mathbb P\left(
\left|
\frac1L\sum_{\ell=1}^L
\vy_\ell^H(\hat{\mSigma}-\tilde{\mSigma})\vy_\ell
-
\Tr\bigl(
(\hat{\mSigma}-\tilde{\mSigma})\mSigma_o^{-1}
\bigr)
\right|>t \right)
\le
2\exp\left(
-cL\min\left\{
\frac{t^2}{\|\mM\|_{HS}^2},\frac{t}{\|\mM\|_2}
\right\}\right).
\label{eq:quadratic-form-concentration}
\end{aligned}
\end{equation}
When $\mM=0$, the centered expression is identically zero, and the right-hand side is interpreted as zero.
\end{lemma}

\begin{proof}
The case $\mM=0$ is immediate. Otherwise, define $\vg_\ell:=\mSigma_o^{1/2}\vy_\ell$ and set
\[
\vg:=(\vg_1^H,\ldots,\vg_L^H)^H
\sim\mathcal{CN}(0,\mI_{Lm}),
\qquad
\bm K:=\frac1L\operatorname{diag}(\mM,\ldots,\mM).
\]
Then
\[
\frac1L\sum_{\ell=1}^L
\vy_\ell^H(\hat{\mSigma}-\tilde{\mSigma})\vy_\ell
-
\Tr\bigl(
(\hat{\mSigma}-\tilde{\mSigma})\mSigma_o^{-1}
\bigr)
=
\vg^H\bm K\vg-\Tr(\bm K).
\]
Since $\bm K$ is Hermitian, $\E[\vg^H\bm K\vg]=\Tr(\bm K)$, and
\[
\|\bm K\|_{HS}^2=\frac{\|\mM\|_{HS}^2}{L},
\qquad
\|\bm K\|_2=\frac{\|\mM\|_2}{L},
\]
the complex Hanson--Wright inequality in Lemma~\ref{Thm:ComplexHW} gives the result.
\end{proof}

\begin{lemma}\label{lem:HS_lower_bound}
Let $\mA\in\mathbb C^{m\times n}$ have i.i.d.\ $\mathcal{CN}(0,1)$ entries, where $m\ge2$ and $n\ge1$.
Let $\mX_o,\mX_1,\ldots,\mX_q$ be fixed diagonal $n\times n$ matrices with diagonal entries in $[x_{\min},x_{\max}]$, where $q\ge1$ and $x_{\min}>0$. Let $\mU_1,\ldots,\mU_q$ be fixed diagonal unitary $m\times m$ matrices. For $r=1,\ldots,q$, define
\[
\mu_r
:=
m(m-1)\|\mX_r^2-\mX_o^2\|_{HS}^2
+
2\bigl(m^2-|\Tr(\mU_r)|^2\bigr)
\Tr(\mX_r^2\mX_o^2).
\]
Then there exist absolute constants $c,C>0$ such that,
for every $t>0$,
\begin{equation}
\begin{aligned}
&\mathbb P\left(
\exists\,r:
\left\|
\mA\mX_r^2\mA^H
-
\mU_r\mA\mX_o^2\mA^H\mU_r^H
\right\|_{HS}^2
<
\frac{\mu_r}{2}-t
\right)
\le
Cq\exp\left(
-\frac{ct}{x_{\max}^4(m+n)}
\right)
+
C\exp\bigl(-c(m+n)\bigr).
\end{aligned}
\label{eq:HS-lower-finite-family}
\end{equation}
\end{lemma}

\begin{proof}
We first consider a fixed $r$ and suppress the index $r$ in the following calculations.
Let $\va_i^H$ denote the $i$-th row of $\mA$, and write $\mU=\operatorname{diag}(u_1,\ldots,u_m)$. We first expand the squared HS norm into its diagonal and off-diagonal terms:
\begin{equation}\label{eq:app-expand-HSnorm}
\left\|
\mA\mX^2\mA^H-\mU\mA\mX_o^2\mA^H\mU^H
\right\|_{HS}^2
=
\sum_{i\ne j}|\va_i^H (\mX^2-u_i\overline{u_j}\mX_o^2) \va_j|^2
+
\sum_{i=1}^m
|\va_i^H(\mX^2-\mX_o^2)\va_i|^2.
\end{equation}
Since the diagonal term is nonnegative, it is enough to lower bound the off-diagonal term, which we denote by
\[
\Delta_{\mathrm{off}}
:=
\sum_{i\ne j}|\va_i^H\mB_{ij}\va_j|^2,
\qquad
\mB_{ij}:=\mX^2-u_i\overline{u_j}\mX_o^2.
\]
Since $\va_i$ and $\va_j$ are independent standard complex Gaussian vectors for $i\ne j$, we have
\begin{equation*}
\E|\va_i^H\mB_{ij}\va_j|^2
=\|\mB_{ij}\|_{HS}^2
=\|\mX^2-\mX_o^2\|_{HS}^2
+|1-u_i\overline{u_j}|^2\Tr(\mX^2\mX_o^2),
\end{equation*}
where the second equality uses $|u_i|=|u_j|=1$. Summing over all $i\ne j$ and using the identity $\sum_{i\ne j}|1-u_i\overline{u_j}|^2 = 2\bigl(m^2-|\Tr(\mU)|^2\bigr)$,
we obtain the expectation
\begin{equation}
\mu = \E[\Delta_{\mathrm{off}}]
=
\sum_{i\ne j} \E|\va_i^H\mB_{ij}\va_j|^2
=
m(m-1)\|\mX^2-\mX_o^2\|_{HS}^2
+
2\bigl(m^2-|\Tr(\mU)|^2\bigr)
\Tr(\mX^2\mX_o^2).
\label{eqn:expect_trace}
\end{equation}

We next show that $\Delta_{\mathrm{off}}$ concentrates around its mean $\mu$. The main difficulty is that the summands in $\Delta_{\mathrm{off}}$ are not independent, as different pairs may share a row of $\mA$. To handle this dependence, we introduce $\tilde{\mA}$ as an independent copy of $\mA$, with rows $\tilde{\va}_j^H$, and define
\begin{align}\label{eq:app-def-tildeDelta}
\tilde{\Delta}_{\mathrm{off}}
:= \sum_{i\ne j}|\va_i^H\mB_{ij}\tilde{\va}_j|^2
= \sum_{i=1}^m \va_i^H
\big(
\sum_{j\ne i}\mB_{ij}\tilde{\va}_j\tilde{\va}_j^H\mB_{ij}^H
\big)\va_i.
\end{align}
Then we apply the decoupling inequality to relate the fluctuations of $\Delta_{\mathrm{off}}$ to those of $\tilde{\Delta}_{\mathrm{off}}$. 
To handle all $q$ comparisons simultaneously, we restore the index $r$ and use the same independent copy $\tilde{\mA}$ for all comparisons. For each $r=1,\ldots,q$, let $\mB_{ij}^{(r)} :=\mX_r^2-u_{r,i}\overline{u_{r,j}}\mX_o^2$, where $u_{r,i}$ denotes the $i$-th diagonal entry of $\mU_r$.
We define $\Delta_{\mathrm{off},r}$ and $\tilde{\Delta}_{\mathrm{off},r}$ accordingly, and introduce the centered kernels
\[
h_{ij}^{(r)}(\bm a,\bm b)
:=
|\bm a^H\mB_{ij}^{(r)}\bm b|^2 -\|\mB_{ij}^{(r)}\|_{HS}^2.
\]
Then
\begin{align*}
\Delta_{\mathrm{off},r}-\mu_r &=\sum_{i\ne j} h_{ij}^{(r)}(\va_i,\va_j),\qquad
\tilde{\Delta}_{\mathrm{off},r}-\mu_r =\sum_{i\ne j} h_{ij}^{(r)}(\va_i,\tilde{\va}_j).
\end{align*}
By Lemma~\ref{thm:decoupling-lemma}, there exists an absolute constant $C\ge1$ such that, for every $t>0$,
\begin{equation}
\begin{aligned}
&\mathbb P\left(
\exists r: \Delta_{\mathrm{off},r}<\frac{\mu_r}{2}-t
\right)
\le
\mathbb P\left(
\exists r: |\Delta_{\mathrm{off},r}-\mu_r| >\frac{\mu_r}{2}+t
\right)
\le 
C\,\mathbb P\left(
\exists r: |\tilde{\Delta}_{\mathrm{off},r}-\mu_r|
>\frac{\mu_r/2+t}{C} \right).
\end{aligned}
\label{eq:app-decoupling-off}
\end{equation}
It therefore remains to control $\tilde{\Delta}_{\mathrm{off},r}-\mu_r$. 

We again consider a fixed $r$ and suppress the index $r$ in the following calculations.
Conditioned on $\tilde{\mA}$, $\tilde{\Delta}_{\mathrm{off}}$ in \eqref{eq:app-def-tildeDelta} is a sum of independent Gaussian quadratic forms. Define the block-diagonal matrix
\begin{align}\label{eq:app-def-block-D-matrix}
\mD:&=\operatorname{diag}(\mD_1,\ldots,\mD_m),
\qquad
\mD_i:=
\sum_{j\ne i}
\mB_{ij}\tilde{\va}_j\tilde{\va}_j^H\mB_{ij}^H.
\end{align}
Since each block $\mD_i$ is positive semidefinite, so is $\mD$. Moreover,
\begin{align*}
\E[\tilde{\Delta}_{\mathrm{off}} \mid\tilde{\mA}] =\Tr(\mD).
\end{align*}
The conditional mean $\Tr(\mD)$ is itself random, but it can also be expressed as a sum of quadratic forms with deterministic matrices. Rearranging the sums and then take expectation gives
\begin{equation}
\Tr(\mD)
=
\sum_{j=1}^m \tilde{\va}_j^H
\left( \sum_{i\ne j}\mB_{ij}^H\mB_{ij} \right)
\tilde{\va}_j
= \sum_{j=1}^m \tilde{\va}_j^H \mF_j \tilde{\va}_j ,
\qquad
\E[\Tr(\mD)]= \Tr(\mF)=\mu.
\label{eq:HS-conditional-mean}
\end{equation}
where
\begin{align*}
\mF:&=\operatorname{diag}(\mF_1,\ldots,\mF_m),
\qquad
\mF_j:=\sum_{i\ne j}\mB_{ij}^H\mB_{ij}.
\end{align*}

We can therefore control the fluctuations of $\tilde{\Delta}_{\mathrm{off}}$ around $\mu$ in two steps, using the decomposition
\[
\tilde{\Delta}_{\mathrm{off}}-\mu
=
\bigl(\tilde{\Delta}_{\mathrm{off}}-\Tr(\mD)\bigr)
+
\bigl(\Tr(\mD)-\mu\bigr).
\]
Both terms can be controlled by the Hanson--Wright inequality. 
By stacking the vectors $\va_1,\ldots,\va_m$ and $\tilde{\va}_1,\ldots,\tilde{\va}_m$ into $\bm a,\tilde{\bm a}\in\mathbb C^{mn}$, respectively, we can write
\[
\tilde{\Delta}_{\mathrm{off}}
=\bm a^H\mD\bm a,
\qquad
\Tr(\mD)=\tilde{\bm a}^H\mF\tilde{\bm a},
\]
where $\mD$ and $\mF$ are block-diagonal matrices.
Conditional on $\tilde{\mA}$, $\mD$ is fixed, while $\mF$ is deterministic.
Applying Lemma~\ref{Thm:ComplexHW} to these two quadratic forms gives, for every $t>0$,
\begin{align} \label{eq:app-HS-Delta-trD}
\mathbb P\left(
\left|\tilde{\Delta}_{\mathrm{off}}-\Tr(\mD)\right|>t
\mid\tilde{\mA} \right)
\le
2\exp\left[
-c\min\left\{
\frac{t^2}{\|\mD\|_{HS}^2},\frac{t}{\|\mD\|_2}
\right\}\right],
\end{align}
and
\begin{align}\label{eq:app-HS-TrD-mu}
\mathbb P\left( \left|\Tr(\mD)-\mu\right|>t \right)
\le
2\exp\left[
-c\min\left\{
\frac{t^2}{\|\mF\|_{HS}^2},\frac{t}{\|\mF\|_2}
\right\}\right].
\end{align}
Here the first inequality is conditional on $\tilde{\mA}$, while the second applies directly since $\mF$ is deterministic.

We next bound the spectral and HS norms appearing in these inequalities. For each block $\mD_i$ defined in \eqref{eq:app-def-block-D-matrix}, the inequality $(\bm a-\bm b)(\bm a-\bm b)^H \preceq 2\bm a\bm a^H+2\bm b\bm b^H$ gives
\[
\begin{aligned}
\mD_i
&=\sum_{j\ne i}
\mB_{ij}\tilde{\va}_j\tilde{\va}_j^H\mB_{ij}^H
\preceq
2\mX^2\tilde{\mA}^H\tilde{\mA}\mX^2
+
2\mX_o^2\tilde{\mA}^H\tilde{\mA}\mX_o^2.
\end{aligned}
\]
Since $\|\mX\|_2,\|\mX_o\|_2\le x_{\max}$, we obtain
\begin{equation}
\|\mD\|_2 =\max_i\|\mD_i\|_2 \le 4x_{\max}^4\|\tilde{\mA}\|_2^2.
\label{eq:app-2norm-Dmatrix}
\end{equation}
Similarly, using $\|\mB_{ij}\|_2 = \|\mX^2-u_i\overline{u_j}\mX_o^2\|_2 \le 2x_{\max}^2$, we have
\[
\|\mF\|_2 =\max_j\big\| \sum_{i\ne j}\mB_{ij}^H\mB_{ij} \big\|_2 \le 4(m-1)x_{\max}^4.
\]
Consider the event
\[
\mathcal E_{\tilde A}
:=
\left\{ \|\tilde{\mA}\|_2\le2(\sqrt m+\sqrt n) \right\},
\qquad b_{m,n}:=32x_{\max}^4(m+n).
\]
By Lemma~\ref{lem:svd-gaussian},
\[
\mathbb P(\mathcal E_{\tilde A}^{\,c}) \le 2e^{-(m+n)/2}.
\]
It follows that
\[
\|\mD\|_2\le b_{m,n} \quad\text{on }\mathcal E_{\tilde A},
\qquad \|\mF\|_2\le b_{m,n}.
\]
Finally, by $\|\bm H\|_{HS}^2\le\|\bm H\|_2\Tr(\bm H)$ for any positive semidefinite matrix $\bm H$, we get
\[
\|\mD\|_{HS}^2\le b_{m,n}\Tr(\mD)
\quad \text{on }\mathcal E_{\tilde A},
\qquad \|\mF\|_{HS}^2\le b_{m,n}\mu.
\]

Note that $C\ge 1$ is an absolute constant in \eqref{eq:app-decoupling-off}, then for any $s>0$, taking $t=\mu/(8C)+s$ in the second Hanson--Wright inequality in \eqref{eq:app-HS-TrD-mu} and using $\|\mF\|_{HS}^2\le b\mu$ and $\|\mF\|_2\le b$,
we obtain
\begin{align}\label{eq:app-HS-TrD-mu-concentration}
&\mathbb P\left( |\Tr(\mD)-\mu|>\frac{\mu}{8C}+s \right)
\le
2\exp\left[
-c\min\left\{
\frac{(\mu/(8C)+s)^2}{b_{m,n}\mu},\frac{\mu/(8C)+s}{b_{m,n}}
\right\}\right]
\le 2\exp\left(-\frac{cs}{b_{m,n}}\right),
\end{align}
where $C$ is an absolute constant absorbed into $c$. This inequality also holds when $\mu=0$ since $\mF=0$. 
Similarly, conditional on any $\tilde{\mA}\in\mathcal E_{\tilde A}$, taking $t=\Tr(\mD)/(8C)+s$ in \eqref{eq:app-HS-Delta-trD} and using $\|\mD\|_{HS}^2\le b_{m,n}\Tr(\mD)$ and $\|\mD\|_2\le b_{m,n}$, we have
\[
\begin{aligned}
&\mathbb P\left(\left|\tilde{\Delta}_{\mathrm{off}}-\Tr(\mD)\right|
>\frac{\Tr(\mD)}{8C}+s \;\middle|\;\tilde{\mA} \right)
\le
2\exp\left[ -c\min\left\{
\frac{(\Tr(\mD)/(8C)+s)^2}{b_{m,n}\Tr(\mD)}, \frac{\Tr(\mD)/(8C)+s}{b_{m,n}}
\right\} \right]
\le 2e^{-\frac{cs}{b_{m,n}}}.
\end{aligned}
\]
When $\Tr(\mD)=0$, the deviation is again identically zero. Averaging over $\tilde{\mA}$ restricted to $\mathcal E_{\tilde A}$ yields
\begin{equation}
\mathbb P\left(
\left\{
\left|\tilde{\Delta}_{\mathrm{off}}-\Tr(\mD)\right|
>\frac{\Tr(\mD)}{8C}+s \right\}
\cap\mathcal E_{\tilde A} \right)
\le 2e^{-\frac{cs}{b_{m,n}}}.
\label{eq:app-concentration-conditional}
\end{equation}
Combining these two bounds, whenever neither deviation exceeds its threshold, we have
\begin{align*}
|\tilde{\Delta}_{\mathrm{off}}-\mu| 
&\le
|\tilde{\Delta}_{\mathrm{off}}-\Tr(\mD)|
+
|\Tr(\mD)-\mu| \\
&\le \frac{\Tr(\mD)}{8C}+s + \frac{\mu}{8C}+s \\
&\le \frac{\mu}{2C} + 3s.
\end{align*}
where the last inequality uses $C\ge1$ and inequality $\Tr(\mD)\le\left(1+\frac{1}{8C}\right)\mu+s$ from \eqref{eq:app-HS-TrD-mu-concentration}.
Taking $s=t/(3C)$ and applying a union bound
therefore gives
\begin{equation}
\mathbb P\left(
\left\{
|\tilde{\Delta}_{\mathrm{off}}-\mu|
>\frac{\mu/2+t}{C}
\right\}
\cap\mathcal E_{\tilde A}
\right)
\le 4e^{-\frac{ct}{b_{m,n}}},
\label{eq:app-concentration-decoupled}
\end{equation}
after adjusting the absolute constant $c$.

We now return to all $q$ comparisons. Since the same independent copy $\tilde{\mA}$ is used for every comparison, the event $\mathcal E_{\tilde A}$ is common to all of them. Applying \eqref{eq:app-concentration-decoupled} and a union bound, we obtain
\[
\begin{aligned}
&\mathbb P\left(
\exists\,r:
|\tilde{\Delta}_{\mathrm{off},r}-\mu_r|
>\frac{\mu_r/2+t}{C}
\right)
\le
\mathbb P(\mathcal E_{\tilde A}^{\,c})
+
\sum_{r=1}^q
\mathbb P\left(
\left\{
|\tilde{\Delta}_{\mathrm{off},r}-\mu_r|
>\frac{\mu_r/2+t}{C}
\right\}
\cap\mathcal E_{\tilde A}
\right)
\le
2e^{-\frac{m+n}{2}}+4q e^{-\frac{ct}{b_{m,n}}}.
\end{aligned}
\]
Combining this with \eqref{eq:app-decoupling-off}
and recalling that
$b_{m,n}=32x_{\max}^4(m+n)$ gives
\[
\begin{aligned}
&\mathbb P\left(
\exists\,r:
\Delta_{\mathrm{off},r}<\frac{\mu_r}{2}-t
\right)
\le
2C e^{-\frac{m+n}{2}}
+4Cq e^{-\frac{ct}{b_{m,n}}}
\le
C\exp\bigl(-c(m+n)\bigr)
+Cq\exp\left(
-\frac{ct}{x_{\max}^4(m+n)}
\right),
\end{aligned}
\]
where $c,C>0$ are absolute constants. Since the full squared HS norm in \eqref{eq:app-expand-HSnorm} is at least $\Delta_{\mathrm{off},r}$ for every $r$, this proves \eqref{eq:HS-lower-finite-family}.

\end{proof}

\subsection{Proof of Lemma~\ref{lem:quartic-y-upbound-1}}\label{PF:quartic-y-upbound-1}

We first control the fluctuation $|\mathbb G(\tilde{\vtheta})|$ at a fixed net point and then apply a union bound to obtain a simultaneous bound over $\mathcal N$.

\paragraph{Pointwise concentration.}
For each $\vtheta\in\Theta$, define
\[
\mM(\vtheta) := \mSigma_o^{-1/2}\bigl(\mSigma(\vtheta)-\mSigma_o\bigr) \mSigma_o^{-1/2}.
\]
Fix a deterministic point $\tilde{\vtheta}\in\mathcal N$, let $\tilde{\mSigma}:=\mSigma(\tilde{\vtheta})$. By the definition of $\mathbb G$ in \eqref{eq:def-G}, we have
\[
|\mathbb G(\tilde{\vtheta})|
=
\left|
\frac{1}{L}\sum_{\ell=1}^L
\vy_\ell^H\bigl(\tilde{\mSigma}-\mSigma_o\bigr)\vy_\ell
-
\Tr\bigl((\tilde{\mSigma}-\mSigma_o)\mSigma_o^{-1}\bigr)
\right|.
\]
Conditional on $\mA$, the observations are independent with distribution $\mathcal{CN}(0,\mSigma_o^{-1})$, and $\tilde{\mSigma}-\mSigma_o$ is a deterministic Hermitian matrix. In particular,
\[
\E\!\left[
\vy_\ell^H \bigl(\tilde{\mSigma}-\mSigma_o\bigr)
\vy_\ell \,\middle|\,\mA
\right]
=
\Tr\bigl( (\tilde{\mSigma}-\mSigma_o)\mSigma_o^{-1} \bigr).
\]
Thus, the quantity inside the absolute value is an average of independent centered quadratic forms. If $\mM(\tilde\vtheta)=0$, then $\mathbb G(\tilde\vtheta)=0$. Otherwise, we can use the Lemma~\ref{lem:quad-yt}, which gives, for every $t>0$,
\begin{equation}
\mathbb P\left( |\mathbb G(\tilde\vtheta)|>t \,\middle|\,\mA \right)
\le
2\exp\left( -c\min\left\{
\frac{Lt^2}{\|\mM(\tilde\vtheta)\|_{HS}^2},
\frac{Lt}{\|\mM(\tilde\vtheta)\|_2}
\right\} \right),
\label{eq:empirical error upper bound}
\end{equation}
where $c>0$ is an absolute constant.

\paragraph{The choice of threshold.}
To remove the operator norm from the final bound, we first bound $\|\mM(\tilde\vtheta)\|_2$ uniformly over the parameter class. By the definitions of $\Lambda_{\min}$ and $\Lambda_{\max}$ in \eqref{eq:def-Lambda-of-Sigma},
\[
\Lambda_{\min}\mI
\preceq
\tilde{\mSigma}^{-1}
\preceq
\Lambda_{\max}\mI
\qquad
\text{for every }\vtheta\in\Theta.
\]
Since $\vtheta_o\in\Theta$, these inequalities imply
\[
\frac{\Lambda_{\min}}{\Lambda_{\max}}\mI
\preceq
\mI+\mM(\tilde\vtheta)
=
\mSigma_o^{-1/2}\tilde{\mSigma}\mSigma_o^{-1/2}
\preceq
\frac{\Lambda_{\max}}{\Lambda_{\min}}\mI.
\]
Consequently, using $\Lambda_{\max}/\Lambda_{\min}\ge1$ and the definition of $C_{\sigma,m,n}$ in \eqref{eq:def-Lambda-of-Sigma}, we obtain
\begin{equation}
\|\mM(\tilde\vtheta)\|_2
\le
\frac{\Lambda_{\max}}{\Lambda_{\min}}-1
\le
C_{\sigma,m,n}
\qquad
\text{for every }\tilde\vtheta\in\Theta.
\label{eq:M-uniform-spectral-bound}
\end{equation}

Recall from \eqref{eq:number-of-logN} that
\[
\log|\mathcal N| \le \kappa_\delta
:=
k_x\log\left(\frac{a_xn^{b_x}}{\delta_x}\right)
+ k_\phi\log\left(\frac{a_\phi m^{b_\phi}}{\delta_\phi}\right).
\]
To obtain a simultaneous bound over the net, we choose the threshold so that both terms in the concentration exponent are at least a sufficiently large multiple of $\kappa_\delta$.
Set
\[
\mathcal R_\delta
:=
\sqrt{\frac{c_1\kappa_\delta}{L}},
\qquad
t(\vtheta)
:=
\mathcal R_\delta\|\mM(\tilde\vtheta)\|_{HS}
+
C_{\sigma,m,n}\mathcal R_\delta^2,
\]
where $c_1>0$ will be chosen sufficiently large.
It suffices to consider $\kappa_\delta>0$, since the claimed probability lower bound is trivial when $\kappa_\delta=0$. Then we get
\[
\frac{Lt(\tilde\vtheta)^2}{\|\mM(\tilde\vtheta)\|_{HS}^2}
\ge L\mathcal R_\delta^2 = c_1\kappa_\delta, 
\quad
\frac{Lt(\tilde\vtheta)}{\|\mM(\tilde\vtheta)\|_2}
\ge \frac{LC_{\sigma,m,n}\mathcal R_\delta^2}{\|\mM(\tilde\vtheta)\|_2}
\ge L\mathcal R_\delta^2 = c_1\kappa_\delta.
\]
Substituting these inequalities into \eqref{eq:empirical error upper bound} yields
\begin{equation}
\mathbb P\left(
|\mathbb G(\tilde\vtheta)|
>
\mathcal R_\delta\|\mM(\tilde\vtheta)\|_{HS}
+
C_{\sigma,m,n}\mathcal R_\delta^2
\,\middle|\,\mA
\right)
\le
2e^{-cc_1\kappa_\delta}.
\label{eq:empirical-R-2net}
\end{equation}
The same bound holds when $\mM(\tilde\vtheta)=0$, because the failure event is then empty.

\paragraph{Uniform control over the net.}
Applying a union bound to \eqref{eq:empirical-R-2net}, we obtain
\begin{align*}
&\mathbb P\left(
\exists\,\tilde\vtheta\in\mathcal N:
|\mathbb G(\tilde\vtheta)|
>
\mathcal R_\delta\|\mM(\tilde\vtheta)\|_{HS}
+ C_{\sigma,m,n}\mathcal R_\delta^2
\,\middle|\,\mA
\right) \\
\le&
\sum_{\vtheta\in\mathcal N}
\mathbb P\left(
|\mathbb G(\tilde\vtheta)|>t(\tilde\vtheta)
\,\middle|\,\mA
\right)
\le
2|\mathcal N|e^{-cc_1\kappa_\delta}
\le 2e^{-(cc_1-1)\kappa_\delta}.
\end{align*}
Choose $c_1$ sufficiently large that $c_0:=cc_1-1>0$. Then, conditional on $\mA$, with probability at least $1-2e^{-c_0\kappa_\delta}$,
\[
|\mathbb G(\tilde\vtheta)|
\le \mathcal R_\delta\|\mM(\tilde\vtheta)\|_{HS}
+
C_{\sigma,m,n}\mathcal R_\delta^2
\qquad \text{for every }\tilde\vtheta\in\mathcal N.
\]
Since this inequality holds simultaneously over the net, it also holds at the data-dependent nearest net point $\tilde\vtheta$. This completes the proof.

\subsection{Proof of Lemma~\ref{lem:quartic-y-upbound-2}}\label{PF:quartic-y-upbound-2}

We condition on $\mA$ throughout the proof and write $\hat{\mSigma}:=\mSigma(\hat{\vtheta})$ and $\tilde{\mSigma}:=\mSigma(\tilde{\vtheta})$.
By the the definition of $\mathbb G(\vtheta)$ in \eqref{eq:def-G} and triangle inequality, we have
\begin{align}
    |\mathbb{G}(\hat{\vtheta}) - \mathbb{G}(\tilde{\vtheta})|  
    &= \left| \frac{1}{L}\sum_{\ell=1}^L\vy_\ell^H(\hat{\mSigma} - \tilde{\mSigma})\vy_\ell - \Tr((\hat{\mSigma} - \tilde{\mSigma})\mSigma_o^{-1}) \right|  \nonumber\\
    &\le \left| \frac{1}{L}\sum_{l=1}^{L} \vy_l^H (\hat{\mSigma} - \tilde{\mSigma}) \vy_l \right|
    + \left| \text{Tr}((\hat{\mSigma} - \tilde{\mSigma})\mSigma_o^{-1}) \right|.
    \label{eqn: lem error Ghat Gtilde 1}
\end{align}
We first control the trace term, using the inequality \(\Tr(M)\le m\|M\|_2\), we have
\begin{align}
\left| \Tr((\hat{\mSigma} - \tilde{\mSigma})\mSigma_o^{-1}) \right|
\le
m\Bigl\|(\hat{\mSigma} - \tilde{\mSigma})\mSigma_o^{-1}\Bigr\|_2
\le
m\lambda_{\max}(\mSigma_o^{-1})\,\|\hat{\mSigma} - \tilde{\mSigma}\|_2.
\label{eqn: lem error Ghat Gtilde 2}
\end{align}
Next we bound the empirical quadratic term. 
Conditional on $\mA$, the observations are independent with distribution $\mathcal{CN}(\vec{0},\mSigma_o^{-1})$. Define $\vg_\ell:=\mSigma_o^{1/2}\vy_\ell$.
Then, conditional on $\mA$, $\vg_1,\ldots,\vg_L$ are i.i.d.\ $\mathcal{CN}(\vec{0},\mI_m)$, and $\vy_\ell=\mSigma_o^{-1/2}\vg_\ell$. Therefore,
\begin{align}
    \left| \frac{1}{L}\sum_{l=1}^{L} \vy_l^H (\hat{\mSigma} - \tilde{\mSigma}) \vy_l \right| 
    \le \|\mSigma^{-\frac{1}{2}}_o(\hat{\mSigma} - \tilde{\mSigma})\mSigma^{-\frac{1}{2}}_o\|_2 \cdot \frac{1}{L}\sum_{l=1}^{L} \|\vg_l\|_2^2  
    \le  \lambda_{\max}(\mSigma_o^{-1})\,\|\hat{\mSigma} - \tilde{\mSigma}\|_2 \cdot \frac{1}{L}\sum_{l=1}^{L} \|\vg_l\|_2^2.   \label{eqn: lem error Ghat Gtilde 3}
\end{align}
By Lemma~\ref{lem: bound of chi2}, taking $t=m$, there exists an absolute constant $c_2>0$ such that, with probability at least $1-\exp(-c_2mL)$,
\[
\frac{1}{L}\sum_{\ell=1}^L \|\vg_\ell\|_2^2 \le 2m.
\]
Substituting this into \eqref{eqn: lem error Ghat Gtilde 3} yields, on the same event,
\begin{align}
\left| \frac{1}{L}\sum_{l=1}^{L} \vy_l^H (\hat{\mSigma} - \tilde{\mSigma}) \vy_l \right| 
\le
2m\,\lambda_{\max}(\mSigma_o^{-1})\|\hat\mSigma-\tilde\mSigma\|_2.
\label{eqn: lem error Ghat Gtilde 4}
\end{align}
Combining \eqref{eqn: lem error Ghat Gtilde 1}, \eqref{eqn: lem error Ghat Gtilde 2}, and \eqref{eqn: lem error Ghat Gtilde 4}, we conclude that, with probability at least \(1-\exp(-c_2mL)\),
\[
|\mathbb{G}(\hat{\vtheta}) - \mathbb{G}(\tilde{\vtheta})| 
\le
3m\,\lambda_{\max}(\mSigma_o^{-1})\|\hat\mSigma-\tilde\mSigma\|_2.
\]

It remains to bound the precision-matrix discrepancy. Applying \eqref{eqn:lem preci error upbound 1} in Lemma~\ref{lem: cov trace and 2 upbound} to the pair $(\hat{\vtheta},\tilde{\vtheta})$ gives
\begin{align}
\|\hat\mSigma-\tilde\mSigma\|_2
\le
\frac{E_{\max}}{\Lambda_{\min}^2}
\left(
\|\hat{\vx}^2-\tilde{\vx}^2\|_\infty
+
2x_{\max}^2\|\hat\mPhi-\tilde\mPhi\|_2
\right).\label{eqn: lem error Ghat Gtilde 5}
\end{align}
Since $\hat{\vx}$ and $\tilde{\vx}$ have entries bounded by $x_{\max}$, the covering property implies
\begin{align}\label{eq:bound-x2-by-delta}
\|\hat{\vx}^2-\tilde{\vx}^2\|_\infty
\le
2x_{\max}\|\hat{\vx}-\tilde{\vx}\|_\infty
\le
2x_{\max}\|\hat{\vx}-\tilde{\vx}\|_2
\le
2x_{\max}\delta_x.
\end{align}
Similarly, using $|e^{ia}-e^{ib}|\le |a-b|$ for real $a,b$, we obtain
\begin{align}
\label{eq:bound-Phi-by-delta}
\|\hat{\mPhi}-\tilde{\mPhi}\|_2
=
\max_{1\le j\le m}
|e^{i\hat{\phi}_j}-e^{i\tilde{\phi}_j}|
\le
\|\hat{\bm{\phi}}-\tilde{\bm{\phi}}\|_2
\le
\delta_\phi.
\end{align}
Substituting these estimates into \eqref{eqn: lem error Ghat Gtilde 5} and using $\lambda_{\max}(\mSigma_o^{-1})\le\Lambda_{\max}$, we conclude that, with probability at least \(1-\exp(-c_2mL)\),
\begin{align*}
|\mathbb G(\hat{\vtheta})-\mathbb G(\tilde{\vtheta})|
&\le
C_2\frac{\Lambda_{\max}E_{\max}}{\Lambda_{\min}^2}
m(\delta_x+\delta_\phi)
\le
C_2C_{\sigma,m,n} m(\delta_x+\delta_\phi),
\end{align*}
where $C_2:=6\max\{x_{\max},x_{\max}^2\}$ depends only on $x_{\max}$ and we use $\Lambda_{\max}E_{\max}/\Lambda_{\min}^2\le C_{\sigma,m,n}$ from \eqref{eqn:def-C-sigma-m-n}. This completes the proof.

\subsection{Proof of Lemma~\ref{lem:population-gap-lower-bound}}
\label{PF:lowBd-fbar-trace}

We decompose the expected likelihood difference as
\begin{align}\label{eqn:decompose fbar hat minus o}
\overline{\mathcal L}(\hat{\vtheta})-\overline{\mathcal L}(\vtheta_o)
&=
\bigl(\overline{\mathcal L}(\tilde{\vtheta})-\overline{\mathcal L}(\vtheta_o)\bigr)
+
\bigl(\overline{\mathcal L}(\hat{\vtheta})-\overline{\mathcal L}(\tilde{\vtheta})\bigr).
\end{align}
We first lower bound $\overline{\mathcal L}(\tilde{\vtheta}) -\overline{\mathcal L}(\vtheta_o)$. Let $t_1,\ldots,t_m$ be the eigenvalues of $\mSigma_o^{-1/2}\tilde{\mSigma}\mSigma_o^{-1/2}$.
By the definition of expected likelihood in \eqref{def:exp-liklihood},
\begin{align}
\overline{\mathcal L}(\tilde{\vtheta})
-\overline{\mathcal L}(\vtheta_o)
&=
\Tr\left(\mSigma_o^{-1/2}\tilde\mSigma\mSigma_o^{-1/2}\right) 
-\log\det\left(\mSigma_o^{-1/2}\tilde\mSigma\mSigma_o^{-1/2}\right) - m 
=\sum_{i=1}^m(t_i-\log t_i-1).
\label{eqn:population-gap-T}
\end{align}
The definition of $\Lambda_{\max},\Lambda_{\min}$, $C_{\sigma,m,n}$ and the spectral bounds in \eqref{eq:def-Lambda-of-Sigma} imply
\[
C_{\sigma,m,n}^{-1}
\leq \frac{\Lambda_{\min}}{\Lambda_{\max}}
\leq t_i
\leq \frac{\Lambda_{\max}}{\Lambda_{\min}}
\leq C_{\sigma,m,n}.
\]
For every $t\in[C_{\sigma,m,n}^{-1},C_{\sigma,m,n}]$, since $0<1+s(t-1)\leq C_{\sigma,m,n}$ for $s\in[0,1]$, we get
\begin{equation}
t-\log t-1
=(t-1)^2\int_0^1\frac{s}{1+s(t-1)}\,ds
\geq \frac{(t-1)^2}{2C_{\sigma,m,n}}.
\label{eqn:scalar-log-lower-bound}
\end{equation}
Applying this inequality and using cyclicity of trace gives
\begin{align}
\overline{\mathcal L}(\tilde{\vtheta})
-\overline{\mathcal L}(\vtheta_o)
&\geq \frac{1}{2C_{\sigma,m,n}}
\sum_{i=1}^m(t_i-1)^2
=\frac{1}{2C_{\sigma,m,n}}
\Tr\!\left[
\mSigma_o^{-1}(\tilde{\mSigma}-\mSigma_o)
\mSigma_o^{-1}(\tilde{\mSigma}-\mSigma_o)
\right].
\label{eqn:net-population-gap-lower}
\end{align}

For the second term, by \eqref{def:exp-liklihood},
\begin{align}
\overline{\mathcal L}(\hat{\vtheta})
-\overline{\mathcal L}(\tilde{\vtheta})
&=
\Tr\!\left[
(\hat{\mSigma}-\tilde{\mSigma})\tilde{\mSigma}^{-1}
\right]
-\log\det\!\left(
\tilde{\mSigma}^{-1/2}\hat{\mSigma}\tilde{\mSigma}^{-1/2}
\right)
+
\Tr\!\left[
(\hat{\mSigma}-\tilde{\mSigma})
(\mSigma_o^{-1}-\tilde{\mSigma}^{-1})
\right].
\label{eqn:population-gap-net-decomposition}
\end{align}
The matrix inequality established in \eqref{eqn:scalar-log-lower-bound} gives
\[
\Tr\!\left[
(\hat{\mSigma}-\tilde{\mSigma})\tilde{\mSigma}^{-1}
\right]
-\log\det\!\left(
\tilde{\mSigma}^{-1/2}
\hat{\mSigma}
\tilde{\mSigma}^{-1/2}
\right)
\geq 0.
\]
Hence, \eqref{eqn:population-gap-net-decomposition} can be lower bounded by
\begin{align*}
\overline{\mathcal L}(\hat{\vtheta})
-\overline{\mathcal L}(\tilde{\vtheta})
\ge
\Tr\!\left[
(\hat{\mSigma}-\tilde{\mSigma})
(\mSigma_o^{-1}-\tilde{\mSigma}^{-1})
\right]
\geq
-m\|\hat{\mSigma}-\tilde{\mSigma}\|_2
\|\mSigma_o^{-1}-\tilde{\mSigma}^{-1}\|_2
\geq
-2m\Lambda_{\max}
\|\hat{\mSigma}-\tilde{\mSigma}\|_2.
\end{align*}
By Lemma~\ref{lem: cov trace and 2 upbound} and the net bounds, we have
\[
\|\hat{\mSigma}-\tilde{\mSigma}\|_2
\leq
\frac{E_{\max}}{\Lambda_{\min}^2}
\left(2x_{\max}\delta_x+2x_{\max}^2\delta_\phi\right).
\]
Consequently, using the definition of $\Lambda_{\max},\Lambda_{\min}$ and $C_{\sigma,m,n}$ in \eqref{eq:def-Lambda-of-Sigma}, we get
\begin{align}
\overline{\mathcal L}(\hat{\vtheta})
-\overline{\mathcal L}(\tilde{\vtheta})
\geq
-\frac{4m\Lambda_{\max}E_{\max}}{\Lambda_{\min}^2}
\left(x_{\max}\delta_x+x_{\max}^2\delta_\phi\right)
\geq
-C C_{\sigma,m,n}m(\delta_x+\delta_\phi),
\label{eqn:population-gap-net-logdet}
\end{align}
where $C>0$ depends only on $x_{\max},x_{\min}$.

Combining this bound with \eqref{eqn:net-population-gap-lower} yields
\begin{align*}
\overline{\mathcal L}(\hat{\vtheta})
-\overline{\mathcal L}(\vtheta_o)
&=
\overline{\mathcal L}(\tilde{\vtheta})
-\overline{\mathcal L}(\vtheta_o)
+\overline{\mathcal L}(\hat{\vtheta})
-\overline{\mathcal L}(\tilde{\vtheta})\\
&\geq
\frac{1}{2C_{\sigma,m,n}}
\Tr\!\left[
\mSigma_o^{-1}(\tilde{\mSigma}-\mSigma_o)
\mSigma_o^{-1}(\tilde{\mSigma}-\mSigma_o)
\right]
- C C_{\sigma,m,n}m(\delta_x+\delta_\phi).
\end{align*}
This completes the proof.

\subsection{Proof of Lemma~\ref{lem:lowerbound-HSnorm}}\label{PF:lowerbound-HSnorm}
\begin{proof}
Let $\tilde{\mU}:=\tilde{\mPhi}^H\mPhi_o$, by unitary invariance of the HS norm, we have
\begin{align*}
\left\|
\tilde{\mPhi}\mA\tilde{\mX}^2\mA^H\tilde{\mPhi}^H - \mPhi_o\mA\mX_o^2\mA^H\mPhi_o^H
\right\|_{HS}^2
=
\left\|
\mA\tilde{\mX}^2\mA^H - \tilde{\mU}\mA\mX_o^2\mA^H\tilde{\mU}^H
\right\|_{HS}^2.
\end{align*}
For each $\tilde\vtheta\in\mathcal N$, define
\[
\mu(\tilde\vtheta)
:=
m(m-1)\|\tilde{\mX}^2-\mX_o^2\|_{HS}^2
+
2\bigl(m^2-|\Tr(\tilde{\mU})|^2\bigr) \Tr(\tilde{\mX}^2\mX_o^2).
\]
Then applying Lemma~\ref{lem:HS_lower_bound} to the finite family
$\{(\tilde{\mX},\tilde{\mU}):
\tilde\vtheta\in\mathcal N\}$ gives, for every $t>0$,
\begin{align}\label{eq:app-bound-HS-applied}
&\mathbb P\left( \exists\,\tilde\vtheta\in\mathcal N:
\left\| \mA\tilde{\mX}^2\mA^H
-
\tilde{\mU}\mA\mX_o^2\mA^H\tilde{\mU}^H\right\|_{HS}^2
< \frac{\mu(\tilde\vtheta)}{2}-t
\right)
\le
C|\mathcal N|\exp\left(-\frac{ct}{x_{\max}^4(m+n)}\right) 
+ C\exp\bigl(-c(m+n)\bigr).
\end{align}
We choose
\[
t=C_0x_{\max}^4(m+n)\kappa_\delta,
\]
where $C_0>0$ is a sufficiently large absolute constant. Since $|\mathcal N|\le e^{\kappa_\delta}$, taking $C_0$ large enough yields
\[
Ce^{\kappa_\delta-cC_0\kappa_\delta} + Ce^{-c(m+n)}
\le Ce^{-\kappa_\delta}+Ce^{-c(m+n)}.
\]
Therefore, with probability at least $1-Ce^{-\kappa_\delta}-Ce^{-c(m+n)}$, simultaneously for all $\tilde\vtheta\in\mathcal N$,
\[
\left\|
\tilde{\mPhi}\mA\tilde{\mX}^2\mA^H\tilde{\mPhi}^H
-
\mPhi_o\mA\mX_o^2\mA^H\mPhi_o^H
\right\|_{HS}^2
\ge
\frac12\mu(\tilde\vtheta)
-
C_0x_{\max}^4(m+n)\kappa_\delta.
\]

It remains to lower bound $\mu(\tilde\vtheta)$. Since $\tilde\vx,\vx_o\in[x_{\min},x_{\max}]^n$ and $m\ge2$,
\[
\|\tilde{\mX}^2-\mX_o^2\|_{HS}^2
\ge
4x_{\min}^2\|\tilde\vx-\vx_o\|_2^2,
\qquad 
\Tr(\tilde{\mX}^2\mX_o^2) \ge nx_{\min}^4,
\]
and
\[
m(m-1)\ge \frac{m^2}{2},
\qquad m^2-|\Tr(\tilde{\mU})|^2
\ge m\bigl(m-|\Tr(\tilde{\mU})|\bigr).
\]
Therefore,
\begin{align*}
\frac12\mu(\tilde\vtheta)
\ge x_{\min}^2m^2\|\tilde\vx-\vx_o\|_2^2
+ x_{\min}^4nm \bigl(m-|\Tr(\tilde{\mU})|\bigr)
\ge \min\{x_{\min}^2,x_{\min}^4\}\, E(\tilde\vtheta).
\end{align*}
\end{proof}

\subsection{Proof of Lemma~\ref{lem:noise-comparison}}
\label{app:noise-comparison}

\begin{proof}
The case $s_1=s_2$ is immediate by taking $\widehat\vtheta_1=\widehat\vtheta_2$. Assume $s_1<s_2$. Write $\ell(a,b):=d_{\rm joint}^2(a,b)$ for $a,b\in\Theta$. This loss is symmetric and satisfies
\begin{equation}
\ell(a,c)\le2\ell(a,b)+2\ell(b,c), \qquad a,b,c\in\Theta.
\label{eq:noise-comparison-triangle}
\end{equation}
The signal component satisfies this inequality by $\|u+v\|_2^2\le2\|u\|_2^2+2\|v\|_2^2$.
For the phase component, let $\mPhi_a,\mPhi_b,\mPhi_c$ denote the corresponding phase matrices, and let $\alpha$ and $\beta$ be optimal alignment angles for $(\mPhi_a,\mPhi_b)$ and $(\mPhi_b,\mPhi_c)$, respectively. Then we have
\begin{align*}
d_\Phi^2(\mPhi_a,\mPhi_c)
&\le\|\mPhi_a-e^{i(\alpha+\beta)}\mPhi_c\|_{\rm HS}^2
=\|(\mPhi_a-e^{i\alpha}\mPhi_b)
 +e^{i\alpha}(\mPhi_b-e^{i\beta}\mPhi_c)\|_{\rm HS}^2
\le2d_\Phi^2(\mPhi_a,\mPhi_b)
 +2d_\Phi^2(\mPhi_b,\mPhi_c).
\end{align*}
Combining the normalized amplitude and phase components proves \eqref{eq:noise-comparison-triangle}.

Stack the observations at variance $s_1$ into $\mY\in\mathbb C^{mL}$, and define
\[
\vec{\eta}\sim\mathcal{CN}\bigl(0,(s_2-s_1)\mI_{mL}\bigr)
\]
be independent of $(\mA,\mY)$. For each observed $(\mA,\mY)$, choose
\begin{equation}
\widehat\vtheta_1(\mA,\mY)
\in\operatorname*{arg\,min}_{\vtheta\in\Theta}
\E_{\vec{\eta}}\!\left[ \ell\bigl(\vtheta,\hat\vtheta_2(\mA,\mY+\vec{\eta})\bigr) \right].
\label{eq:noise-comparison-estimator}
\end{equation}
Here the expectation is only over the known distribution of $\vec{\eta}$, with $(\mA,\mY)$ fixed. Since $\widehat\vtheta_2$ is measurable and $\ell$ is bounded and continuous, the expected loss in \eqref{eq:noise-comparison-estimator} is measurable in $(\mA,\mY)$ and continuous in $\vtheta$.
Since $\Theta$ is a non-empty compact metric space, the measurable maximum theorem~\cite[Theorem~18.19]{aliprantis2006infinite} guarantees the existence of a measurable minimizer in \eqref{eq:noise-comparison-estimator}.
Fix one such selection to define $\widehat\vtheta_1$. The resulting estimator depends only on $(\mA,\mY)$ and the fixed measurable estimator $\widehat\vtheta_2$.

Fix $(\mA,\mY)$ and $\vtheta_o\in\Theta$, by \eqref{eq:noise-comparison-triangle} and symmetry of the loss, we obtain
\begin{align*}
\ell(\widehat\vtheta_1,\vtheta_o)
&\le2\E_\eta\ell(\widehat\vtheta_1,\hat\vtheta_2(\mA,\mY+\vec{\eta}))
 +2\E_\eta\ell(\hat\vtheta_2(\mA,\mY+\vec{\eta}),\vtheta_o)\\
&\le2\E_\eta\ell(\vtheta_o,\hat\vtheta_2(\mA,\mY+\vec{\eta}))
 +2\E_\eta\ell(\hat\vtheta_2(\mA,\mY+\vec{\eta}),\vtheta_o)\\
&=4\E_\eta\ell(\hat\vtheta_2(\mA,\mY+\vec{\eta}),\vtheta_o).
\end{align*}
The second inequality follows from the minimizing property in \eqref{eq:noise-comparison-estimator}.

Conditional on $\mA$, adding $\vec{\eta}$ increases the additive-noise variance from $s_1$ to $s_2$. Thus $(\mA,\mY+\vec{\eta})$ under variance $s_1$ has the same joint law as $(\mA,\mY)$ under variance $s_2$. Taking expectation over $(\mA,\mY)$ therefore gives
\begin{align*}
\E_{\vtheta_o,s_1}
\bigl[\ell(\widehat\vtheta_1,\vtheta_o)\bigr]
\le4\,\E_{\vtheta_o,s_1}\E_\eta
\bigl[\ell(\widehat\vtheta_2(\mA,\mY+\vec{\eta}),\vtheta_o)\bigr]
=4\,\E_{\vtheta_o,s_2}
\bigl[\ell(\widehat\vtheta_2,\vtheta_o)\bigr],
\end{align*}
which proves \eqref{eq:noise-comparison-risk}.
\end{proof}

\section{Proof of Lower Bound Results}

\begin{lemma}[Generalized Fano lemma,{\cite[Lemma~3]{yu1997assouad}}]
\label{lem:generalized_fano}
Let $\mathcal P$ be a family of probability measures, each associated with a parameter $\theta(P)$ in a parameter space $\mathcal T$. Let $d$ be a pseudo-metric on $\mathcal T$. Suppose there exist an integer $r\ge2$, distributions $P_1,\ldots,P_r\in\mathcal P$, and $\alpha_r,\beta_r>0$ such that
\begin{enumerate}
\item[(i)]
The parameters $\{\theta(P_1),\ldots,\theta(P_r)\}$ form an $\alpha_r$-separated subset of $(\mathcal T,d)$, namely,
\[
d(\theta(P_i),\theta(P_j))\ge\alpha_r, \qquad 1\le i\ne j\le r.
\]
\item[(ii)]
The pairwise KL divergences satisfy
\[
D_{\mathrm{KL}}(P_i\|P_j)\le\beta_r, \qquad 1\le i\ne j\le r.
\]
\end{enumerate}
Then, for any measurable $\mathcal T$-valued estimator $\widehat\theta$,
\[
\max_{1\le j\le r}
\mathbb E_{P_j}\left[d(\widehat\theta,\theta(P_j))\right]
\ge
\frac{\alpha_r}{2}\left(1-\frac{\beta_r+\log2}{\log r}\right).
\]
\end{lemma}

\subsection{Proof of Lemma~\ref{lem:signal-oracle-lb}}\label{pf:lem:signal-oracle-lb}
\begin{proof}
Let $R_{\mathbb R}(\mathcal C_x;M,n,J,\sigma_z)$ denote the minimax risk in the real-valued fixed-sensing-matrix model of \cite{xing2025minimax}, with $M$ rows, $J$ looks, and loss $n^{-1}\|\widehat \vx-\vx_o\|_2^2$.
We compare this risk with $R_X$ by constructing the complex observations from a real experiment with $M=2m$ and $J=2L$.

Let
\[
 \mB=\begin{pmatrix}\mG_1\\ \mG_2\end{pmatrix}\in\mathbb R^{2m\times n}
\]
have independent $N(0,1)$ entries, where $\mG_1,\mG_2\in\mathbb R^{m\times n}$. For each $\ell=1,\ldots,L$, consider the two real observations
\[
 \mU_\ell=\mB\mX_o\vu_\ell+\vec{\xi}_\ell,\qquad
 \mV_\ell=\mB\mX_o\vv_\ell+\veta_\ell,
\]
where $\vu_\ell,\vv_\ell\sim N(0,I_n)$ and $\vec{\xi}_\ell,\veta_\ell\sim N(0,\sigma_z^2I_{2m})$.
All these vectors are mutually independent and independent of $\mB$. Thus $(\mU_\ell,\mV_\ell)_{\ell=1}^L$ constitutes $2L$ looks sharing the same real sensing matrix $\mB$. For any vector $\vq\in\mathbb R^{2m}$, write $\vq^{(1)}$ and $\vq^{(2)}$ for its first and last $m$ coordinates, respectively. Define
\[
 \mA=\frac{\mG_1+i\mG_2}{\sqrt2},\qquad
 \widetilde \vy_\ell= \frac{\mU_\ell^{(1)}-\mV_\ell^{(2)}
       +i\bigl(\mU_\ell^{(2)}+\mV_\ell^{(1)}\bigr)}2.
\]
A direct expansion gives
\[
 \widetilde \vy_\ell=\mA\mX_o\vw_\ell+\widetilde \vz_\ell,
 \qquad \vw_\ell=\frac{\vu_\ell+i\vv_\ell}{\sqrt2},
 \qquad \widetilde \vz_\ell=
 \frac{\vec{\xi}_\ell^{(1)}-\veta_\ell^{(2)}+i\bigl(\vec{\xi}_\ell^{(2)}+\veta_\ell^{(1)}\bigr)}2.
\]
The entries of $\mA$ are independent $\CN(0,1)$ variables, and $\vw_\ell\sim\CN(0,I_n)$. The real and imaginary parts of $\widetilde \vz_\ell$ are independent $N(0,\sigma_z^2I_m/2)$ vectors, so $\widetilde \vz_\ell\sim\CN(0,\sigma_z^2I_m)$.
Moreover, all the speckle and additive-noise vectors are mutually independent and independent of $\mA$. Hence the transformed data have exactly the joint distribution specified in~\eqref{eq:model-oracle-signal}, including when $\sigma_z=0$.

The map from $(\mB,(\mU_\ell,\mV_\ell)_{\ell=1}^L)$ to $(\mA,(\widetilde \vy_\ell)_{\ell=1}^L)$ is deterministic and does not depend on $\vx_o$. Therefore, every estimator in the complex model induces an estimator in the real model with the same risk at every $\vx_o\in\mathcal C_x$. Taking the supremum over $\vx_o$ and then the infimum over estimators gives
\begin{equation}
 R_X(\mathcal C_x;m,n,L,\sigma_z)
 \ge R_{\mathbb R}(\mathcal C_x;2m,n,2L,\sigma_z).
 \label{eq:l6-risk-comparison}
\end{equation}

It remains to apply \cite[Theorem~2.15]{xing2025minimax} to the right-hand side. The rate assumption in Lemma~\ref{lem:signal-oracle-lb} implies
\[
 \max\{\sigma_z^4,4m^2,n^2\}k_x\log n
 \le4\max\{\sigma_z^4,m^2,n^2\}k_x\log n
 \le4m^2n^{1-\varepsilon}L
 \le(2m)^2n^{1-\varepsilon}(2L).
\]
All the remaining assumptions of that theorem are preserved under
$m\mapsto2m$ and $L\mapsto2L$, and the admissible signal-class family
is unchanged. Writing $\mathcal F_x=\mathcal F_x(a_x,b_x,k_x,n)$,
we conclude from \eqref{eq:l6-risk-comparison} that
\begin{align*}
 \sup_{\mathcal C_x\in\mathcal F_x}
 R_X(\mathcal C_x;m,n,L,\sigma_z)
 &\ge\sup_{\mathcal C_x\in\mathcal F_x}
 R_{\mathbb R}(\mathcal C_x;2m,n,2L,\sigma_z) \\
 &\ge c\,\frac{\max\{\sigma_z^4,4m^2,n^2\}k_x\log n}{(2m)^2n(2L)} \\
 &\ge\frac c8\,\frac{\max\{\sigma_z^4,m^2,n^2\}k_x\log n}{m^2nL},
\end{align*}
where $c>0$ depends only on $\varepsilon,x_{\min},x_{\max}$.
This proves the claimed complex-valued extension.
\end{proof}

\subsection{Proof of Lemma~\ref{lem:phase_class_construction}}\label{pf:lem:phase_class_construction}

\begin{proof}
We first construct a compact phase class satisfying the polynomial-complexity condition of order $k_\phi$. Let
\[
h:=\left\lfloor\frac{m-1}{N_{\mathrm{div}}^\phi}\right\rfloor,
\qquad
I_\ell:=\{2+(\ell-1)h,\ldots,1+\ell h\},
\quad \ell=1,\ldots,N_{\mathrm{div}}^\phi.
\]
Here $h\asymp m/N_{\mathrm{div}}^\phi$ is the block size, and the sets $I_\ell$ form disjoint consecutive blocks within $\{2,\ldots,m\}$. Let
\[
\mathcal S_{\mathrm{fin}}
:=
\left\{ S\subseteq\{1,\ldots,N_{\mathrm{div}}^\phi\}:|S|=k_\phi \right\}
\]
be the family of supports selecting exactly $k_\phi$ blocks. For each $S\in\mathcal S_{\mathrm{fin}}$, define the corresponding binary vector $\vec b^{\,S}\in\{0,1\}^m$ by
\[
b_j^S :=
\begin{cases}
1, & j\in\displaystyle\bigcup_{\ell\in S} I_\ell,\\
0, & \text{otherwise},
\end{cases}
\qquad j=1,\ldots,m.
\]
We then construct the phase class
\[
\Gamma_\phi
:=
\bigcup_{S\in\mathcal S_{\mathrm{fin}}}
\left\{
t\vec b^{\,S}:t\in[0,\pi/8]
\right\}.
\]
Thus, each phase vector has exactly $k_\phi$ active blocks sharing a common angle $t$, with all remaining coordinates set to zero. By construction, $\Gamma_\phi$ is compact and satisfies the constraints $\vphi\in[0,\pi/8]^m$ and $ \phi_1=0$.

We next verify the polynomial-complexity condition. For any fixed $S\in\mathcal S_{\mathrm{fin}}$, the disjointness of the blocks implies $\|\vec b^{\,S}\|_2^2=k_\phi h$. Hence, for any $s,t\in[0,\pi/8]$,
\[
\|t\vec b^{\,S}-s\vec b^{\,S}\|_2 = \sqrt{k_\phi h}\,|t-s|.
\]
Since the set $\{t\vec b^{\,S}:t\in[0,\pi/8]\}$ is a line segment of length $\pi\sqrt{k_\phi h}/8$, an equally spaced grid on $[0,\pi/8]$ with spacing $\epsilon/\sqrt{k_\phi h}$ therefore provides an $\epsilon$-cover of $\{t\vec b^{\,S}:t\in[0,\pi/8]\}$, giving
\[
N_\epsilon\left(
\{t\vec b^{\,S}:t\in[0,\pi/8]\}
\right)
\le
1+\frac{\pi\sqrt{k_\phi h}}{8\epsilon}.
\]
Taking the union of these covers over all $\binom{N_{\mathrm{div}}^\phi}{k_\phi}$ possible supports, we obtain
\[
\begin{aligned}
N_\epsilon(\Gamma_\phi)
&\le
\binom{N_{\mathrm{div}}^\phi}{k_\phi}
\left(
1+\frac{\pi\sqrt{k_\phi h}}{8\epsilon}
\right)
\le
\left(
\frac{eN_{\mathrm{div}}^\phi}{k_\phi}
\right)^{k_\phi}
\left(
\frac{2\sqrt{k_\phi h}}{\epsilon}
\right)^{k_\phi}
\le
\left(\frac{2em}{\epsilon}\right)^{k_\phi}.
\end{aligned}
\]
Here we used the standard binomial bound $\binom{N}{k}\le (eN/k)^k$, together with $\frac{N_{\mathrm{div}}^\phi}{k_\phi}\sqrt{k_\phi h} =\frac{N_{\mathrm{div}}^\phi h}{\sqrt{k_\phi h}}\le m-1.$ Thus, $\Gamma_\phi$ satisfies the polynomial-complexity condition of order $k_\phi$.

We next use a greedy packing argument to construct a large family of well-separated supports. Let $\mathcal S_{\mathrm{sep}}\subseteq\mathcal S_{\mathrm{fin}}$ be a maximal family satisfying
\[
|S\setminus T|\ge k_\phi/2, \qquad \text{for all distinct }S,T\in\mathcal S_{\mathrm{sep}},
\]
where the set $S\setminus T$ consists of the block indices that belong to $S$ but not to $T$. Such a family can be obtained by a greedy selection procedure. By maximality, for every $S\in\mathcal S_{\mathrm{fin}}$, there exists $T\in\mathcal S_{\mathrm{sep}}$ such that
\[
|S\setminus T|<k_\phi/2.
\]
We then bound the cardinality $r:=|\mathcal S_{\mathrm{sep}}|$. For each $T\in\mathcal S_{\mathrm{sep}}$, define $B(T):= \left\{ S\in\mathcal S_{\mathrm{fin}}: |S\setminus T|<k_\phi/2 \right\}.$ Since all supports have cardinality $k_\phi$, a support $S$ with $|S\setminus T|=q$ can be obtained by removing $q$ elements from $T$ and adding $q$ elements from its complement. Therefore,
\[
\begin{aligned}
|B(T)|
&=
\sum_{q=0}^{\lceil k_\phi/2\rceil-1}
\binom{k_\phi}{q}\binom{N_{\mathrm{div}}^\phi-k_\phi}{q}
\le
2^{k_\phi}
\left(\frac{2eN_{\mathrm{div}}^\phi}{k_\phi}\right)^{k_\phi/2}.
\end{aligned}
\]
Here we used $\sum_{q=0}^{k_\phi}\binom{k_\phi}{q}=2^{k_\phi}$ and, for $1\le q<k_\phi/2$,
\[
\begin{aligned}
\binom{N_{\mathrm{div}}^\phi-k_\phi}{q}
&\le
\left(\frac{eN_{\mathrm{div}}^\phi}{q}\right)^q
\le
\left(\frac{2eN_{\mathrm{div}}^\phi}{k_\phi}\right)^{k_\phi/2},
\end{aligned}
\]
where the last inequality follows from the monotonicity of $t\mapsto t\log(eN_{\mathrm{div}}^\phi/t)$ on $(0,N_{\mathrm{div}}^\phi]$.

By maximality, the sets $B(T)$, $T\in\mathcal S_{\mathrm{sep}}$, cover $\mathcal S_{\mathrm{fin}}$. Define the cardinality $r:=|\mathcal S_{\mathrm{sep}}|$, we get
\[
\binom{N_{\mathrm{div}}^\phi}{k_\phi}
=
|\mathcal S_{\mathrm{fin}}|
\le \sum_{T\in\mathcal S_{\mathrm{sep}}}|B(T)|
\le r\max_{T\in\mathcal S_{\mathrm{sep}}}|B(T)|.
\]
Combining this with the standard binomial bounds $(p/s)^s\le\binom{p}{s}\le(ep/s)^s$ gives
\[
\left(\frac{N_{\mathrm{div}}^\phi}{8ek_\phi}\right)^{k_\phi/2}
\le r
\le
\left( \frac{eN_{\mathrm{div}}^\phi}{k_\phi} \right)^{k_\phi}.
\]
Since $N_{\mathrm{div}}^\phi/k_\phi\asymp m^\eta$, we conclude that $\log r=\Theta(k_\phi\log m)$, which proves (i).

We now use the separated supports to construct a phase packing. Write $\mathcal S_{\mathrm{sep}}=\{S_1,\ldots,S_r\}$. For any $\delta_r\in(0,\pi/8)$, define $\vec{\phi}^{(i)}:=\delta_r\vec b^{\,S_i}$ for $i=1,\ldots,r$, and let
\[
U_{\mathrm{sep}}
:=\{\vec{\phi}^{(1)},\ldots,\vec{\phi}^{(r)}\}
\subset\Gamma_\phi.
\]
For any distinct $i,j$, let $q:=|S_i\setminus S_j|=|S_j\setminus S_i|$. By the separation property of $\mathcal S_{\mathrm{sep}}$, we have $k_\phi/2\le q\le k_\phi$. Since each block contains $h$ coordinates, the phase differences between $\vec{\phi}^{(i)}$ and $\vec{\phi}^{(j)}$ are $\delta_r$ on $qh$ coordinates, $-\delta_r$ on another $qh$ coordinates, and zero elsewhere. Consequently,
\[
\begin{aligned}
\Tr(\mPhi_i^H\mPhi_j)
&=(m-2qh)+qh\,e^{\mathrm i\delta_r}+qh\,e^{-\mathrm i\delta_r}
=m-2qh(1-\cos\delta_r).
\end{aligned}
\]
Since $2qh\le m$ and $\delta_r\in(0,\pi/8)$, this quantity is positive. Therefore,
\[
\begin{aligned}
d_\Phi^2(\mPhi_i,\mPhi_j)
=2\bigl(m-|\Tr(\mPhi_i^H\mPhi_j)|\bigr)
=4qh(1-\cos\delta_r)
=\|\mPhi_i-\mPhi_j\|_F^2.
\end{aligned}
\]
Finally, using $k_\phi/2\le q\le k_\phi$, $h\asymp m/N_{\mathrm{div}}^\phi$, and $\delta_r^2/4 \le 1-\cos\delta_r \le \delta_r^2/2$ for $|\delta_r|<1$, we obtain
\[
c\frac{mk_\phi}{N_{\mathrm{div}}^\phi}\delta_r^2
\le d_\Phi^2(\mPhi_i,\mPhi_j)
\le C\frac{mk_\phi}{N_{\mathrm{div}}^\phi}\delta_r^2.
\]
The equality with the squared Frobenius distance also gives
\[
\|\mPhi_i-\mPhi_j\|_F^2
\le C\frac{mk_\phi}{N_{\mathrm{div}}^\phi}\delta_r^2.
\]
This proves (ii) and (iii).

\end{proof}

\subsection{Proof of Lemma~\ref{lem:phase-oracle-kl}}
\label{Pf:lem-phase-oracle-kl}

\begin{proof}
Condition on $\mA$, for $\ell\in\{i,j\}$, let $\mC_\ell:=\mPhi_\ell(\sigma_z^2\mI_m+\mS_o)\mPhi_\ell^H$ denote the covariance matrix of the one-look distribution $Q_{\mPhi_\ell} \sim \mathcal{CN}(0,\mC_\ell)$, then by Lemma~\ref{lem:KL-div-Gaussian},
\begin{align}\label{eq:phase-one-look-kl}
D_{\rm KL}(Q_{\mPhi_i}\|Q_{\mPhi_j})
&=
\log\frac{\det\mC_j}{\det\mC_i}
+\Tr(\mC_j^{-1}\mC_i)-m 
= \Tr(\mC_j^{-1}\mC_i)-m
\end{align}
where the log term is zero due to $\mC_i$ and $\mC_j$ are unitarily similar. By the nonnegativity of KL divergence, we have
\[
\begin{aligned}
D_{\rm KL}(Q_{\mPhi_i}\|Q_{\mPhi_j})
&\le
D_{\rm KL}(Q_{\mPhi_i}\|Q_{\mPhi_j})
+D_{\rm KL}(Q_{\mPhi_j}\|Q_{\mPhi_i})\\
&=
\Tr(\mC_j^{-1}\mC_i)
+\Tr(\mC_i^{-1}\mC_j)-2m\\
&=
\left\|
\mC_j^{-1/2}
(\mC_i-\mC_j)
\mC_i^{-1/2}
\right\|_{\rm HS}^2\\
&\le
\frac{
\|\mPhi_i\mS_o\mPhi_i^H
-\mPhi_j\mS_o\mPhi_j^H\|_{\rm HS}^2
}{
\lambda_{\min}^2(\sigma_z^2\mI_m+\mS_o)
}.
\end{aligned}
\]
Here the third equality follows by expanding the squared HS norm and using the cyclicity of trace, while the last inequality uses $\|ABC\|_{\rm HS} \le\|A\|_2\|B\|_{\rm HS}\|C\|_2$ and the unitary similarity of $\mC_i$, $\mC_j$, and $\sigma_z^2\mI_m+\mS_o$.

We next bound the HS norm in the numerator. Let $\va_p^H$ denote the $p$-th row of $\mA$, and write $\mathbf U=\operatorname{diag}(u_1,\ldots,u_m):=\mPhi_j^H\mPhi_i$. By unitary invariance,
\[
\begin{aligned}
&\|\mPhi_i\mS_o\mPhi_i^H
-\mPhi_j\mS_o\mPhi_j^H\|_{\rm HS}^2
= \|\mathbf U\mS_o\mathbf U^H-\mS_o\|_{\rm HS}^2
= \sum_{p\ne q}|u_p-u_q|^2|\va_p^H\mX_o^2\va_q|^2.
\end{aligned}
\]
Since $|u_p-u_q|^2 \le2|1-u_p|^2+2|1-u_q|^2$, we obtain
\[
\begin{aligned}
&\|\mPhi_i\mS_o\mPhi_i^H
-\mPhi_j\mS_o\mPhi_j^H\|_{\rm HS}^2
\le 4\|\mPhi_i-\mPhi_j\|_F^2 \max_{1\le p\le m} \sum_{q\ne p} |\va_p^H\mX_o^2\va_q|^2.
\end{aligned}
\]
Fix $p$ and conditional on $\va_p$, the random variables $\{\va_p^H\mX_o^2\va_q:q\ne p\}$ are independent and satisfy
\[
\va_p^H\mX_o^2\va_q\mid\va_p
\sim
\mathcal{CN}
\left(0,\|\mX_o^2\va_p\|_2^2\right).
\]
Consequently,
\[
\frac{2\sum_{q\ne p} |\va_p^H\mX_o^2\va_q|^2}{\|\mX_o^2\va_p\|_2^2}
\,\bigg|\,\va_p
\sim\chi^2_{2(m-1)}.
\]
Applying the Lemma~\ref{lem:laurent-massart} with $t=m$ gives
\[
\mathbb P\left(
 \sum_{q\ne p} |\va_p^H\mX_o^2\va_q|^2
>C_0m\|\mX_o^2\va_p\|_2^2
\,\middle|\,\va_p
\right)
\le e^{-m}
\]
for an absolute constant $C_0>0$. Moreover,
\[
\|\mX_o^2\va_p\|_2^2
\le x_{\max}^4\|\va_p\|_2^2,
\qquad
2\|\va_p\|_2^2\sim\chi^2_{2n}.
\]
Applying Lemma~\ref{lem:laurent-massart} again with $t=n$ yields
\[
\mathbb P\left(\|\mX_o^2\va_p\|_2^2>C_0x_{\max}^4n\right)
\le e^{-n},
\]
after increasing $C_0$ if necessary.

Combining these two bounds and taking a union bound over $p=1,\ldots,m$, we obtain
\[
\mathbb P\left(\max_{1\le p\le m}\sum_{q\ne p} |\va_p^H\mX_o^2\va_q|^2>C_0^2x_{\max}^4mn\right)
\le me^{-m}+me^{-n}.
\]
Hence, with probability at least $1-me^{-m}-me^{-n}$,
\[
\|\mPhi_i\mS_o\mPhi_i^H
-\mPhi_j\mS_o\mPhi_j^H\|_{\rm HS}^2
\le
Cmn\|\mPhi_i-\mPhi_j\|_F^2,
\]
simultaneously for all diagonal unitary matrices $\mPhi_i,\mPhi_j$, where $C>0$ depends only on $x_{\max}$. 
Finally, since the $L$ looks are conditionally independent, we get
\[
D_{\rm KL}(P_{\mPhi_i}\|P_{\mPhi_j})
= L D_{\rm KL}(Q_{\mPhi_i}\|Q_{\mPhi_j})
\le \frac{CLmn}{\lambda_{\min}^2(\sigma_z^2\mI_m+\mS_o)}\|\mPhi_i-\mPhi_j\|_F^2.
\]
This completes the proof.
\end{proof}

\subsection{Proof of Lemma~\ref{lem:refined-kl-tall-phase}}\label{pf:lem:refined-kl-tall-phase}

\begin{proof}
We condition on $\mA$ and assume $\sigma_z^2>0$. For notation convience, we set $\mathbf H:=\mPhi_j^H\mPhi_i-\mI_m$. The calculation in the proof of Lemma~\ref{lem:phase-oracle-kl} in \eqref{eq:phase-one-look-kl} gives
\[
D_{\rm KL}(Q_{\mPhi_i}\|Q_{\mPhi_j})
=
\Tr\bigl((\sigma_z^2\mI_m+\mS_o)^{-1}(\mI_m+\mathbf H)(\sigma_z^2\mI_m+\mS_o)(\mI_m+\mathbf H)^H\bigr)-m.
\]
Expanding the KL divergence gives
\[
\begin{aligned}
D_{\rm KL}(Q_{\mPhi_i}\|Q_{\mPhi_j})
&=
\Tr(\mathbf H+\mathbf H^H) +\Tr((\sigma_z^2\mI_m+\mS_o)^{-1}\mathbf H(\sigma_z^2\mI_m+\mS_o)\mathbf H^H)\\
&\overset{(a)}{=}
-\Tr((\sigma_z^2\mI_m+\mS_o)^{-1}\mS_o\mathbf H\mathbf H^H)
+\Tr((\sigma_z^2\mI_m+\mS_o)^{-1}\mathbf H\mS_o\mathbf H^H)\\
&\overset{(b)}{\le}
-\Tr((\sigma_z^2\mI_m+\mS_o)^{-1}\mS_o\mathbf H\mathbf H^H)
+\frac{1}{\sigma_z^2}\Tr(\mathbf H\mS_o\mathbf H^H)\\
&\overset{(c)}{=}
\frac{1}{\sigma_z^2}
\Tr((\sigma_z^2\mI_m+\mS_o)^{-1}\mS_o^2\mathbf H\mathbf H^H).
\end{aligned}
\]
where $(a)$ follows from $\mathbf H+\mathbf H^H=-\mathbf H\mathbf H^H$. The inequality in $(b)$ uses $(\sigma_z^2\mI_m+\mS_o)^{-1}\preceq\sigma_z^{-2}\mI_m$ and $\mathbf H\mS_o\mathbf H^H\succeq0$, and $(c)$ follows from $\mathbf H^H\mathbf H=\mathbf H\mathbf H^H$.
Since $\mS_o\succeq0$, we obtain
\begin{align}\label{eq:app-KLbound-phase}
D_{\rm KL}(Q_{\mPhi_i}\|Q_{\mPhi_j})
&\le
\frac{1}{\sigma_z^2}
\Tr((\sigma_z^2\mI_m+\mS_o)^{-1}\mS_o^2\mathbf H\mathbf H^H) \nonumber\\
&\le
\frac{\|\mS_o\|_2}
{\sigma_z^2(\sigma_z^2+\|\mS_o\|_2)}
\Tr(\mS_o\mathbf H\mathbf H^H) \nonumber\\
&\le
\frac{\|\mS_o\|_2
\max_{1\le p\le m}(\mS_o)_{pp}}
{\sigma_z^2(\sigma_z^2+\|\mS_o\|_2)}
\|\mPhi_i-\mPhi_j\|_F^2.
\end{align}
where the last inequality uses the diagonal structure of $\mathbf H$ and $\|\mathbf H\|_F=\|\mPhi_i-\mPhi_j\|_F$.

It remains to control $\|\mS_o\|_2$ and its diagonal entries. Let $\va_p^H$ denote the $p$-th row of $\mA$. By Lemma~\ref{lem:svd-gaussian}, taking $t=\sqrt m$ and using $n\le m$ gives
\[
\mathbb P\left(\|\mA\|_2^2>9m\right) \le 2e^{-m/2}.
\]
Moreover, since $2\|\va_p\|_2^2\sim\chi^2_{2n}$, the Lemma~\ref{lem:laurent-massart} and a union bound yield
\[
\mathbb P\left(\max_{1\le p\le m}\|\va_p\|_2^2>4n\right)
\le me^{-n}.
\]
Thus, on the event
\[
\mathcal E:=
\left\{ \|\mA\|_2^2\le9m,\; \max_{1\le p\le m}\|\va_p\|_2^2\le4n
\right\},
\]
which holds with probability $\mathbb P(\mathcal E) \ge 1 - 2e^{-m/2}-me^{-n}$, we have
\[
\|\mS_o\|_2\le9x_{\max}^2m,
\qquad
\max_{1\le p\le m}(\mS_o)_{pp} \le4x_{\max}^2n.
\]
Substituting these estimates into \eqref{eq:app-KLbound-phase} and using the monotonicity of $t\mapsto t/(\sigma_z^2+t)$ gives
\[
\begin{aligned}
D_{\rm KL}(Q_{\mPhi_i}\|Q_{\mPhi_j})
&\le
\frac{36x_{\max}^4mn}
{\sigma_z^2(\sigma_z^2+9x_{\max}^2m)}
\|\mPhi_i-\mPhi_j\|_F^2
\le
C\frac{mn}
{\sigma_z^2\max\{\sigma_z^2,m\}}
\|\mPhi_i-\mPhi_j\|_F^2,
\end{aligned}
\]
where $C>0$ depends only on $x_{\max}$. Finally, since the $L$ looks are conditionally independent,
\[
\begin{aligned}
D_{\rm KL}(P_{\mPhi_i}\|P_{\mPhi_j})
&=
L D_{\rm KL}(Q_{\mPhi_i}\|Q_{\mPhi_j})
\le
CL\frac{mn}
{\sigma_z^2\max\{\sigma_z^2,m\}}
\|\mPhi_i-\mPhi_j\|_F^2.
\end{aligned}
\]
Since $\mathcal E$ depends only on $\mA$,
the bound holds simultaneously for all
diagonal unitary matrices $\mPhi_i,\mPhi_j$.
\end{proof}
\section{Proof of Theorem~\ref{thm:signal-only} (Signal-Only Upper Bound)}
\label{sec:signal-appendix} 
To start, we note the main technical tools in the proof. 
We use the real Gaussian counterpart of
Lemma~\ref{lem:svd-gaussian}, i.e., the same bounds for independent
$\mathcal N(0,1)$ entries follow from
\cite[Theorem~II.13]{davidson2001local}.
We then note the real Hanson--Wright inequality in
\cite[Theorem~1.1]{rudelson2013}, which is the real-valued counterpart of Lemma~\ref{Thm:ComplexHW}.

We will also use the real-valued counterpart of
Lemma~\ref{lem:quad-yt}. Specifically, let
$\vu_1,\ldots,\vu_L\stackrel{\mathrm{i.i.d.}}{\sim}
\mathcal N(0,\mK)$ with $\mK\succ0$, and let $\mH$ be a real
symmetric matrix. Then, for every $t>0$,
\begin{align}
\Pp\left(
\left|
\frac{1}{L}\sum_{\ell=1}^L \vu_\ell^\top \mH\vu_\ell
-\Tr(\mH\mK)
\right|>t
\right)
\le
2\exp\left[
-cL\min\left\{
\frac{t^2}{\|\mK^{1/2}\mH\mK^{1/2}\|_{\mathrm{HS}}^2},
\frac{t}{\|\mK^{1/2}\mH\mK^{1/2}\|_2}
\right\}
\right].
\label{eq:real-quadratic-form-concentration}
\end{align}
This follows from the real Hanson--Wright inequality by the same
whitening argument as in the proof of Lemma~\ref{lem:quad-yt}.

We first prove exact recovery of $\mD_o$ up to a sign. 

\begin{proposition}[Exact sign recovery]
\label{prop:signal-signs}
For every $B>0$, there is a constant $C_B>0$, depending only on
$B,x_{\min},x_{\max}$, such that, under \eqref{eq:signal-model} and
\(
 m\ge C_B,~n\ge C_B\log(em),~L\ge C_B\frac{V}{mn}\log(em),\)
we have 
\begin{equation}
 \Pp\bigl(\widehat{\mD}\notin\{\mD_o,-\mD_o\}\bigr)
 \le C_Bm^{-B}.
 \label{eq:signal-sign-failure}
\end{equation}
\end{proposition}

\begin{proof}
We write $\va_i^\top$ for the $i$th row of $\mA$ and define
\(
 \mK_o=\sigma_z^2\mI_m+\mA\mX_o^2\mA^\top,~
 \Lambda=\sigma_z^2+m+n,~ \Lambda^2\asymp V.\)
Constants denoted by $c,C$ below may change from line to line, and constants
with a subscript $B$ may also depend on the prescribed probability exponent.

\smallskip\noindent
\textit{Step 1: reduction to positive cut weights.}
In the proof, remove the true signs by setting
$\vq_\ell=\mD_o\vy_\ell$.
Conditional on $\mA$, these vectors are independent $\mathcal N(0,\mK_o)$.
Let $\widehat{\mQ}=\mD_o\widehat{\mW}\mD_o$, and for
$T\subset[m]$ define
$\widehat H(T)=\sum_{i\in T,j\notin T}\widehat Q_{ij}$.
For any $\vv\in\{\pm1\}^m$, we let $T(\vv)=\{i:v_i=-1\}$. By symmetry, we have 
\begin{equation}
 \vone^\top\widehat{\mQ}\vone-\vv^\top\widehat{\mQ}\vv
 =4\widehat H(T(\vv)).
 \label{eq:signal-cut}
\end{equation}
To see this, note that only edges crossing the cut contribute, and each is counted in
both orientations. Since $\vv=\vec{\varepsilon}\circ\vec{\varepsilon}_o$,
strict positivity of every nontrivial cut makes
$\pm\vec{\varepsilon}_o$ the only maximizers of \eqref{eq:signal-sync}.
The corresponding population cut weight is
\begin{equation}
 H_{\mA}(T):=\E[\widehat H(T)\mid\mA]
 =\sum_{\substack{i\in T\\j\notin T}}
  (\va_i^\top\mX_o^2\va_j)(\va_i^\top\va_j).
 \label{eq:signal-population}
\end{equation}
The additive-noise covariance contributes only to the diagonal and hence
does not appear here.

\smallskip\noindent
\textit{Step 2: a uniform population gap.}
By symmetry it suffices to consider $r=|T|\le m/2$.
Write $S=T^c$, $s=m-r$, and $\ell_r=\log(em/r)$.
The real Gaussian singular-value and Hanson--Wright bounds described above,
followed by union bounds, give an event $\mathcal G$ with
\begin{equation}
 \Pp(\mathcal G^c)\le C_Bm^{-B}+Cm e^{-cn},\qquad
 \begin{cases}
  cn\le\|\va_i\|_2^2\le Cn,&i\in[m],\\
  \|\mA_T\|_2^2\le C_B(n+r\ell_r),&
  |T|=r,\ 1\le r\le m/2.
 \end{cases}
 \label{eq:signal-design-event}
\end{equation}
For the spectral-norm bound, fix $r$ and apply
Lemma~\ref{lem:svd-gaussian} to each $\mA_T$ with
$t=\sqrt{K_Br\ell_r}$. There are at most
$\binom mr\le e^{r\ell_r}$ subsets $T$ of size $r$, so a union bound
controls the bound simultaneously over all such $T$. Choosing $K_B$
sufficiently large and then summing over $1\le r\le m/2$ gives a total
failure probability at most $C_Bm^{-B}$.

For a fixed $T$, condition on $\mA_T$ and set
\(
 \mB_T=\mA_T^\top\mA_T,~
 \mQ_T=\tfrac12(\mX_o^2\mB_T+\mB_T\mX_o^2),~
 \tau_T=\Tr(\mQ_T)=\Tr(\mX_o^2\mB_T).
\)
The remaining rows are independent standard Gaussian vectors, and
\begin{equation}
 H_{\mA}(T)=\sum_{j\in S}\va_j^\top\mQ_T\va_j.
 \label{eq:signal-conditional-cut}
\end{equation}
Let $\mathcal F_T$ be the event that the row bounds in
\eqref{eq:signal-design-event} hold for $i\in T$ and that its bound for
$\mA_T$ holds. Then $\mathcal F_T$ depends only on $\mA_T$ and
$\mathcal G\subset\mathcal F_T$. On $\mathcal F_T$,
\begin{equation}
 \tau_T\asymp nr,\qquad
 \|\mQ_T\|_2\le C_B(n+r\ell_r),\qquad
 \|\mQ_T\|_{\mathrm{HS}}^2
 \le x_{\max}^4\Tr(\mB_T)\|\mB_T\|_2
 \le C_Bnr(n+r\ell_r).
 \label{eq:signal-local-norms}
\end{equation}

Conditional on $\mA_T$, the matrix $\mQ_T$ is fixed, while
$\{\va_j\}_{j\in S}$ are independent $\mathcal N(0,\mI_n)$ vectors.
Therefore, by \eqref{eq:signal-conditional-cut},
\[
 \E\left[\frac{1}{s}H_{\mA}(T)\,\middle|\,\mA_T\right]
 =\Tr(\mQ_T)=\tau_T.
\]
Consequently,
\[
 \left\{H_{\mA}(T)\le \frac12s\tau_T\right\}
 \subseteq
 \left\{
 \left|\frac{1}{s}H_{\mA}(T)-\tau_T\right|
 \ge \frac{\tau_T}{2}
 \right\}.
\]
Applying \eqref{eq:real-quadratic-form-concentration} conditionally on
$\mA_T$, with $L=s$, $\mK=\mI_n$, $\mH=\mQ_T$, and
$t=\tau_T/2$, gives
\begin{align}
 \Pp\left(H_{\mA}(T)\le\tfrac12s\tau_T\mid\mA_T\right)
 &\le
 2\exp\left[-c\min\left\{
    \frac{s\tau_T^2}{\|\mQ_T\|_{\mathrm{HS}}^2},
    \frac{s\tau_T}{\|\mQ_T\|_2}
 \right\}\right]\notag\\
 &\le
 2\exp\left[-c_B\frac{snr}{n+r\ell_r}\right].
 \label{eq:signal-pop-tail}
\end{align}
where the last inequality follows from
\eqref{eq:signal-local-norms}.
Since $r\ell_r^2\le Cm$ and $s\ge m/2$,
\begin{equation}
 \frac{sn}{\ell_r(n+r\ell_r)}
 \ge c\min\left\{\frac{m}{\log(em)},\,n\right\}.
 \label{eq:signal-entropy}
\end{equation}
After enlarging the constants in the assumptions on $m,n$, the exponent
in \eqref{eq:signal-pop-tail} is at least $(B+3)r\ell_r$.
The conditioning remains on $\mA_T$. Specifically, if
$\mathcal B_T=\{H_{\mA}(T)\le s\tau_T/2\}$, then
\begin{align*}
 \Pp\left(\mathcal G\cap\bigcup_{\substack{T\subset[m]\\1\le|T|\le m/2}}
                 \mathcal B_T\right)
 &\le\sum_{r=1}^{\lfloor m/2\rfloor}\ \sum_{|T|=r}
       \E\!\left[\mathbf1_{\mathcal F_T}
           \Pp(\mathcal B_T\mid\mA_T)\right]\\
 &\le 2\sum_{r=1}^{\lfloor m/2\rfloor}e^{-(B+2)r\ell_r}
 \le C_Bm^{-B}.
\end{align*}
Consequently, with failure probability at most $C_Bm^{-B}+Cm e^{-cn}$,
\begin{equation}
 H_{\mA}(T)\ge cn|T|(m-|T|)
 \quad\text{for every nonempty proper }T\subset[m].
 \label{eq:signal-pop-gap}
\end{equation}

\smallskip\noindent
\textit{Step 3: preservation of the gaps with finitely many looks.}
The real version of Lemma~\ref{lem:svd-gaussian} and the real, unweighted row-sum
calculation in the proof of Lemma~\ref{lem:phase-oracle-kl} give a further design event on which
\begin{equation}
 E_{\max}:=\|\mA\|_2^2\le C(m+n),\qquad
 \max_i\sum_{j\ne i}G_{ij}^2\le Cmn,
 \label{eq:signal-gram-event}
\end{equation}
with failure probability at most $Cm e^{-cn}+Cm e^{-cm}$.
Here $E_{\max}$ is the real analogue of \eqref{eq:def-Emax-norm2A}.
For the row-sum bound, the identity used in the proof of Lemma~\ref{lem:phase-oracle-kl} becomes
$\sum_{j\ne i}(\va_i^\top\va_j)^2\mid\va_i
 \overset d=\|\va_i\|_2^2\chi^2_{m-1}$.

For $r=|T|\le m/2$, let $\mP_T$ be the coordinate projector onto $T$ and set
$\mR_T=\mP_T\mG\mP_S+\mP_S\mG\mP_T$. Then
\begin{equation}
 \widehat H(T)=\frac1{2L}\sum_{\ell=1}^L
       \vq_\ell^\top\mR_T\vq_\ell,\qquad
 H_{\mA}(T)=\tfrac12\Tr(\mR_T\mK_o).
 \label{eq:signal-empirical-form}
\end{equation}
On \eqref{eq:signal-gram-event}, with
$\mM_T=\mK_o^{1/2}\mR_T\mK_o^{1/2}$, the needed norm bounds are
\begin{align}
 \|\mR_T\|_{\mathrm{HS}}^2
 &=2\sum_{i\in T,j\in S}G_{ij}^2\le Cnmr\le Cnrs,
 &\|\mR_T\|_2&=\|\mG_{T,S}\|_2\le E_{\max},\notag\\
 \|\mK_o\|_2&\le\sigma_z^2+x_{\max}^2E_{\max}\le C\Lambda,
 &\|\mM_T\|_{\mathrm{HS}}^2&\le C\Lambda^2nrs,\quad
   \|\mM_T\|_2\le C\Lambda^2.
 \label{eq:signal-empirical-norms}
\end{align}
Conditional on $\mA$, the vectors $\vq_1,\ldots,\vq_L$ are independent
$\mathcal N(0,\mK_o)$ vectors. By \eqref{eq:signal-empirical-form},
\[
 \widehat H(T)
 =\frac{1}{2L}\sum_{\ell=1}^L \vq_\ell^\top\mR_T\vq_\ell,
 \qquad
 H_{\mA}(T)=\frac12\Tr(\mR_T\mK_o).
\]
Thus, applying \eqref{eq:real-quadratic-form-concentration} conditionally
on $\mA$, with $\mK=\mK_o$, $\mH=\mR_T/2$, and
$t=H_{\mA}(T)/2$, gives
\begin{align}
 \Pp\left(\left|\widehat H(T)-H_{\mA}(T)\right|
           \ge\tfrac12H_{\mA}(T)\mid\mA\right)
 &\le 2\exp\left[-c\min\left\{
       \frac{LH_{\mA}(T)^2}{\|\mM_T\|_{\mathrm{HS}}^2},
       \frac{LH_{\mA}(T)}{\|\mM_T\|_2}\right\}\right]\notag\\
 &\le 2\exp\left[-c\frac{Lnrs}{\Lambda^2}\right].
 \label{eq:signal-empirical-tail}
\end{align}
Here the last inequality follows from the population gap
\eqref{eq:signal-pop-gap} and the norm bounds
\eqref{eq:signal-empirical-norms}.

Since $\Lambda^2\asymp V$ and $s\ge m/2$, the assumed lower bound on
$L$ makes the exponent in \eqref{eq:signal-empirical-tail} at least
$(B+3)r\log(em)$. The same cut-counting argument as above then shows
that every nontrivial empirical cut is strictly positive, except on an
event of conditional probability at most $C_Bm^{-B}$.

Combining the design and observation events with \eqref{eq:signal-cut} yields
a total  failure probability at most
$C_Bm^{-B}+Cm e^{-cn}+Cm e^{-cm}$. Under $n\ge C_B\log(em)$ and $m\ge C_B$,
the exponential terms are absorbed into $C_Bm^{-B}$ after increasing the
constants. This proves \eqref{eq:signal-sign-failure}.
\end{proof}

\paragraph{Reduction to the signal-only estimation problem.}
Let $\rho$ denote the minimum on the right-hand side of
\eqref{eq:signal-bound}. Under the condition on $mL$,
\cite[Theorem~2.14]{xing2025minimax} provides a signal estimator
$\mathcal T^0_{\mA}$
whose worst-case risk in the model without distortion is at most $C\rho$. Fix such an estimator and define
\begin{equation}
 \mathcal T_{\mA}(\mZ)
 =\frac12\left[
     \Pi_{\mathrm{box}}\mathcal T^0_{\mA}(\mZ)
    +\Pi_{\mathrm{box}}\mathcal T^0_{\mA}(-\mZ)\right],
 \label{eq:signal-symmetrization}
\end{equation}
where $\Pi_{\mathrm{box}}$ clips each coordinate to $[x_{\min},x_{\max}]$.
In the model without distortion,
$(\mA,\mZ)\overset d=(\mA,-\mZ)$ for every fixed $\vx_o$.
Clipping and Jensen's inequality therefore show that
\eqref{eq:signal-symmetrization} does not increase any fixed-parameter risk.
It is measurable, takes values in the box, and is exactly invariant under
$\mZ\mapsto-\mZ$.

Define the oracle estimate
$\widehat{\vx}_{\mathrm{or}}=\mathcal T_{\mA}(\mD_o\mY)$ and the event
$\mathcal E=\{\widehat{\mD}\in\{\mD_o,-\mD_o\}\}$.
On $\mathcal E$, the invariance just established gives
$\widehat{\vx}=\widehat{\vx}_{\mathrm{or}}$ pointwise.
Both estimates lie in $[x_{\min},x_{\max}]^n$, so
\begin{align}
 \E d_x^2(\widehat{\vx},\vx_o)
 &\le \E d_x^2(\widehat{\vx}_{\mathrm{or}},\vx_o)
      +(x_{\max}-x_{\min})^2\Pp(\mathcal E^c)\notag\\
 &\le C\rho+C_Bm^{-B}.
 \label{eq:signal-oracle-risk}
\end{align}
The first term has the required oracle law because
$\mD_o\vy_\ell=\mA\mX_o\vw_\ell+\mD_o\vz_\ell$ and
$\mD_o\vz_\ell\sim \mathcal N(0,\sigma_z^2\mI_m)$.
The second term follows from Proposition~\ref{prop:signal-signs}.
This is an unconditional comparison, so the same data may be used in both
stages. 
Finally, if $V\ge nL\log m$, the second term in the minimum is bounded by the
first, proving the simplified form.

\newpage

\bibliographystyle{IEEEtran}
\bibliography{refs}

\end{document}